%% file: Bolder_Fusion_main_I.tex
\documentclass[11pt]{article}
\pdfoutput=1
\usepackage[margin=1in,marginparwidth=.95in]{geometry}

\usepackage[utf8x]{inputenc}
\usepackage{amsmath,amssymb,amsthm}
\usepackage[pdftex]{graphicx}   
\usepackage{color}
\usepackage{ifthen}
\usepackage{mathtools}
\usepackage{todonotes}
\usepackage{verbatim}

\usepackage{yfonts}

\usepackage{tikz}
\usetikzlibrary{arrows,matrix}
\tikzset{arrow/.style={semithick,>=stealth',shorten >=1pt,shorten <=1pt}}

\usepackage[all]{xy}

\usepackage[pdftex,     
            plainpages=false,   
            breaklinks=true,    
            pdftitle={Transfer for free loop spaces of fusion systems and covering maps},
            pdfauthor={Sune Precht Reeh, Tomer M. Schlank, Nathaniel Stapleton}
           ]{hyperref}
\usepackage{thumbpdf}

\usepackage[abbrev,msc-links]{amsrefs}

\numberwithin{equation}{section}
\numberwithin{figure}{section}
\numberwithin{table}{section}
\allowdisplaybreaks

\theoremstyle{plain}
\newtheorem{theorem}{Theorem}[section]
\newtheorem{prop}[theorem]{Proposition}
\newtheorem{lemma}[theorem]{Lemma}
\newtheorem{cor}[theorem]{Corollary}

\theoremstyle{definition}
\newtheorem{definition}[theorem]{Definition}
\newtheorem{convention}[theorem]{Convention}

\newtheoremstyle{remark}%
{5pt}
{11pt}
{}
{}
{\itshape}
{.}
{ }
{}

\theoremstyle{remark}
\newtheorem{remark}[theorem]{Remark}

\newtheorem{example}[theorem]{Example}

\newcommand{\ph}{\varphi}
\newcommand{\e}{\varepsilon}
\renewcommand{\hat}{\widehat}
\renewcommand{\tilde}{\widetilde}

\DeclareMathOperator{\Hom}{Hom}
\DeclareMathOperator{\Map}{Map}

\DeclareMathOperator{\Ho}{Ho}
\DeclareMathOperator*{\colim}{colim}
\newcommand{\tprod}{{\textstyle\prod}}

\newcommand{\CC}{\mathbb{C}}

\newcommand{\ZZ}{\mathbb{Z}}

\newcommand{\x}{\times}

\newcommand{\into}{\hookrightarrow}

\newcommand{\xto}{\xrightarrow}

\newcommand{\tto}{\duorightarrow}
\newcommand{\ffrom}{\duoleftarrow}
\newcommand{\xtto}{\xduorightarrow}
\newcommand{\xffrom}{\xduoleftarrow}
\newcommand{\cmp}{\mathbin{\odot}}

\newcommand{\ntuples}[2]{{#2}^{(#1)}}
\newcommand{\cntuples}[2]{{#2}^{[#1]}}

\newcommand{\freeO}[1]{\ensuremath \Loops^{#1}}
\newcommand{\tup}{\underline}

\newcommand{\change}{\zeta}
\newcommand{\rep}{\PackageError{main}{\\rep not currently implemented}{Use \\reptup instead}\textswab}
\newcommand{\reptup}[1]{\tup{\textswab{#1}}}
\newcommand{\cover}[1]{\widetilde{#1}}
\newcommand{\PB}[1]{PB^{#1}}
\newcommand{\PBlift}[1]{PB_{\mathrm{lift}}^{#1}}
\newcommand{\PBbroke}[1]{PB_{\mathrm{broken}}^{#1}}
\DeclareMathOperator{\wind}{wind}

\newcommand{\Z}{\mathbb{Z}}

\newcommand{\R}{\mathbb{R}}

\newcommand{\Loops}{\mathcal{L}}
\newcommand{\twistO}[1]{L^{\dagger}_{#1}}
\newcommand{\untwistO}[1]{L^{\dagger}_{#1}|_{\mathrm{lift}}}

\newcommand{\mmod}[1][]{%
    \mathrel{\!\ooalign{$#1/$\cr%
        $\mkern6mu #1/$}\!}%
}

\DeclareMathOperator{\Id}{Id}

\DeclareMathOperator{\Span}{Span}

\DeclareMathOperator{\id}{id}
\DeclareMathOperator{\incl}{incl}

\DeclareMathOperator{\tr}{tr}
\DeclareMathOperator{\ev}{ev}
\DeclareMathOperator{\pev}{\partial ev}
\DeclareMathOperator{\cC}{\mathcal{C}}

\DeclarePairedDelimiter{\gen}{\langle}{\rangle}
\DeclarePairedDelimiter{\abs}{\lvert}{\rvert}

\newcommand{\AG}{\mathbb{AG}}

\newcommand{\Sp}{\mathord{\mathrm{Sp}}}

\newcommand{\Top}{\mathord{\mathrm{Top}}}
\newcommand{\Grp}{\mathord{\mathrm{Grp}}}
\DeclareMathOperator{\SubGrp}{SubGrp}
\newcommand{\Cov}{\mathord{\mathrm{Cov}}}

\newcommand{\Sinfp}{\Sigma^{\infty}_{+}}
\newcommand{\Sinfpp}{\hat\Sigma^{\infty}_{+}}

\newcommand{\lc}[1]{\prescript{#1}{}\!}

\newcommand\duorightarrow{%
  \mathrel{\ooalign{$\rightarrow$\cr%
  $\mkern4mu\rightarrow$}}%
}
\newcommand\xduorightarrow[2][]{%
  \mathrel{\ooalign{$\xrightarrow[#1\mkern4mu]{#2\mkern4mu}$\cr%
  \hidewidth$\rightarrow\mkern4mu$}}%
}
\newcommand\duoleftarrow{%
  \mathrel{\ooalign{$\leftarrow$\cr%
  $\mkern4mu\leftarrow$}}%
}
\newcommand\xduoleftarrow[2][]{%
  \mathrel{\ooalign{$\xleftarrow[\mkern4mu #1]{\mkern4mu #2}$\cr%
  $\mkern4mu\leftarrow$\hidewidth}}%
}

\title{Evaluation maps and transfers for free loop spaces I}
\author{Sune Precht Reeh \and Tomer M. Schlank \and Nathaniel Stapleton}
\date{}

\begin{document}

\maketitle

\begin{abstract}
We construct and study a functorial extension of the evaluation map $S^1 \times \Loops X \to X$ to transfers along finite covers. For finite covers of classifying spaces of finite groups, we provide algebraic formulas for this extension in terms of bisets. In the sequel \cite{RSS_Bold2}, we show that this induces a natural evaluation map on the full subcategory of the homotopy category of spectra consisting of $p$-completed classifying spectra of finite groups.
\end{abstract}

\tableofcontents

\section{Introduction}
The free loop space functor 
\[
X \mapsto \Loops(X) 
\]
on the category of topological spaces comes equipped with a natural evaluation map
\[
S^1 \times \Loops(X) \to X.
\]
This natural map arises as the counit of the adjunction between the functors $S^1 \times -$ and $\Loops(-)$ on the category of topological spaces. There are good reasons stemming from chromatic homotopy theory to wonder if this evaluation map can be naturally extended to (certain maps between) suspension spectra. That is: Given a map of spectra $f \colon \Sinfp X \to \Sinfp Y$, is there a map of spectra $L^{\dagger}(f) \colon \Sinfp(S^1 \times \Loops X) \to \Sinfp(S^1 \times \Loops Y)$ that is compatible with $f$ through the suspension spectrum of the evaluation map? 

The maps between suspension spectra that are most accessible are maps induced by maps of spaces and transfer maps along finite covers. Thus it is natural to consider the extension of the category of topological spaces obtained by taking the category of spans of spaces, $\Cov$, in which the backward map is required to be a finite cover. In this case, category theory does not formally provide natural evaluation maps: the category $\Cov$ does not admit many limits or colimits. However, in this paper we show that the evaluation map does admit a natural extension to $\Cov$. 

The extension of the evaluation map natural transformation to objects and forward maps is predetermined. The difficulty lies in constructing an appropriate extension of $S^1 \times \Loops(-)$ to finite covers.  Given a span $f \colon X \ffrom Z \to Y$ in $\Cov$, we construct a span $L^{\dagger}(f) \colon S^1 \times \Loops X \ffrom \mathord{?} \to S^1 \times \Loops Y$ that is not necessarily given by a product of spans. In fact, if $f$ is a span of the form $X \ffrom Y \xto\id Y$ (in which the forward map is the identity), then $L^{\dagger}(f)$ is a span in which the forward map is not necessarily the identity. We provide group-theoretic formulas for $L^{\dagger}(-)$ when restricted to the homotopy category of the full subcategory of $\Cov$ on coproducts of classifying spaces of finite groups. Further, making use of naturality, in \cite{RSS_Bold2} we will extend our construction to the Burnside category of fusion systems. 

The main results of the paper can be stated as follows:

\begin{theorem}\label{thmIntroMain}
We construct a family of endofunctors $\twistO n \colon \Cov\to \Cov$ for $n\geq 0$ with the following list of properties:
\begin{enumerate}
\renewcommand{\theenumi}{$(\roman{enumi})$}\renewcommand{\labelenumi}{\theenumi}
\item\label{itemIntroLnZero} $\twistO 0$ is the identity functor on $\Cov$.
\item\label{itemIntroLnObjects} On objects, $\twistO n$ takes a space $X$ to 
\[\twistO n(X) = (S^1)^n\x \freeO n (X).\]
\item\label{itemIntroEquivariant}The symmetric group $\Sigma_n$ acts on $(S^1)^n\x \freeO n(-)$ diagonally by permuting the coordinates of both $(S^1)^n$ and $\freeO n(-)$. For every $\sigma\in\Sigma_n$ the diagonal action of $\sigma$ on $(S^1)^n\x \freeO n(-)$ induces a natural isomorphism $\sigma\colon \twistO n\underset\cong\Rightarrow \twistO n$.
\item\label{itemIntroLnForwardMaps} On forward maps, i.e. spans $X\xffrom\id X\xto f Y$, the functor $\twistO n$ coincides with the torus times the usual free loop space functor in $\Top$ so that $\twistO n(\id,f)=(S^1)^n\x \freeO n(f)$.
\item\label{itemIntroEvalSquare} For each $n \geq 0$, the evaluation maps $\ev_X\colon (S^1)^n\x \freeO n(X)\to X$ form a natural transformation $\ev\colon \twistO n \Rightarrow \Id_{\Cov}$.
\item\label{itemIntroLnPartialEvaluation} For each $n \geq 0$, the partial evaluation maps $\pev_X\colon S^1\x \freeO {n+1}(X) \to \freeO n (X)$ given by 
\[
\pev(t,f)= (\tup s\mapsto f(\tup s,t)), \qquad \text{for } t\in S^1, \tup s\in (S^1)^{n},
\] form a natural transformation $(S^1)^n\x \pev\colon \twistO {n+1} \Rightarrow \twistO n$.
\item\label{itemIntroIterateLn} For all $n,m\geq 0$, and any space $X$, the space $(S^1)^{n+m}\x \freeO {n+m} X$ embeds into $(S^1)^m\x \freeO m((S^1)^n\x \freeO n X)$ as the collection of those components corresponding to $m$-fold loops in $(S^1)^n\x \freeO n X$ that are constant in the $(S^1)^n$-coordinate.

These embeddings form a natural transformation $\twistO {n+m}(-)\Rightarrow \twistO m(\twistO n(-))$.
\end{enumerate}
In Section \ref{secGroups}, we give explicit group-theoretic formulas for $\twistO n$ on the homotopy category of the full subcategory of $\Cov$ consisting of classifying spaces of finite groups. In \cite{RSS_Bold2}, we use the preceding results for $p$-groups to extend both $\twistO n$ and the theorem to the Burnside category of saturated fusion systems.
In Sections \ref{secGroups} and \cite{RSS_Bold2}, both the category $\Cov$ and the circle $S^1$ are replaced by suitable algebraic models, but otherwise the theorem still holds as stated.
\end{theorem}

This theorem resulted from our (successful) attempt to apply the various flavors of character theory that arise in chromatic homotopy theory to fusion systems. Generally speaking, character maps are a composite of two maps: an evaluation map and a ``change of coefficients." To extend character theory from finite groups to fusion systems, one must contend with the characteristic idempotent, which is built from both transfer maps and group homomorphisms. Thus, from this perspective, it is reasonable to attempt to extend the evaluation map to be natural with respect to transfer maps. Character theory does not demand this level of generality because it applies a finite height cohomology theory (such as $K$-theory or Morava $E$-theory) to the evaluation map -- and this is a destructive procedure. Yet we were surprised to discover that there is a global solution to the problem of extending the evaluation map to transfers and that is what is presented here. 

\subsection{Evaluation maps and finite covers}\label{secIntroCov}

Let $\Top$ be the category of topological spaces that are homotopy equivalent to a CW complex and let $\Cov$ be the category of spans in $\Top$, $X \ffrom E \to Y$, in which the map $X \ffrom E$ is a finite cover. There is a canonical faithful functor $\Top \hookrightarrow \Cov$.

For a space $X$, there is an evaluation map $(S^1)^n \times \Loops^n(X) \to X$ that is natural for $X$ in $\Top$. The first goal is to extend this natural transformation to $\Cov$. That is, we would like a functor $\twistO n \colon \Cov \to \Cov$ and a natural transformation $\twistO n \Rightarrow \id_{\Cov}$ compatible with the inclusion of $\Top$ into $\Cov$ and the evaluation map on $\Top$.
Since $\Top$ is naturally a subcategory with the same objects, the value of $\twistO n$ on objects and forward maps is already determined: on these objects and maps $\twistO n(-) = (S^1)^n \times \Loops^n(-)$.

Since each span decomposes as
\[
X \xffrom{\pi} E \xto{f} Y = (E \xffrom{\id_E} E \xto{f} Y) \circ (X \xffrom{\pi} E \xto{\id_E} E),
\]
it suffices to describe the construction of $\twistO n$ on a span of the form $X \xffrom{\pi} E \xto{\id_E} E$ (i.e. on a finite cover). We will denote this span by $X \xffrom{\pi} E$. Naively, one might guess that the value of $\twistO n$ on $X \xffrom{\pi} E$ should be the finite cover 
\[
(S^1)^n \times \Loops^n(X) \xffrom{(S^1)^n \times \Loops^n(\pi)} (S^1)^n \times \Loops^n(E). 
\]
Although this does produce a functor $(S^1)^n \times \freeO n(-)\colon \Cov\to \Cov$, it does not interact properly with respect to the evaluation map: the diagram
\[
\xymatrix{(S^1)^n \times \Loops^n(E) \ar@{->>}[d]_{(S^1)^n\x \freeO n(\pi)} \ar[r] & E \ar@{->>}[d]^{\pi} \\ (S^1)^n \times \Loops^n(X) \ar[r] & X}
\]
does not necessarily commute in $\Cov$. In order for the square to commute in $\Cov$, $(S^1)^n \times \Loops^n(E)$ must be the pullback, but in most cases it is not. Although this naive attempt does not provide us with a solution, it does indicate that the pullback must play a critical role in the definition of $\twistO n$. 

Denote the pullback of the finite cover $X \xffrom{\pi} E$ along the evaluation map by $\PB{n}(\pi)$. We see from the example above that $\twistO n(X \xffrom{\pi} E)$ must have the form
\[
(S^1)^n \times \Loops^n(X) \ffrom \PB{n}(\pi) \to (S^1)^n \times \Loops^n(E).
\]
Thus the goal is to construct a map of spaces $\PB{n}(\pi) \to (S^1)^n \times \Loops^n(E)$ and show that the resulting construction interacts properly with evaluation.

Assume $n=1$ for simplicity. The space $\PB{}(\pi) = \PB{1}(\pi)$ then consists of triples of the form $(s, l, e)$, where $s$ is a point on the circle $S^1$, $l \colon S^1 \to X$ is a loop in $X$, and $e \in E$ is a point in the fiber over $l(s)$. View $S^1$ as $\R/\Z$ and let $[k] \colon S^1 \to S^1$ be multiplication by $k$. Given a loop $l \in \Loops X$, let $k$ be the smallest positive natural number such that $l[k]$ lifts along $\pi$ to a map $\widetilde{l[k]} \colon S^1 \to E$ going through $e\in E$ in the fiber over $l[k](\tfrac sk)$. We define
\[
\wind(\pi) \colon \PB{}(\pi) \to S^1 \times \Loops(E)
\]
by sending a triple $(s,l,e)$ to $(s,\widetilde{l[k]}(-+t_s))$, where 
\[
t_s = \frac{s-ks}{k}.
\]
Note that $[k](-+t_s)\colon S^1\to S^1$ wraps the circle around itself $k$-times and at the same time sends $s\in S^1$ to itself. Hence $l[k](-+t_s)$ winds the loop $l$ around itself $k$ times while still mapping $s\in S^1$ to the point $l(s)$.
By construction, the loop $\widetilde{l[k]}(-+t_s)\colon S^1\to E$ then has the following two properties: if we evaluate the loop at $s\in S^1$ we get the point $e\in E$ in the fiber over $l(s)$, and if we postcompose the loop with $\pi$ we recover the winded loop $l[k](-+t_s)$. Of course, the loop $\widetilde{l[k]}(-+t_s)\in \freeO{}(E)$ is uniquely determined by these two properties.

This construction can be generalized to produce a map
\[
\wind^n(\pi) \colon \PB{n}(\pi) \to (S^1)^n \times \Loops^n(E).
\]
We show that $\wind^n(\pi)$ is continuous and we define $\twistO n(X \xffrom{\pi} E \xrightarrow{f} Y)$ to be the span
\[
(S^1)^n \times \Loops^n(X) \ffrom \PB{n}(\pi) \xto{((S^1)^n \times \Loops^n(f)) \circ \wind^n(\pi)} (S^1)^n \times \Loops^n(Y).
\]

We prove that $\twistO n \colon \Cov \to \Cov$ is a functor and that there is a natural transformation $\twistO n \Rightarrow \Id_{\Cov}$ that extends the evaluation map natural transformation on $\Top$. Further, we study the action of $\Sigma_n$ on $\twistO n$, the partial evaluation maps, and other naturally occurring relationships between the different functors as $n$ varies. All of these properties are collected in Theorem \ref{thmIntroMain}, which we prove for $\Cov$ as Theorem \ref{thmCovMain}.

\subsection{Evaluation maps and the global Burnside category} \label{intro2}

There is an algebraically controlled subcategory of the homotopy category of $\Cov$. Let $G$ and $H$ be finite groups and let $\AG(G,H)$ be the Grothendieck group of isomorphism classes of finite $H$-free $(G,H)$-bisets under disjoint union. Further, let $\AG_+(G,H)$ denote the commutative monoid of isomorphism classes of finite $H$-free $(G,H)$-bisets under disjoint union. There is a fully faithful functor from $\AG_{+}$ to the homotopy category of $\Cov$ that on objects takes $G$ to the classifying space $BG$. The class of objects of the categories $\AG$ and $\AG_+$ can be extended to formal unions of finite groups by linearity (i.e. letting the sets of morphisms consist of matrices). 

The restrictions of $\Loops^n$ and $\twistO n$ to $\AG_+$ both extend to $\AG$ by linearity. We describe these functors when $n=1$. Given a finite group $G$, we set 
\[
\Loops G = \coprod_{[g]} C_G(g), 
\]
where the formal coproduct is over conjugacy classes of elements in $G$ and $C_G(g)$ denotes the centralizer.

Let $X \in \AG_+(G,H)$ so that $X$ is represented by a finite $H$-free $(G,H)$-biset. The naive functor $\Loops(-)\colon \Cov \to \Cov$ gives rise to a matrix $\Loops X \in \AG_+(\Loops G, \Loops H)$. The entry in this matrix corresponding to the pair of conjugacy classes $[g] \subseteq G$ and $[h] \subseteq H$ is represented by the $(C_G(g),C_H(h))$-biset
\[
(\Loops X)_{[g],[h]} = \{x \in X \mid gx = xh\}.
\]
As $G$ is finite, let $\ell = |G|$. The evaluation map in this setting $\ev \colon \Z \times \Loops G \to G$ is equivalent to a map
\[
\ev \colon \coprod_{[g]} \Z/\ell \times C_G(g) \to G.
\]
The component corresponding to the conjugacy class $[g]$ is the group homomorphism $\ev_g \colon \Z/\ell \times C_G(g) \to G$ sending $(t,z)$ to $g^tz$.

The free abelian group $\AG(G,H)$ has a canonical basis consisting of the isomorphism classes of transitive $H$-free $(G,H)$-bisets. These are determined by a pair $[K, \phi]$, where $K \subseteq G$ is a subgroup taken up to conjugacy and $\phi \colon K \to H$ is a homomorphism that is also taken up to conjugacy.

Now we describe $\twistO 1$ in this setting. Let $X|_{\Z/\ell \times C_G(g)} \in \AG_{+}(\Z/\ell \times C_G(g),H)$ be the restriction of $X$ along the evaluation map $\ev_g \colon \Z/\ell \times C_G(g) \to G$. In terms of the canonical basis, this biset has a decomposition of the form
\[
X|_{\Z/\ell \times C_G(g)} = \sum_{\substack{R\leq C_G(g)\\\ph\colon R\to H}} c_{R,\ph} \cdot [\ev_g^{-1}(R), \ph\circ \ev_g],
\]
where each $c_{R,\ph}$ is a non-negative integer. Furthermore, the decomposition informs us that the pullback $\PB{}(\pi)$ of Section \ref{secIntroCov} is equivalent to $c_{R,\ph}$ copies of $B(\ev_g^{-1}(R))$, for each pair $(R,\ph)$.
Similarly, the map $\wind(\pi)$ of Section \ref{secIntroCov}, when restricted to $B(\ev_g^{-1}(R))$, corresponds to the map
\[
\wind(R,g) \colon \ev_g^{-1}(R) \xto{\cong} \Z/\ell \x R
\]
given by
\[
\wind(R,g)(t,z) = (t,g^{t-kt}z),
\]
where $k$ is the smallest positive natural number such that $g^k$ lies in the subgroup $R$.

We restrict attention to a basis element of the form $[\ev_g^{-1}(R), \ph\circ \ev_g]$.
To understand the group theoretic description of $\twistO{}(X)$, we would like to factor 
\[
[\ev_g^{-1}(R), \ph\circ \ev_g] \colon \Z/\ell \times C_G(g) \xrightarrow{} H
\]
through $\Z/\ell \times \Loops (H)$. The factorization is given by
\[
\Z/\ell \times C_G(g) \xto{\tr} \ev_g^{-1}(R) \xto[\cong]{\wind(R,g)} \Z/\ell \x R \xto{\Z/\ell \x \ph} \Z/\ell \x C_H(\ph(g^k)) \subseteq \Z/\ell \x \Loops H,
\]
where $g^k$ is the smallest power of $g$ lying in $R\leq C_G(g)$.

By additivity, these constructions extend to the category $\AG$. The resulting functors $\twistO n \colon \AG \to \AG$ satisfy the corresponding version of Theorem \ref{thmIntroMain}, which we state and prove as Theorem \ref{thmGroupMain}. Curiously, we discovered this group theoretic description of $\twistO n$ prior to understanding its relation to $\Cov$.

Let $\Sinfp BG$ be the classifying spectrum of the finite group $G$. The Segal conjecture, proved by Carlsson, asserts that 
\[
[\Sinfp BG, \Sinfp BH] \cong \AG(G,H)^{\wedge}_{I},
\]
where $I \subset A(G)$ is the augmentation ideal of the Burnside ring. Thus the full subcategory of the homotopy category of spectra consisting of (finite coproducts of) suspension spectra of finite groups is purely algebraic. 

However, we prove that $\twistO n$ cannot be extended from the Burnside category $\AG$ to the full subcategory of the category of spectra consisting of classifying spectra of finite groups. In fact, we show that no functorial extension of the evaluation map natural transformation on $\Top$ to $\Cov$ can be extended further to classifying spectra of finite groups (see Example \ref{exampleDoesNotPassToSpectra} and Remark \ref{remarkDoesNotPassToSpectra}). Any functorial extension of the evaluation map natural transformation fails to be continuous with respect to the $I$-adic topology. However, in \cite{RSS_Bold2} we show that this obstruction disappears after $p$-completion. This allows us to extend $\twistO n$ to the Burnside category of fusion systems, which contains the homotopy category of $p$-completed classifying spectra of finite groups as a full subcategory.

\subsection*{Acknowledgements} It is a pleasure to thank Mike Hopkins and Haynes Miller for their suggestions. Their helpful comments led to a complete revision of these papers. Further, we thank Cary Malkiewich for generously sharing his thoughts regarding transfers and free loop spaces after his work with Lind in \cite{LindMalkiewich}. Section 2 is based on ideas coming from discussions with him. Finally, we thank Matthew Gelvin and Erg\"un Yal\c{c}{\i}n for their comments.

While working on this paper the first author was funded by the Independent Research Fund Denmark (DFF–4002-00224) and later on by BGSMath and the María de Maeztu Programme (MDM–2014-0445). The second and third author were jointly supported by the US-Israel Binational Science Foundation under grant 2018389. The second author was partially supported by the Alon Fellowship and ISF 1588/18. The third author was partially supported by NSF grant DMS-1906236 and the SFB 1085 \emph{Higher Invariants} at the University of Regensburg. All three authors thank the HIM and the MPIM for their hospitality.

\input{"covering_maps_I.tex"}

\input{"groups_and_bisets_I.tex"}

\appendix

\input{"covering_appendix_I.tex"}

\makeatletter
\def\eprint#1{\@eprint#1 }
\def\@eprint #1:#2 {%
    \ifthenelse{\equal{#1}{arXiv}}%
        {\href{http://front.math.ucdavis.edu/#2}{arXiv:#2}}%
        {\href{#1:#2}{#1:#2}}%
}
\makeatother

\end{document}

%% file: covering_maps_I.tex
\section{Covering maps and free loop spaces}\label{secCoveringMaps}

Consider the evaluation map
\[
(S^1)^n \times \freeO n(-) \xto{\ev} \Id
\]  
as a natural transformation of endofunctors on the category of spaces. The category of spaces embeds in the category of spans of spaces in which the backward maps are finite covers. In this section we construct an extension of the evaluation map to the category of spans. After producing this extension, we establish various identities satisfied by the extension 

Let $\Top$ be the category of topological spaces homotopy equivalent to a CW complex. This category is enriched in topological spaces by putting the compact open topology on the set of continuous maps between two spaces. We will make use of a few properties of this category. First, a finite cover of a space homotopy equivalent to a CW complex is also homotopy equivalent to a CW complex. This follows from the fact that a finite cover of a CW complex naturally has the structure of a CW complex. Secondly, the free loop space of a space homotopy equivalent to a CW complex is homotopy equivalent to a CW complex by \cite{Milnor}. Finally, the free loop space of a finite sheeted covering map in $\Top$ is a finite sheeted covering map between the mapping spaces. We prove this result in Proposition \ref{propMappingSpaceCoveringMap} of Appendix \ref{appendixCovering}.

It is possible to carry out the constructions and prove the results in this section using other categories of spaces. For instance, one may use the category of topological spaces or the category of compactly generated weak Hausdorff spaces. In Appendix \ref{appendixCovering}, we describe the relationship between the free loop space functor and finite sheeted covering maps in these categories. We have chosen to work with the category $\Top$ because it contains all classifying spaces of finite groups and because it allows us to give a simple proof of Proposition \ref{propKContinuousWithPathsProof}.

The category $\Span(\Top)$ of spans of spaces has objects spaces in $\Top$ and morphisms given by isomorphism classes of spans. As usual in these situations, composition is given by taking the (isomorphism class of the) pullback. We define the subcategory $\Cov$ of $\Span(\Top)$ to be the category with objects spaces and morphisms spans for which the backward map is a finite cover and the forward map is a map of spaces. 

We will always take ``finite cover'' to mean a continuous map $p\colon E\to B$ that is locally trivializable with $p^{-1}(U) = V_1\sqcup \dotsb \sqcup V_n$, where each $V_i$ maps homeomorphically to the neighborhood $U$. In particular, we do not require all fibers to have the same size (they necessarily have constant size over each component of $B$), and we allow empty fibers as well.

Since the pullback of a finite cover is a finite cover, composition makes sense in $\Cov$. Note that there is a faithful embedding of $\Top$ into $\Cov$ given by sending a map $f \colon X \to Y$ to the span $X \xffrom\id X \xto f Y$.

There is also a homotopy category associated to $\Cov$. A map in the homotopy category $\Ho(\Cov)$ between spaces $X$ and $Y$ is a finite cover of $X$ and a homotopy class of maps from the finite cover to $Y$. Note that the homotopy category of spaces faithfully embeds in $\Ho(\Cov)$.

Given a finite cover $Y \tto X$ and a torus $(S^1)^n$, the induced map $(S^1)^n \times Y \tto (S^1)^n \times X$ is also a finite cover. It follows that $(S^1)^n \times (-)$ is an endofunctor of $\Cov$. In Appendix \ref{appendixCovering}, we prove that $\freeO n(Y) \tto \freeO n(X)$ is a finite cover, where $\freeO n(-)$ is the $n$-fold free loop space functor. The functor $\freeO n(-)$ preserves pullbacks and hence composition of spans in $\Cov$. It follows that $\freeO n(-)$ is an endofunctor of $\Cov$.

Now one may wonder if the evaluation map $(S^1)^n \times \freeO n(-) \xto{\ev} \Id$ is a natural transformation of endofunctors of $\Cov$. There is no reason to expect this to be the case. In Example \ref{exampleDoesNotCommute}, we shall see that $(S^1)^n \x \freeO n(-) \xto{\ev} \Id$ is \emph{not} natural with respect to the backward maps in $\Cov$. The main purpose of this section is to replace $(S^1)^n \x \freeO n(-)$ by a functor $\twistO n(-)$ that is equal to $(S^1)^n \x \freeO n(-)$ on objects and forward maps in $\Cov$, differs from $(S^1)^n \x \freeO n(-)$ on finite sheeted coverings, and comes equipped with a natural transformation $\twistO n(-) \xto{\ev} \Id$.

\subsection{The construction of $\twistO n$}

Since the evaluation map is natural on $\Top \subset \Cov$, we will begin by studying the case of covering maps. Let $\pi \colon Y \tto X$ be a finite cover. Consider the pullback
\[
\xymatrix{\PB n(\pi) \ar[r]^{\cover\ev} \ar@{->>}[d]^{\ev^*\pi} & Y \ar@{->>}[d]^{\pi} \\ (S^1)^n \times \freeO n(X) \ar[r]^-{\ev}  & X.}
\]
Since $Y \tto X$ is a finite cover, the homotopy pullback and the pullback have the same homotopy type. A point in $\PB n(\pi)$ is given by a point $y \in Y$ and a pair $(\tup s,f) \in (S^1)^n \times \freeO n(X)$ such that $\pi(y) = f(\tup s)$ (i.e. the point $y \in Y$ lies in the fiber over $f(\tup s)$). Generally speaking, $\PB n(\pi)$ consists of a number of connected components. 

Assume $\tup s=(s_1,\dotsc,s_n)$ with $s_i\in S^1$, and let $q_i \colon I \to S^1$ for $1\leq i\leq n$ be the quotient map sending $0$ and $1$ to $s_i$. By path-lifting we can lift each map 
\[
f(s_1,\dotsc,s_{i-1},q_i(-),s_{i+1},\dotsc,s_n)\colon I \to X
\]
to a path in $Y$ starting at $y$ and ending at a (possibly different) point in the fiber over $f(\tup s)$.

\begin{definition}\label{defCovKMap}
We define a map $k\colon \PB n(\pi) \to (\ZZ_{>0})^n$ as follows. Given a triple $(y,\tup s,f)\in \PB n(\pi)$ the homotopy class $[f(s_1,\dotsc,s_{i-1},-,s_{i+1},\dotsc,s_n)]\in \pi_1(X,f(\tup s))$ in the fundamental group acts on the fiber over $f(\tup s)$. The positive integer $k(y,\tup s,f)_i$ is defined to be the size of the orbit of $y\in \pi^{-1}(f(\tup s))$ under the action of $[f(s_1,\dotsc,s_{i-1},-,s_{i+1},\dotsc,s_n)]$.
\end{definition}

\begin{remark}\label{remarkKMap}
Equivalently, $k(y,\tup s,f)_i$ is the smallest positive integer $k_i$ such that loop 
\[
f(s_1,\dotsc,s_{i-1},-,s_{i+1},\dotsc,s_n),
\] 
viewed as a loop with base point $f(\tup s)$, composed with itself $k_i$ times lifts to a loop in $Y$ when we start lifting the loop at $y\in \pi^{-1}(f(\tup s))$.

Yet another interpretation of $k(y,\tup s,f)_i$ is that this is the smallest positive integer $k$ such that the $k$th power of $[f(s_1,\dotsc,s_{i-1},-,s_{i+1},\dotsc,s_n)]\in \pi_1(X,f(\tup s))$ lies in the image of the injective group homomorphism $\pi_1(Y,y) \into \pi_1(X,f(\tup s))$.
\end{remark}

\begin{prop}\label{propKContinuousWithPathsProof}
The map $k\colon \PB n(\pi) \to (\ZZ_{>0})^n$ defined above is continuous.
\end{prop}
\begin{proof}
Since $\PB n(\pi)$ is homotopy equivalent to a CW complex, the path components map $\PB n(\pi) \to \pi_0(\PB n(\pi))$ is continuous. We will show that $k$ factors through this map.

Let $f \colon I \to \PB n(\pi)$ be a path. It suffices to show that $k$ gives the same value on the endpoints of this path. We have an induced map $\gamma\colon I \to X$. Let $x_0$ and $x_1$ be the endpoints of $\gamma$ and let $y_0$ and $y_1$ be the chosen lifts in $Y$. The path $\gamma$ identifies $\pi_1(X,x_0)$ and $\pi_1(X,x_1)$ and the different lifts of $\gamma$ to paths in $Y$ identifies the fibers $\pi^{-1}(x_0)$ and $\pi^{-1}(x_1)$ as $\pi_1(X,x_0)$ and $\pi_1(X,x_1)$-sets, respectively, by a bijection sending $y_0$ to $y_1$. The positive integers $k(f(0))_i$ and $k(f(1))_i$ are, by definition, the size of the orbit of $y_0$ and $y_1$ under the action of elements in the fundamental groups that are identified. Thus these numbers are identical. 
\end{proof}

\begin{definition}
We define the subspace $\PBlift n(\pi)\subseteq \PB n(\pi)$ to consist of all the triples $(y,\tup s,f)\in \PB n(\pi)$ such that $k(y,\tup s, f)_i=1$ for all $1\leq i\leq n$, i.e. the preimage of the tuple $(1,\dotsc,1)\in \Z^n$ under the map $k$. Note that $\PBlift n(\pi)$ is then a collection of connected components in $\PB n(\pi)$.

Let $\PBbroke n(\pi)\subseteq \PB n(\pi)$ be the collection of the remaining components of $\PB n(\pi)$.
\end{definition}

\begin{remark}
In view of Remark \ref{remarkKMap}, a triple $(y,\tup s,f)\in \PB n(\pi)$ lies in $\PBlift n(\pi)$ if and only if the map $f\colon (S^1)^n\to X$ lifts to a map $(S^1)^n\to Y$ when we lift by pathlifting starting at $\tup s\in (S^1)^n$ and from the point $y$ in the fiber of $f(\tup s)$.
\end{remark}

When working with the group structure on $S^1$, we will work additively (ie. we will view $S^1$ as $\R/\Z$).

\begin{prop}\label{propPBliftHomeomorphic}
The space $\PBlift n(\pi)$ is homeomorphic to $(S^1)^n\times \freeO n(Y)$ via the map that takes a triple $(y,\tup s, f)\in \PBlift n(\pi)$ to the pair $(\tup s,g)\in (S^1)^n\times \freeO n(Y)$ with $g$ the unique lift of $f$ to $Y$ with $g(\tup s)=y$.
\end{prop}

\begin{proof}
Define a map $\alpha\colon (S^1)^n\x \freeO nY\to \PBlift n(\pi)$ by taking a pair $(\tup s,g)\in (S^1)^n\x \freeO nY$ to the triple $(g(\tup s),\tup s,\pi\circ g)\in\PBlift n(\pi)$. The map $\alpha$ is a continuous bijection since each triple $(y,\tup s,f)\in\PBlift n(\pi)$ satisfies that $f\in \freeO nX$ lifts to a unique $g\in\freeO nY$ starting at the point $y$ in the fiber of $f(\tup s)$.

Let $\beta_1\colon \PBlift n(\pi)\to (S^1)^n$ be the projection $(y,\tup s,f)\mapsto \tup s$ which is clearly continuous. Next consider the following homotopy lifting problem:
\[
\begin{tikzpicture}
\node (M) [matrix of math nodes] {
    \PBlift n(\pi) &[3cm] Y \\[2cm]
    \PBlift n(\pi)\x I^n & X. \\
};
\path[auto,arrow,->]
    (M-1-1) edge node{$(y,\tup s, f)\mapsto y$} (M-1-2)
            edge node[swap]{$(y,\tup s, f)\mapsto (y,\tup s, f, \tup 0)$} (M-2-1)
    (M-1-2) edge[->>] node{$\pi$} (M-2-2)
    (M-2-1) edge node[swap]{$(y,\tup s, f, \tup t)\mapsto f(\tup s + \tup t)$} (M-2-2)
            edge[dashed] node{$\exists \hat\beta_2$} (M-1-2)
;
\end{tikzpicture}
\]
The diagram commutes, so the homotopy lifting property of the covering map $\pi$ gives us a continuous map $\hat\beta_2\colon \PBlift n(\pi)\x I^n\to Y$. By adjunction we get a continuous map $\beta_2\colon \PBlift n(\pi)\to \Map(I^n,Y)$ that takes a triple $(y,\tup s,f)$ to the unique path $g$ in $Y$ satisfying $\pi(g(\tup t))=f(\tup s+\tup t)$ and with $g(\tup 0)=y$. Since $(y,\tup s,f)\in \PBlift n(\pi)$, $f$ lifts to an $n$-fold loop in $Y$ starting at $y$. Thus the map $\beta_2\colon \PBlift n(\pi)\to \Map(I^n,Y)$ lands in $\freeO nY$.

The torus $(S^1)^n$ acts continuously on $\freeO nY$ by the map $\mu\colon (S^1)^n\x \freeO nY \to \freeO nY$ given by
\[
\mu(\tup s,g) = (\tup t\mapsto g(\tup t-\tup s) ).
\]
The composite $\mu\circ (\beta_1,\beta_2)\colon \PBlift n(\pi)\to (S^1)^n\x \freeO n Y$ is continuous and takes a triple $(y,\tup s,f)\in \PBlift n(\pi)$ to the pair $(\tup s, g)\in (S^1)^n\x \freeO n Y$ where $g\in \freeO nY$ is the unique $n$-fold loop with $g(\tup s)=y$ and $\pi\circ g = f$. Hence $\mu\circ (\beta_1,\beta_2)$ is the inverse homeomorphism to $\alpha$ from the start of the proof.
\end{proof}

Given any point $s\in \R/\Z$ and an integer $\ell$, the map $t\mapsto \ell(t-s)+s$ wraps $S^1\cong \R/\Z$ around itself $\ell$ times with $s$ as a fixed point.
\begin{definition}\label{defWind}
We construct a map $\wind(\pi)\colon \PB n(\pi)\to \PBlift n(\pi)\cong (S^1)^n\x \freeO n (Y)$ as follows. Given a triple $(y,\tup s, f)\in \PB n(\pi)$, let $k(y,\tup s,f)=(k_1,\dotsc,k_n)\in (\Z_{>0})^n$ be the associated $n$-tuple of positive integers and define $\hat f\colon (S^1)^n \to X$ to be the continuous map
\[
\hat f(t_1, \dotsc, t_n) = f(k_1(t_1-s_1)+s_1, \dotsc, k_n(t_n-s_n)+s_n).
\]
Then $\hat f(\tup s)=f(\tup s)$ so $(y,\tup s, \hat f)$ is a point of $\PB n(\pi)$.
Furthermore, each loop\linebreak $\hat f(s_1,\dotsc, s_{i-1}, -, s_{i+1}, \dotsc, s_n)$ in $X$ is the corresponding loop for $f$ composed with itself $k_i$ times. Hence $k(y,\tup s, \hat f)_i=1$ for all $1\leq i\leq n$, and $(y,\tup s,\hat f)\in \PBlift n(\pi)$.

We define $\wind(\pi)(y,\tup s, f) = (y,\tup s,\hat f)$. Since $k(y,\tup s,f)$ depends continuously on $(y,\tup s,f)\in \PB n(\pi)$, so does $\hat f$, and thus $\wind(\pi)\colon \PB n(\pi)\to \PBlift n(\pi)\cong (S^1)^n\x \freeO n(Y)$ is continuous.
\end{definition}

\begin{remark}\label{remarkWindRespectsEval}
The map $\wind(\pi)\colon \PB n(\pi) \to (S^1)^n\x \freeO n(Y)$ commutes with the projection $\cover{\ev_X}\colon \PB n(\pi)\to Y$ and evaluation $\ev_Y\colon (S^1)^n\x \freeO n(Y)\to Y$.
\end{remark}

Given a span $X \xffrom{\pi} Y \xto{f} Z$ in $\Cov$ we can use $\wind(\pi)$ to define a new span $\twistO n(\pi, f) = (\ev^*(\pi), ((S^1)^n\x \freeO n(f)) \circ \wind(\pi))$ by the composite
\[
(S^1)^n\x \freeO n(X) \xffrom{\ev^*(\pi)} \PB n(\pi) \xrightarrow{\wind(\pi)}  (S^1)^n\x \freeO n(Y) \xrightarrow{(S^1)^n\x \freeO n(f)} (S^1)^n\x \freeO n(Z). 
\]
It turns out that this construction gives an endofunctor on the category $\Cov$.

\begin{prop}\label{propCovTwistedLoopFunctor}
The construction above is an endofunctor $\twistO n \colon \Cov \to \Cov$. 
\end{prop}

\begin{proof}
Because every morphism in $\Cov$ factors as a covering map $X \xffrom\pi Y \xto{\id} Y$ followed by a map in $\Top$ (ie. maps of the form $X \xffrom\id X \xto f Y$), it is sufficient to show that $\twistO n$ respects composition for the four different ways of combining these generators. The first two cases are easy, the last two cases require some work.

\emph{Covering map followed by forward map:} The composition in $\Cov$ of a covering $(\pi,\id)\colon X \xffrom\pi Y \xto{\id} Y$ followed by a forward map $(\id,f)\colon Y \xffrom\id Y \xto f Z$ is simply the span $(\pi,f)\colon X \xffrom\pi Y \xto f Z$. The construction $\twistO n$ respects this composition essentially by construction:
\begin{align*}
 \twistO n(\id,f)\circ \twistO n(\pi,\id) 
 ={}& 
\begin{tikzpicture}[ampersand replacement=\&, baseline=(BASE).base]
\scriptsize
\node (M) [matrix of math nodes] {
    \&\& ? \&\& \\[1cm]
    \&\PB n(\pi) \&\& (S^1)^n \x \freeO n(Y) \& \\[1cm]
    (S^1)^n\x \freeO n(X) \&\& (S^1)^n \x \freeO n(Y) \&\& (S^1)^n \x \freeO n(Z) \\
};
\path [auto, arrow, ->] 
    (M-1-3) edge[->>,dashed] (M-2-2)
            edge[dashed] (M-2-4)
    (M-2-2) edge[->>] node{$\ev^*\pi$} (M-3-1)
            edge node{$\wind(\pi)$} (M-3-3)
    (M-2-4) edge[->>] node(BASE){$\id$} (M-3-3)
            edge node{$(S^1)^n \x \freeO n(f)$} (M-3-5)
;
\end{tikzpicture}
\\ ={}& 
\begin{tikzpicture}[ampersand replacement=\&, baseline=(BASE.base)]
\scriptsize
\node (M) [matrix of math nodes] {
    \&\PB n(\pi) \& \\[1cm]
    (S^1)^n\x \freeO n(X) \&\& (S^1)^n \x \freeO n(Z) \\
};
\path [auto, arrow, ->] 
    (M-1-2) edge[->>] node(BASE){$\ev^*\pi$} (M-2-1)
            edge node{$((S^1)^n \x \freeO n(f))\circ \wind(\pi)$} (M-2-3)
;
\end{tikzpicture}
\\ ={}& \twistO n(\pi,f).
\end{align*}

\emph{Composition of forward maps:}
Given two forward maps $(\id,f)\colon X \xffrom\id X \xto f Y$ and $(\id,g)\colon Y \xffrom\id Y \xto g Z$, we just have 
\[\twistO n(\id,g)\circ \twistO n(\id,f) =  \Bigl(
\begin{tikzpicture}[ampersand replacement=\&, baseline=(BASE.base)]
\scriptsize
\node (M) [matrix of math nodes] {
    \&(S^1)^n\x \freeO n(X) \& \\[1cm]
    (S^1)^n\x \freeO n(X) \&\& (S^1)^n \x \freeO n(Z) \\
};
\path [auto, arrow, ->] 
    (M-1-2) edge[->>] node(BASE){$\id$} (M-2-1)
            edge node{$(S^1)^n \x \freeO n(g\circ f)$} (M-2-3)
;
\end{tikzpicture}
\Bigr) = \twistO n(\id,g\circ f).\]

\emph{Composition of covering maps:} Given two spans of the form $X\xffrom\pi Y\xto\id Y$ and $Y\xffrom \rho Z\xto\id Z$ in $\Cov$, we consider the composition of spans $\twistO n(\rho,\id)\circ\twistO n(\pi,\id)$:
\[
\begin{tikzpicture}[ampersand replacement=\&, baseline=(BASE).base]
\scriptsize
\node (M) [matrix of math nodes] {
    \&\& ? \&\& \\[1cm]
    \&\PB n(\pi) \&\& \PB n(\rho) \& \\[1cm]
    (S^1)^n\x \freeO n(X) \&\& (S^1)^n \x \freeO n(Y) \&\& (S^1)^n \x \freeO n(Z). \\
};
\path [auto, arrow, ->] 
    (M-1-3) edge[->>,dashed] (M-2-2)
            edge[dashed] (M-2-4)
    (M-2-2) edge[->>] node{$\ev_X^*\pi$} (M-3-1)
            edge node{$\wind(\pi)$} (M-3-3)
    (M-2-4) edge[->>] node(BASE){$\ev_Y^*\rho$} (M-3-3)
            edge node{$\wind(\rho)$} (M-3-5)
;

\end{tikzpicture}
\]
We claim that the pullback in the middle gives us the space $\PB n(\pi\rho)$. To confirm this, we study the following diagram:
\[
\begin{tikzpicture}[ampersand replacement=\&, baseline=(BASE).base]

\node (M) [matrix of math nodes] {
    ? \&[1cm] \PB n(\rho) \&[1cm] Z \\[1cm]
    \PB n(\pi) \&  (S^1)^n \x \freeO n(Y)\& Y \\[1cm]
    (S^1)^n\x \freeO n(X) \&\& X. \\
};
\path [auto, arrow, ->] 
    (M-1-1) edge[dashed] (M-1-2)
            edge[->>,dashed] (M-2-1)
    (M-2-1) edge[->>] node{$\ev_X^*\pi$} (M-3-1)
            edge node(BASE){$\wind(\pi)$} (M-2-2)
    (M-1-2) edge[->>] node{$\ev_Y^*\rho$} (M-2-2)
            edge node{$\cover{\ev_Y}$} (M-1-3)
    (M-1-3) edge[->>] node{$\rho$} (M-2-3)
    (M-2-2) edge node{$\ev_Y$} (M-2-3)
    (M-3-1) edge node{$\ev_X$} (M-3-3)
    (M-2-3) edge[->>] node{$\pi$} (M-3-3)
;

\end{tikzpicture}
\]
According to Remark \ref{remarkWindRespectsEval}, the composite $\ev_Y\circ \wind(\pi)\colon \PB n(\pi) \to Y$ is just the usual projection $\cover{\ev_X}\colon \PB n(\pi) \to Y$. Both the bottom square and the top square of the diagram above are therefore pullback squares. Hence the big outer square has to be a pullback, so the top-left corner must be $\PB n(\pi\rho)$ as claimed.

Furthermore, the map $\PB n(\pi\rho)\tto \PB n(\pi)$ takes a triple $(z,\tup s, f)\in \PB n(\pi\rho)$ to the triple $(\rho(z),\tup s,f)\in \PB n(\pi)$, and the composite $\PB n(\pi\rho)\tto \PB n(\pi)\tto (S^1)^n\x\freeO n(X)$ is just $\ev_X^*(\pi\rho)$.

Similarly, we can describe the map $\PB n(\pi\rho)\to \PB n(\rho)$ if we recall Definition \ref{defWind}, the definition of $\wind(\pi)$.
Given a triple $(z,\tup s, f)\in \PB n(\pi\rho)$ we first have $(\rho(z),\tup s,f)\in \PB n(\pi)$. Let $f'\in \freeO n(Y)$ be the (unique) loop such that $f'(\tup s)=\rho(z)$ and
\[
\pi(f'(\tup t)) = f(k(\rho(z),\tup s, f)_1\cdot (t_1-s_1)+s_1, \dotsc, k(\rho(z),\tup s,f)_n\cdot (t_n-s_n)+s_n).
\]
Then $\wind(\pi)(\rho(z),\tup s, f) = (\tup s,f')\in (S^1)^n\x \freeO n(Y)$, and the map $\PB n(\pi\rho)\to \PB n(\rho)$ takes $(z,\tup s, f)$ to $(z,\tup s, f')$.

The composite $\twistO n(\rho,\id)\circ\twistO n(\pi,\id)$ is therefore given by the span \[
(S^1)^n\x \freeO n(X) \xffrom{\ev_X^*(\pi\rho)} \PB n(\pi\rho) \xto{(z,\tup s,f)\mapsto (z,\tup s, f')} \PB n(\rho) \xto{\wind(\rho)} (S^1)^n \x \freeO n(Z).
\]
It remains to check that the composite map $\PB n(\pi\rho) \to \PB n(\rho) \to (S^1)^n \x \freeO n(Z)$ coincides with $\wind(\pi\rho)$. For this recall that $\wind(\rho)$ takes the triple $(z,\tup s, f')\in \PB n(\rho)$ to the pair $(\tup s, f'')\in (S^1)^n\x \freeO n(Z)$, where $f''$ satisfies $f''(\tup s)=z$ and
\begin{align*}
\rho(f''(\tup t)) &= f'(k(z,\tup s, f')_1\cdot (t_1-s_1)+s_1, \dotsc, k(z,\tup s,f')_n\cdot (t_n-s_n)+s_n),
\end{align*}
or equivalently
\begin{align*}
\pi\rho(f''(\tup t)) ={}& \pi \bigl(f'(k(z,\tup s, f')_1\cdot (t_1-s_1)+s_1, \dotsc, k(z,\tup s,f')_n\cdot (t_n-s_n)+s_n)\bigr)
\\ ={}& f(k(\rho(z),\tup s,f)_1\cdot k(z,\tup s,f')_1\cdot (t_1-s_1)+s_1,
\dotsc, 
\\ &\hspace{3.5cm} k(\rho(z),\tup s,f)_n\cdot k(z,\tup s,f')_n\cdot (t_n-s_n)+s_n).
\end{align*}
The positive integer $k(\rho(z),\tup s,f)_i$ is the smallest exponent such that the corresponding power of $[f(s_1,\dotsc,s_{i-1},-,s_{i+1},\dotsc,s_n)]\in \pi_1(X,\pi\rho (z))$ lies in the subgroup $\pi_1(Y,\rho (z))\into\pi_1(X,\pi\rho (z))$, and the resulting power is then the class 
\[
[f'(s_1,\dotsc,s_{i-1},-,s_{i+1},\dotsc,s_n)]\in\pi_1(Y,\rho (z)). 
\]
Similarly, $k(z,\tup s,f')_i$ is the smallest exponent such that the corresponding power of the class $[f'(s_1,\dotsc,s_{i-1},-,s_{i+1},\dotsc,s_n)]$ in turn lies in the even smaller subgroup $\pi_1(Z,z)\into \pi_1(Y,\rho (z))\into \pi_1(X,\pi\rho (z))$.

Hence the product $k(\rho(z),\tup s,f)_i\cdot k(z,\tup s,f')_i$ is the smallest exponent such that the power of $[f(s_1,\dotsc,s_{i-1},-,s_{i+1},\dotsc,s_n)]\in \pi_1(X,\pi\rho (z))$ lies in the subgroup $\pi_1(Z,z)\into \pi_1(X,\pi\rho (z))$, i.e.
\[k(\rho(z),\tup s,f)_i\cdot k(z,\tup s,f')_i = k(z,\tup s,f)_i\]
for each $1\leq i\leq n$. Consequently,
\[\pi\rho(f''(\tup t)) = f(k(z,\tup s, f)_1\cdot (t_1-s_1)+s_1, \dotsc, k(z,\tup s,f)_n\cdot (t_n-s_n)+s_n)\]
so the map $(z,\tup s,f)\mapsto (\tup s,f'')$ from $\PB n(\pi\rho) \to (S^1)^n \x \freeO n Z$ equals $\wind(\pi\rho)$ as claimed.

\emph{Forward map followed by covering map:}
Consider the spans $X \xffrom\id X \xto f Y$ and $Y\xffrom\pi Z \xto\id Z$. The composite of these is given by the pullback
\[
\begin{tikzpicture}[ampersand replacement=\&, baseline=(BASE).base]
\scriptsize
\node (M) [matrix of math nodes] {
    \&[.5cm]\&[.5cm] W \&[.5cm]\&[.5cm] \\[1cm]
    \& X \&\& Z \& \\[1cm]
    X \&\& Y \&\& Z. \\
};
\path [auto, arrow, ->] 
    (M-1-3) edge[->>] node{$f^*\pi$} (M-2-2)
            edge node{$\cover f$} (M-2-4)
    (M-2-2) edge[->>] node{$\id$} (M-3-1)
            edge node{$f$} (M-3-3)
    (M-2-4) edge[->>] node(BASE){$\pi$} (M-3-3)
            edge node{$\id$} (M-3-5)
;
\end{tikzpicture}
\]
Hence $\twistO n((\pi,\id)\circ (\id,f)) = \twistO n(f^*\pi,\cover f)$ is the span
\begin{equation}\label{eqTwistedLnOfSpan}
(S^1)^n\x\freeO n X \xffrom{\ev_X^*f^*\pi} \PB n(f^*\pi) \xto{((S^1)^n\x \freeO n(\cover f))\circ \wind(f^*\pi)} (S^1)^n\x\freeO nZ.
\end{equation}
On the other hand, the composite $\twistO n(\pi,\id)\circ \twistO n(\id,f)$ is given by
\[
\begin{tikzpicture}[ampersand replacement=\&, baseline=(BASE).base]
\scriptsize
\node (M) [matrix of math nodes] {
    \&[.5cm]\&[.5cm] ? \&[.5cm]\&[.5cm] \\[1cm]
    \& (S^1)^n\x\freeO nX \&\& \PB n(\pi) \& \\[1cm]
    (S^1)^n\x\freeO nX \&\& (S^1)^n\x\freeO nY \&\& (S^1)^n\x\freeO nZ. \\
};
\path [auto, arrow, ->] 
    (M-1-3) edge[->>] node[swap]{$((S^1)^n\x\freeO n(f))^*\ev_Y^*\pi$} (M-2-2)
            edge node{$\cover{(S^1)^n\x\freeO n(f)}$} (M-2-4)
    (M-2-2) edge[->>] node{$\id$} (M-3-1)
            edge node{$(S^1)^n\x\freeO n(f)$} (M-3-3)
    (M-2-4) edge[->>] node(BASE){$\ev_Y^*\pi$} (M-3-3)
            edge node{$\wind(\pi)$} (M-3-5)
;
\end{tikzpicture}
\]
The forward map $f$ commutes with the evaluations $\ev_X$ and $\ev_Y$, i.e. the square
\[
\begin{tikzpicture}
\node (M) [matrix of math nodes] {
    (S^1)^n\x\freeO nX &[1cm] X \\[1cm]
    (S^1)^n\x\freeO nY & Y \\
};
\path[auto,arrow,->]
    (M-1-1) edge node{$\ev_X$} (M-1-2)
            edge node[swap]{$(S^1)^n\x\freeO n(f)$} (M-2-1)
    (M-1-2) edge node{$f$} (M-2-2)
    (M-2-1) edge node{$\ev_Y$} (M-2-2)
;
\end{tikzpicture}
\]
commutes. It follows that the different ways of pulling back the cover $\pi$ over $Y$ agree to give the same cover over $(S^1)^n\x\freeO nX$, i.e. $((S^1)^n\x\freeO n(f))^*\ev_Y^*\pi = \ev_X^*f^*\pi$, so the composite $\twistO n(\pi,\id)\circ \twistO n(\id,f)$ can instead be written as
\[
\begin{tikzpicture}[ampersand replacement=\&, baseline=(BASE).base]
\scriptsize
\node (M) [matrix of math nodes] {
    \&[.5cm]\&[.5cm] \PB n(f^*\pi) \&[.5cm]\&[.5cm] \\[1cm]
    \& (S^1)^n\x\freeO nX \&\& \PB n(\pi) \& \\[1cm]
    (S^1)^n\x\freeO nX \&\& (S^1)^n\x\freeO nY \&\& (S^1)^n\x\freeO nZ. \\
};
\path [auto, arrow, ->] 
    (M-1-3) edge[->>] node[swap]{$\ev_X^*f^*\pi$} (M-2-2)
            edge node{$\cover{(S^1)^n\x\freeO n(f)}$} (M-2-4)
    (M-2-2) edge[->>] node{$\id$} (M-3-1)
            edge node{$(S^1)^n\x\freeO n(f)$} (M-3-3)
    (M-2-4) edge[->>] node(BASE){$\ev_Y^*\pi$} (M-3-3)
            edge node{$\wind(\pi)$} (M-3-5)
;
\end{tikzpicture}
\]
This matches the span $\twistO n(f^*\pi,\cover f)$ in equation \eqref{eqTwistedLnOfSpan} assuming that we can show that the following diagram of forward maps commutes:
\[
\begin{tikzpicture}
\node (M) [matrix of math nodes] {
    \PB n(f^*\pi) &[2cm] (S^1)^n\x\freeO nW \\[1cm]
    \PB n(\pi) & (S^1)^n\x\freeO nZ. \\
};
\path[auto,arrow,->]
    (M-1-1) edge node{$\wind(f^*\pi)$} (M-1-2)
            edge node[swap]{$\cover{(S^1)^n\x\freeO n(f)}$} (M-2-1)
    (M-1-2) edge node{$(S^1)^n\x\freeO n(\cover f)$} (M-2-2)
    (M-2-1) edge node{$\wind(\pi)$} (M-2-2)
;
\end{tikzpicture}
\]
The space $W$ consists of those pairs $(x,z)\in X\x Z$ such that $f(x)=\pi(z)$. Similarly the $n$-fold free loop space $\freeO nW$ consists of pairs of loops $(\gamma,\zeta)\in \freeO nX \x \freeO n Z$ such that $f\circ \gamma=\pi\circ \zeta$.

Now $\PB n(f^*\pi)$ consists of quadruples $(x,z,\tup s, \gamma)$ with $(x,z)\in W$, $\tup s\in (S^1)^n$, and $\gamma\in \freeO nX$ satisfying
\[\gamma(\tup s) = (f^*\pi)(x,z) = x.\]
Hence $x$ equals $\gamma(\tup s)\in X$. 

The map $\cover{(S^1)^n\x\freeO n(f)}\colon \PB n(f^*\pi)\to \PB n(f)$ takes a quadruple $(x,z,\tup s,\gamma)\in \PB n(f^*\pi)$ to the triple $(z,\tup s,f\circ \gamma)\in \PB n(\pi)$.

A loop $\gamma\in \freeO n X$ lifts to a pair $(\gamma,\zeta)\in \freeO n W$ precisely whenever the loop $f\circ \gamma\in \freeO n Y$ lifts to a loop $\zeta\in\freeO n Z$. Consequently, the $k$ exponent maps agree
\[k(x,z,\tup s,\gamma)_i = k(z,\tup s,f\circ\gamma)_i\]
for all $1\leq i\leq n$.

The map $\wind(f^*\pi)$ winds $\gamma$ by the exponents $k(x,z,\tup s,\gamma)_i$ to get $\hat \gamma\in \freeO n(X)$ of Definition \ref{defWind} such that $(x,z,\tup s,\hat\gamma)\in \PBlift n (f^*\pi)$ and $f\circ \hat\gamma$ lifts to $\zeta\in\freeO n Z$. Subsequently the map $(S^1)^n\x \freeO n(\cover f)$ maps the triple $(\tup s,\hat \gamma,\zeta)\in (S^1)^n\x \freeO n W$ to the pair $(\tup s,\zeta)\in(S^1)^n\x \freeO n Z$.

At the same time, $\wind(\pi)$ uses the same exponents $k(z,\tup s,f\circ \gamma)_i$ as before to wind $f\circ \gamma$ and get $\widehat{f\circ\gamma}\in \freeO n(Y)$ with $(z,\tup s,\widehat{f\circ\gamma})\in \PBlift n(\pi)$. Since $\wind(\pi)$ and $\wind(f^*\pi)$ twist by the same exponents, we have $\widehat{f\circ \gamma}=f\circ \hat \gamma$ and therefore we get the same lift to $\zeta\in \freeO nZ$. We conclude that
\begin{align*}
{}& (\wind(\pi)\circ(\cover{(S^1)^n\x \freeO n(f)}))(x,z,\tup s,\gamma)
\\ ={}& \wind(\pi)(z,\tup s,f\circ\gamma) 
\\ ={}& (\tup s,\zeta) 
\\ ={}& ((S^1)^n\x \freeO n(\cover f))(\tup s,\hat \gamma,\zeta) 
\\ ={}&   (((S^1)^n\x \freeO n(\cover f))\circ \wind(f^*\pi))(x,z,\tup s,\gamma)
\end{align*}
as required.

This completes the fourth and final case in checking that $\twistO n$ preserves composition of spans in $\Cov$.
\end{proof}

There is a natural alternative to $\twistO n$ that one might consider. Consider a span $X \xffrom{\pi} Y \xto{f} Z$ in $\Cov$. 
Instead of massaging $\PB n(\pi)$ onto the subspace $\PBlift n(\pi)$ by using $\wind(\pi)$, we can instead just throw away all the components in $\PBbroke n(\pi)$ from $\PB n(\pi)$. Concretely this means that we restrict the span $\twistO n(\pi,f)$, which goes through $\PB n(\pi)$, to get a span through $\PBlift n(\pi)$:
\[
\begin{tikzpicture}[ampersand replacement=\&, baseline=(BASE).base]
\scriptsize
\node (M) [matrix of math nodes] {
    \&[1cm]\&[1cm] \PBlift n(\pi) \&[3cm]\&[1cm] \\[1cm]
    \& \&\PB n(\pi)\&  \& \\[1cm]
     (S^1)^n\x \freeO n X\&\&  \&\& (S^1)^n\x\freeO n Z. \\
};
\path [auto, arrow, ->] 
    (M-1-3) edge[->>] node[swap]{$(\ev^*\pi)|_{\PBlift n(\pi)}$} (M-3-1)
            edge node(BASE){$((S^1)^n\x \freeO n f) \circ \wind(\pi)|_{\PBlift n(\pi)}$} (M-3-5)
            edge[left hook->] (M-2-3)
    (M-2-3) edge[->>] node{$\ev^*\pi$} (M-3-1)
            edge node[swap,pos=.8]{$((S^1)^n\x \freeO n f) \circ \wind(\pi)$} (M-3-5)
;
\end{tikzpicture}
\]
The restriction of $\ev^*\pi$ to $\PBlift n(\pi)$ is still a covering map -- the fiber over an $n$-fold loop $f\in \freeO n X$ is just empty if $f$ has no lift to an $n$-fold loop in $Y$.
By construction, $\wind(\pi)|_{\PBlift n(\pi)}=\id_{\PBlift n(\pi)}$ so that the forward map is just $(S^1)^n\x \freeO n f$.

\begin{definition}\label{defCovUntwistedLn}
Given a span $X \xffrom{\pi} Y \xto{f} Z$ in $\Cov$, define $\untwistO n(\pi,f)$ to be the span
\[
(S^1)^n\x \freeO n(X) \xffrom{(\ev^*(\pi))|_{\PBlift n(\pi)}} \PBlift n(\pi)\cong (S^1)^n\x \freeO n(Y) \xrightarrow{(S^1)^n\x \freeO n(f)} (S^1)^n\x \freeO n(Z). 
\]
\end{definition}

\begin{remark}\label{remarkCovUntwisted}
The construction $\untwistO n$ takes covering maps to covering maps and forward maps to forward maps. The covering map $(S^1)^n\x\freeO n(Y)\cong\PBlift n(\pi) \xtto{\ev^*(\pi)} (S^1)^n\x \freeO n(X)$ coincides with the covering map $(S^1)^n\x\freeO n(\pi)$, so the construction $\untwistO n$ turns out to just recover the functor $(S^1)^n\x \freeO n(-)$ on $\Cov$. While $\ev \colon (S^1)^n\x \freeO n(X) \to X$ is not a natural transformation from $\untwistO n$ to $\Id$ (see Example \ref{exampleDoesNotCommute} below), the functor $\untwistO n$ does play a role in character theory.
\end{remark}

\begin{example}\label{exampleDoesNotCommute}
We explain why the functor $(S^1)^n\x\freeO n(-)\colon \Cov \to \Cov$ does not commute with the evaluation map.

Let $\pi\colon Y\to X$ be a non-trivial finite-sheeted covering map (ie. not just the projection of a product $Y\cong X\x F \tto X$, where $F$ is a finite set). Since the cover is non-trivial there exists some loop $\gamma\in \freeO 1 X$ and some point $y$ in the fiber of $x:=f(0)$ such that $\gamma$ does not lift to a loop in $Y$ starting at $y$. Therefore the covering map $\freeO 1(\pi)\colon \freeO 1Y \tto \freeO 1X$ has smaller fiber over $\gamma\in \freeO 1X$ than the fiber of $\pi$ over $x$ because $\freeO 1(\pi)^{-1}(\gamma)$ has no point corresponding to $y\in \pi^{-1}(x)$. Note that $(y,0,\gamma)\in \PB 1(\pi)$, and since $\gamma$ does not lift to a loop starting at $y$, we have $(y,0,\gamma)\not\in \PBlift 1(\pi)$.

We conclude that $\PBlift 1(\pi)$ is a proper subspace of $\PB 1(\pi)$, and as such the square
\[
\begin{tikzpicture}[ampersand replacement=\&, baseline=(BASE).base]
\node (M) [matrix of math nodes] {
    \&[.5cm]\&[.5cm] \PBlift 1(\pi) \&[.75cm]\&[.5cm] \\[1cm]
    \& S^1\x\freeO 1X \&\& Y \& \\[1cm]
    \&\& X \&\&  \\
};
\path [auto, arrow, ->] 
    (M-1-3) edge[->>] node[swap]{$(\ev_X^*\pi)|_{\PBlift 1(\pi)}$} (M-2-2)
            edge node{$\cover{\ev_X}|_{\PBlift 1(\pi)}$} (M-2-4)
    (M-2-2) edge node[swap]{$\ev_X$} (M-3-3)
    (M-2-4) edge[->>] node(BASE){$\pi$} (M-3-3)
;
\end{tikzpicture}
\]
is not a pullback, and hence not commutative in $\Cov$ when $\pi$ and $\ev_X^*\pi$ are considered to be backward maps. By Proposition \ref{propPBliftHomeomorphic}, $\PBlift 1(\pi)$ is homeomorphic to $S^1\x \freeO 1Y$ as coverings over $S^1\x \freeO 1X$, and compatibly with the evaluation maps. The square above can thus also be written as
\[
\begin{tikzpicture}[ampersand replacement=\&, baseline=(BASE).base]
\node (M) [matrix of math nodes] {
    \&[.5cm]\&[.5cm] S^1\x\freeO 1Y \&[.75cm]\&[.5cm] \\[1cm]
    \& S^1\x\freeO 1X \&\& Y \& \\[1cm]
    \&\& X \&\&  \\
};
\path [auto, arrow, ->] 
    (M-1-3) edge[->>] node[swap]{$\id_{S^1}\x\freeO 1(\pi)$} (M-2-2)
            edge node{$\ev_Y$} (M-2-4)
    (M-2-2) edge node[swap]{$\ev_X$} (M-3-3)
    (M-2-4) edge[->>] node(BASE){$\pi$} (M-3-3)
;
\end{tikzpicture}
\]
and is still not a pullback, so does not commute in $\Cov$.
\end{example}

\subsection{Properties of $\twistO n$}

In this section we state and prove Theorem \ref{thmIntroMain}.

\begin{theorem}\label{thmCovLiftingFunctor}\label{thmCovMain}
For $n \geq 0$, the endofunctors $\twistO n \colon \Cov\to \Cov$ of Proposition \ref{propCovTwistedLoopFunctor} have the following properties:
\begin{enumerate}
\renewcommand{\theenumi}{$(\roman{enumi})$}\renewcommand{\labelenumi}{\theenumi}
\item\label{itemCovLnZero} $\twistO 0$ is the identity functor on $\Cov$.
\item\label{itemCovLnObjects} On objects, $\twistO n$ takes a space $X$ to 
\[\twistO n(X) = (S^1)^n\x \freeO n (X).\]
\item\label{itemCovEquivariant} $\twistO n$ is equivariant with respect to the $\Sigma_n$-action on $(S^1)^n\x \freeO n(-)$ that permutes the coordinates of both $(S^1)^n$ and $\freeO n(-)$, i.e. for every $\sigma\in\Sigma_n$ the diagonal action of $\sigma$ on $(S^1)^n\x \freeO n(-)$ induces a natural isomorphism $\sigma\colon \twistO n\overset\cong\Rightarrow \twistO n$.
\item\label{itemCovLnForwardMaps} On forward maps, i.e. spans $X\xffrom\id X\xto f Y$, the functor $\twistO n$ coincides with the torus times the usual free loop space functor in $\Top$ so that $\twistO n(\id,f)=(S^1)^n\x \freeO n(f)$.
\item\label{itemCovEvalSquare} For all $n \geq 0$, the functor $\twistO n$ commutes with evaluation maps, i.e. the evaluation maps $\ev_X\colon (S^1)^n\x \freeO n(X)\to X$ form a natural transformation $\ev\colon \twistO n \Rightarrow \Id_{\Cov}$.
\item\label{itemCovLnPartialEvaluation} For all $n \geq 0$, the partial evaluation maps $\pev_X\colon S^1\x \freeO {n+1}(X) \to \freeO n (X)$ given by 
\[
\pev(t,f)= (\tup s\mapsto f(\tup s,t)), \qquad \text{for } t\in S^1, \tup s\in (S^1)^{n},
\] give natural transformations $(S^1)^n\x \pev\colon \twistO {n+1} \Rightarrow \twistO n$.
\item\label{itemCovIterateLn} For all $n,m\geq 0$, and any space $X$, the space $(S^1)^{n+m}\x \freeO {n+m} (X)$ embeds into $(S^1)^m\x \freeO m((S^1)^n\x \freeO n X)$ as those $m$-fold loops in $(S^1)^n\x \freeO n (X)$ that are constant in the $(S^1)^n$-coordinate, i.e. the embedding is given by
\[((\tup s,\tup r), f) \mapsto \Bigl(\tup r, \tup r'\mapsto \bigl(\tup s, \tup s'\mapsto f(s',r')\bigr)\Bigr)\]
for $\tup s,\tup s'\in (S^1)^n$, $\tup r,\tup r'\in (S^1)^m$ and $f\in \freeO {n+m} X$.

These embeddings then form a natural transformation $\twistO {n+m}(-)\Rightarrow \twistO m(\twistO n(-))$.
\end{enumerate}
\end{theorem}

\begin{proof}
\ref{itemCovLnZero}-\ref{itemCovLnObjects}: True by inspection of the definition of $\twistO n$ leading up to Proposition \ref{propCovTwistedLoopFunctor}.

\ref{itemCovEquivariant}: Let $X\xffrom\pi Y\xto h Z$ be a span in $\Cov$, and let $\sigma\in \Sigma_n$ be any permutation. Permuting the coordinates of $(S^1)^n\x\freeO n(X)$ by $\Sigma_n$ also permutes the loops $f(s_1,\dotsc,s_{i-1},-,s_{i+1},\dotsc, s_n)\colon S^1\to X$ and thus permutes the exponents $k(y,\tup s,f)_i$. It follows that $\wind(\pi)$ respects the action of $\sigma$, and the following diagram commutes as required:
\[
\begin{tikzpicture}[ampersand replacement=\&, baseline=(BASE).base]
\node (M) [matrix of math nodes] {
    (S^1)^n\x \freeO n X \&[1cm] \PB n(\pi) \&[1cm] (S^1)^n\x\freeO n Y \&[2cm] (S^1)^n\x\freeO n Z \\[1cm]
    (S^1)^n\x \freeO n X \& \PB n(\pi) \& (S^1)^n\x\freeO n Y \& (S^1)^n\x\freeO n Z.  \\
};
\path [auto, arrow, ->] 
    (M-1-1) edge node{$\sigma$} node[swap]{$\cong$} (M-2-1)
    (M-1-2) edge[->>] node[swap]{$\ev_X^*\pi$} (M-1-1)
            edge node{$\sigma$} node[swap]{$\cong$} (M-2-2)
            edge node{$\wind(\pi)$} (M-1-3)
    (M-1-3) edge node{$(S^1)^n\x\freeO n(h)$} (M-1-4)
            edge node{$\sigma$} node[swap]{$\cong$} (M-2-3)
    (M-1-4) edge node{$\sigma$} node[swap]{$\cong$} (M-2-4)
    (M-2-2) edge[->>] node{$\ev_X^*\pi$} (M-2-1)
            edge node[swap]{$\wind(\pi)$} (M-2-3)
    (M-2-3) edge node[swap]{$(S^1)^n\x\freeO n(h)$} (M-2-4)
;
\end{tikzpicture}
\]

\ref{itemCovLnForwardMaps}: For a span $X\xffrom\id X \xto f Y$ where the covering map is just $\id\colon X\tto X$, we have $\PB n(\id)=\PBlift n(\id)=(S^1)^n\x\freeO nX$ and $\wind(\id)=\id_{(S^1)^n\x\freeO nX}$. Hence $\twistO n(\id,f)$ is just the forward map $(S^1)^n\x\freeO n(f)$.

\ref{itemCovEvalSquare}: Given a span $X\xffrom\pi Y\xto f X$, consider the following diagram:
\[
\begin{tikzpicture}[ampersand replacement=\&, baseline=(BASE).base]
\node (M) [matrix of math nodes] {
    \PBlift n(\pi) \&[1cm] (S^1)^n\x\freeO n Y \&[2cm] (S^1)^n\x\freeO n Z \\[1cm]
    \PB n(\pi) \&  Y \& Z \\[1cm]
    (S^1)^n\x \freeO n(X) \&X.\& \\
};
\path [auto, arrow, ->] 
    (M-1-1) edge node{$\cong$} (M-1-2)
            edge node[xshift=-.1cm,yshift=-.1cm]{$\cover{\ev_X}$} (M-2-2)
    (M-1-2) edge node{$(S^1)^n\x\freeO n(f)$} (M-1-3)
            edge node{$\ev_Y$} (M-2-2)
    (M-1-3) edge node{$\ev_Z$} (M-2-3)
    (M-2-1) edge node{$\wind(\pi)$} (M-1-1)
            edge node[swap]{$\cover{\ev_X}$} (M-2-2)
            edge[->>] node{$\ev_X^*\pi$} (M-3-1)
    (M-2-2) edge node(BASE){$f$} (M-2-3)
            edge[->>] node{$\pi$} (M-3-2)
    (M-3-1) edge node{$\ev_X$} (M-3-2)
;
\end{tikzpicture}
\]
The diagram commutes in $\Top$ by straightforward calculation of the maps involved. The lower square is a pullback and thus commutes in $\Cov$, hence the entire diagram commutes in $\Cov$ when read as a diagram of backward and forward maps from $(S^1)^n\x\freeO n(X)$ to $Z$. The outer span from $(S^1)^n\x\freeO nX$ to $(S^1)^n\x\freeO n Z$ is simply the span $\twistO n(\pi,f)$. Hence the diagram shows that $\ev_Y\circ \twistO n(\pi,f) = (\pi,f)\circ \ev_X$ in the category $\Cov$.

\ref{itemCovIterateLn}: We postpone the proof of \ref{itemCovLnPartialEvaluation} for the moment and instead proceed with \ref{itemCovIterateLn}. Let $\iota''_3\colon (S^1)^{n+m}\x\freeO {n+m}(X)\x (S^1)^m \x (S^1)^n\to X$ be the continuous map given by
\[\iota''_3((\tup s,\tup r),f,\tup r',\tup s') = f(\tup s',\tup r').\]
The adjoint map $\iota'_3\colon (S^1)^{n+m}\x\freeO {n+m}(X)\x (S^1)^m \to \freeO n(X)$ is then also continuous, hence so is the map $\iota'_2\colon (S^1)^{n+m}\x\freeO {n+m}(X)\x (S^1)^m \to (S^1)^n\x \freeO n(X)$ given by
\[\iota'_2((\tup s,\tup r),f,\tup r') = \bigl(\tup s, \tup s'\mapsto f(\tup s',\tup r')\bigr).\]
Taking adjoints again we have the continuous map $\iota_2\colon (S^1)^{n+m}\x\freeO {n+m}(X) \to \freeO m((S^1)^n\x\freeO n(X))$. We can then conclude that $\iota_X\colon (S^1)^{n+m}\x\freeO {n+m}(X) \to (S^1)^m\x \freeO m((S^1)^n\x\freeO n(X))$ given by
\begin{equation}\label{eqIterationEmbedding}
((\tup s,\tup r), f) \mapsto \Bigl(\tup r, \tup r'\mapsto \bigl(\tup s, \tup s'\mapsto f(s',r')\bigr)\Bigr)
\end{equation}
is continuous.
To see that $\iota_X$ is an embedding, choose any basepoint $\ast\in (S^1)^m$, and let $\text{pr}_1\colon (S^1)^n\x\freeO n(X)\to (S^1)^n$ and $\text{pr}_2\colon (S^1)^n\x\freeO n(X)\to \freeO n(X)$ be the projections. As above, we can then give a continuous map $\rho_X\colon (S^1)^m\x \freeO m((S^1)^n\x\freeO n(X)) \to (S^1)^{n+m}\x\freeO {n+m}(X)$ by the expression
\[
\rho_X(\tup r,g) = \Bigl(\bigl(\text{pr}_1(g(\ast)),\tup r\bigr),\bigl((\tup s',\tup r')\mapsto \text{pr}_2(g(\tup r'))(\tup s')\bigr)\Bigr).
\]
By inspection one can check that $\rho_X\circ \iota_X$ is the identity on $(S^1)^{n+m}\x\freeO {n+m}(X)$, so $\rho_X$ is a retract of $\iota_X$, hence $\iota_X$ is an embedding.

Next we need to show that $\iota\colon \twistO {n+m}(-)\Rightarrow \twistO m(\twistO n(-))$ is a natural transformation. Suppose therefore that $X\xffrom \pi Y \xto h Z$ is an arbitrary span in $\Cov$.
The span $\twistO m(\twistO n(\pi,h))$ then takes the form
\[
\begin{tikzpicture}[ampersand replacement=\&, baseline=(BASE).base]
\scriptsize 
\matrix (M) [matrix of math nodes] {
  \&[-.5cm] \PB m(\ev_X^*\pi) \&[-.5cm]   \&[-.7cm]   \&[-.9cm]  \\[1cm]
 (S^1)^m\x \freeO m((S^1)^n\x\freeO n(X)) \&\& (S^1)^m\x \freeO m(\PB n(\pi)) \&\& \\[1cm]
 \&\&\& (S^1)^m\x \freeO m((S^1)^n \x \freeO n(Y)) \& \\[1cm]
 \&\&\&\& (S^1)^m\x \freeO m((S^1)^n\x \freeO n(Z)). \\
};
\path [auto, arrow, ->]
    (M-1-2) edge[->>] node[swap]{$\ev_{(S^1)^n\x\freeO n(X)}^*\ev_X^*\pi$} (M-2-1)
            edge node(BASE){$\wind(\ev_X^*\pi)$} (M-2-3)
    (M-2-3) edge node{$(S^1)^m\x \freeO m(\wind(\pi))$} (M-3-4)
    (M-3-4) edge node{$(S^1)^m\x \freeO m((S^1)^n\x \freeO n(h))$} (M-4-5)
;
\end{tikzpicture}
\]
We shall prove that the following diagram commutes in $\Cov$. The span at the top of the diagram is $\twistO{n+m}(\pi,h)$, and the span at the bottom is $\twistO m(\twistO n(\pi,h))$.
\begin{equation}\label{eqIterationDiagram}
\hspace{-.5cm}
\begin{tikzpicture}[ampersand replacement=\&, baseline=(M-2-2).base]
\scriptsize 
\matrix (M) [matrix of math nodes] {
   \&[1cm] \PB{n+m}(\pi)  \&[1cm]    \\[1.5cm]
  (S^1)^{n+m}\x \freeO{n+m}(X)\& \PB m((\ev^n_X)^*\pi) \& (S^1)^{n+m}\x \freeO{n+m}(Z) \\[1.5cm]
 (S^1)^m\x \freeO m((S^1)^n\x\freeO n(X)) \&\&  (S^1)^m\x \freeO m((S^1)^n\x \freeO n(Z))  \\
};
\path [auto, arrow, ->]
    (M-1-2) edge[->>] node[swap]{$(\ev^{n+m}_X)^*\pi$} (M-2-1)
            edge (M-2-2)
            edge node{$(S^1)^{n+m}\x\freeO{n+m}(h)\circ \wind^{n+m}(\pi)$} (M-2-3)
    (M-2-1) edge node[swap]{$\iota_X$} (M-3-1)
    (M-2-3) edge node{$\iota_Z$} (M-3-3)
    (M-2-2) edge[->>]  (M-3-1)
            edge (M-3-3)
;
\end{tikzpicture}
\end{equation}
We shall prove that the diagram commutes in three steps. First, the left hand square of \eqref{eqIterationDiagram}
\begin{equation}\label{eqIterationDiagramLeft}
\begin{tikzpicture}[ampersand replacement=\&, baseline=(M-2-1).base]
\scriptsize 
\matrix (M) [matrix of math nodes] {
   \&[0cm] \PB{n+m}(\pi)  \&[0cm]    \\[1.5cm]
  (S^1)^{n+m}\x \freeO{n+m}(X)\&\& \PB m((\ev^n_X)^*\pi) \\[1.5cm]
 \& (S^1)^m\x \freeO m((S^1)^n\x\freeO n(X)) \&  \\
};
\path [auto, arrow, ->]
    (M-1-2) edge[->>] node[swap]{$(\ev^{n+m}_X)^*\pi$} (M-2-1)
            edge (M-2-3)
    (M-2-1) edge node[swap]{$\iota_X$} (M-3-2)
    (M-2-3) edge[->>] node{$(\ev^m_{(S^1)^n\x\freeO n(X)})^*(\ev^n_X)^*\pi$} (M-3-2)
;
\end{tikzpicture}
\end{equation}
commutes in $\Cov$ if and only if it is a pullback square. The right hand square of \eqref{eqIterationDiagram} decomposes into two parts, and commutativity of these will be step two and three:
\begin{equation}\label{eqIterationDiagramRight}
\begin{tikzpicture}[ampersand replacement=\&, baseline=(BASE).base]
\scriptsize 
\matrix (M) [matrix of math nodes] {
   \PB{n+m}(\pi)  \&[-.5cm]   \&[-.7cm]   \&[-.9cm] \\[1cm]
   \PB m((\ev^n_X)^*\pi) \&\& (S^1)^{n+m}\x\freeO{n+m}(Y) \&\\[1cm]
 \& (S^1)^m\x \freeO m(\PB n(\pi)) \&\& (S^1)^{n+m}\x \freeO{n+m}(Z) \\[1cm]
 \&\& (S^1)^m\x \freeO m((S^1)^n \x \freeO n(Y)) \& \\[1cm]
 \&\&\& (S^1)^m\x \freeO m((S^1)^n\x \freeO n(Z)) \\
};
\path [auto, arrow, ->]
    (M-1-1) edge (M-2-1)
            edge node{$\wind^{n+m}(\pi)$} (M-2-3)
    (M-2-3) edge node{$\iota_Y$} (M-4-3)
            edge node{$(S^1)^{n+m}\x\freeO{n+m}(h)$} (M-3-4)
    (M-3-4) edge node{$\iota_Z$} (M-5-4)
    (M-2-1) edge node[swap](BASE){$\wind^m((\ev^n_X)^*\pi)$} (M-3-2)
    (M-3-2) edge node[swap]{$(S^1)^m\x \freeO m(\wind^n(\pi))$} (M-4-3)
    (M-4-3) edge node[swap]{$(S^1)^m\x \freeO m((S^1)^n\x \freeO n(h))$} (M-5-4)
;
\end{tikzpicture}
\end{equation}
As mentioned the leftmost square \eqref{eqIterationDiagramLeft} commutes in $\Cov$ if and only if is a pullback square. Note that the embedding $\iota_X\colon (S^1)^{n+m}\x\freeO {n+m}(X) \to (S^1)^m\x \freeO m((S^1)^n\x\freeO n(X))$ described in \eqref{eqIterationEmbedding} commutes with the evaluation maps of $\freeO{n+m}$, $\freeO n$, and $\freeO m$ since we have equalities
\begin{align*}
    \\ {}& (\ev^n_X\circ \ev^m_{(S^1)^n\x \freeO n(X)}\circ \iota_X)((\tup s,\tup r),f)
    \\ ={}& (\ev^n_X\circ \ev^m_{(S^1)^n\x \freeO n(X)} ) (\tup r,\tup r'\mapsto (\tup s,\tup s'\mapsto f(s',r')))
    \\ ={}& \ev^n_X  (\tup s,\tup s'\mapsto f(s',r))
    \\ ={}& f(\tup s,\tup r)
    \\ ={}& \ev^{n+m}_X((\tup s,\tup r),f).
\end{align*}
Hence if we pull the covering map $(\ev^m_{(S^1)^n\x\freeO n(X)})^*(\ev^n_X)^*\pi$ back along $\iota_X$, then we get the covering map $(\ev^{n+m}_X)^*\pi$, and thus \eqref{eqIterationDiagramLeft} is a pullback. Consequently, $\iota_X$ followed by the backward map $(\ev^m_{(S^1)^n\x\freeO n(X)})^*(\ev^n_X)^*\pi$ equals the span
\[
\begin{tikzpicture}[ampersand replacement=\&, baseline=(BASE).base]
\matrix (M) [matrix of math nodes] {
  \&[-.5cm] \PB{n+m}(\pi)  \&[-.5cm]   \\[1cm]
  (S^1)^{n+m}\x \freeO{n+m}(X)\&\& \PB m((\ev^n_X)^*\pi) \\
};
\path [auto, arrow, ->]
    (M-1-2) edge[->>] node[swap]{$(\ev^{n+m}_X)^*\pi$} (M-2-1)
            edge (M-2-3)
;
\end{tikzpicture}
\]
 in $\Cov$.

Let us next focus on the first square of \eqref{eqIterationDiagramRight}, and let $(y,(\tup s,\tup r),f)\in\PB{n+m}(\pi)$. Recall that $((\tup s,\tup r),f)\in (S^1)^{n+m}\x \freeO{n+m}(X)$ and $y\in Y$ are required to satisfy $y\in \pi^{-1}(f(\tup s,\tup r))$.
The exponents $k(y,(\tup s,\tup r),f)_i$, $1\leq i\leq n+m$, which are used to define $\wind^{n+m}(\pi)$, are determined by ability to lift (powers of) the loops 
\begin{align*}
\gamma_i &=\begin{cases} f((s_1,\dotsc,s_{i-1},-,s_{i+1},\dotsc,s_n),\tup r) &\text{for $1\leq i\leq n$,}
\\ f(\tup s,(r_1,\dotsc,r_{i-n-1}, - ,r_{i-n+1},\dotsc, r_m)) &\text{for $n+1\leq i\leq m$}\end{cases}
\end{align*}
from $X$ to $Y$ starting at $y\in\pi^{-1}(f(\tup s,\tup r))$.

The map between pullbacks $\PB{n+m}(\pi)\to\PB m((\ev^n_X)^*\pi)$ takes an element $(y,(\tup s,\tup r),f)\in \PB{n+m}(\pi)$ to the element $((y,\tup s,\tup s'\mapsto f(\tup s',\tup r)),\iota_X((\tup s,\tup r),f))\in \PB m((\ev^n_X)^*\pi)$, where $(y,\tup s, \tup s'\mapsto f(\tup s',\tup r))\in \PB n(\pi)$ and $\iota_X((\tup s,\tup r),f)\in (S^1)^m\x \freeO m((S^1)^n\x\freeO n(X))$ satisfy
\[
((\ev^n_X)^*\pi)(y,\tup s, \tup s'\mapsto f(\tup s',\tup r)) = (\tup s,\tup s'\mapsto f(\tup s',\tup r)) = \ev^m_X(\iota_X((\tup s,\tup r),f)).
\]
To determine the value of the map $\wind^m((\ev^n_X)^*\pi)$, we need the exponents 
\[
k((y,\tup s,\tup s'\mapsto f(\tup s',\tup r)),\iota_X((\tup s,\tup r),f))_i
\]
for $1\leq i\leq m$, hence we have to take some power of the following loop $S^1\to (S^1)^n\x\freeO n(X)$
\[(\tup s, \tup s' \mapsto f(\tup s',(r_1,\dotsc,r_{i-1}, - ,r_{i+1},\dotsc, r_m)))\]
and lift it along the covering map $(\ev^n_X)^*\pi$ to get a loop $S^1\to \PB n(\pi)$ starting at $(y,\tup s, \tup s'\mapsto f(\tup s',\tup r))\in \PB n(\pi)$. Lifting a loop along $(\ev^n_X)^*\pi$ is equivalent to lifting the loop composed with $\ev^n_X$ along the covering map $\pi$. Hence we are asking to lift some power of the loop
\begin{multline*}
\ev^n_X(\tup s, \tup s' \mapsto f(\tup s',(r_1,\dotsc,r_{i-1}, - ,r_{i+1},\dotsc, r_m))) 
\\= f(\tup s,(r_1,\dotsc,r_{i-1}, - ,r_{i+1},\dotsc, r_m)) = \gamma_{n+i}
\end{multline*}
along the covering map $\pi$ starting at $y$. Consequently we have
\[k((y,\tup s,\tup s'\mapsto f(\tup s',\tup r)),\iota_X((\tup s,\tup r),f))_i = k(y,(\tup s,\tup r),f)_{n+i}\]
for each $1\leq i\leq m$.

Let $u\colon (S^1)^m\to Y$ be the unique lift of the $m$-fold loop
\[
\tup r'\mapsto f(\tup s,(k(y,(\tup s,\tup r),f)_{n+1}\cdot (\tup r'_1-\tup r_1)+\tup r_1,\dotsc, k(y,(\tup s,\tup r),f)_{n+m}\cdot (\tup r'_m-\tup r_m)+\tup r_m))
\]
from $X$ to $Y$ and starting at $u(\tup r)=y$.
Then $\wind^m((\ev^n_X)^*\pi)$ takes the point 
\[
((y,\tup s,\tup s'\mapsto f(\tup s',\tup r)),\iota_X((\tup s,\tup r),f))\in \PB m((\ev^n_X)^*\pi)
\]
to the following point of $(S^1)^m\x\freeO m(\PB n(\pi))$:
\begin{multline} \label{thisone}
    (\tup r, \tup r' \mapsto (u(\tup r'),\tup s,\tup s'\mapsto 
    \\ f(\tup s',(k(y,(\tup s,\tup r),f)_{n+1}\cdot (\tup r'_1-\tup r_1)+\tup r_1,\dotsc, k(y,(\tup s,\tup r),f)_{n+m}\cdot (\tup r'_m-\tup r_m)+\tup r_m)))).
\end{multline}
Next we need to apply $(S^1)^m\x \freeO m(\wind^n(\pi))$ to the element above. This means we have to postcompose the $m$-fold loop $(S^1)^m\to \PB n(\pi)$  with the map $\wind^n(\pi)\colon \PB n(\pi) \to \PBlift n(\pi)\cong (S^1)^n\x\freeO n(Y)$.
The $m$-fold loop $(S^1)^m\to \PB n(\pi)$ lives in a single component of $\PB n(\pi)$, and thus the value of $k\colon \PB n(\pi)\to (\ZZ_{>0})^n$ is constant along the $m$-fold loop. To find the exponents needed to calculate $\wind^n(\pi)$ along the $m$-fold loop, it is therefore enough to determine the exponents at any single point, say for instance at the point $\tup r\in (S^1)^m$. The $m$-fold loop $(S^1)^m\to \PB n(\pi)$ in \eqref{thisone} is
\[
\tup r' \mapsto (u(\tup r'),\tup s,\tup s'\mapsto f(\tup s',(k(y,(\tup s,\tup r),f)_{n+1}\cdot (\tup r'_1-\tup r_1)+\tup r_1,\dotsc, k(y,(\tup s,\tup r),f)_{n+m}\cdot (\tup r'_m-\tup r_m)+\tup r_m)))
\]
 which evaluated at $\tup r$ simply gives us the following point of $\PB n(\pi)$:
\[(y,\tup s,\tup s'\mapsto f(\tup s',\tup r)).\]
The exponent $k(y,\tup s,\tup s'\mapsto f(\tup s',\tup r))_i$ for $1\leq i\leq n$ is determined by lifting a power of the loop
\[f((s_1,\dotsc,s_{i-1}, - ,s_{i+1},\dotsc, s_n),\tup r) = \gamma_{i}\]
from $X$ to $Y$ starting at $y\in \pi^{-1}(f(\tup s,\tup r))$. We therefore conclude that
\[k(y,\tup s,\tup s'\mapsto f(\tup s',\tup r))_i = k(y,(\tup s,\tup r),f)_{i}\]
for each $1\leq i\leq n$.
If we apply $(S^1)^m\x \freeO m(\wind^n(\pi))$ to our element of $(S^1)^m\x\freeO m(\PB n(\pi))$ from \eqref{thisone}, we then arrive at the following element of $(S^1)^m\x \freeO m(\PBlift n(\pi))$:
\begin{multline}\label{eqIteratedWind}
    (\tup r, \tup r' \mapsto (u(\tup r'),\tup s,\tup s'\mapsto f((k(y,(\tup s,\tup r),f)_1\cdot (\tup s'_1-\tup s_1)+\tup s_1,\dotsc, k(y,(\tup s,\tup r),f)_n\cdot (\tup s'_n-\tup s_n)+\tup s_n),
    \\ (k(y,(\tup s,\tup r),f)_{n+1}\cdot (\tup r'_1-\tup r_1)+\tup r_1, \dotsc, k(y,(\tup s,\tup r),f)_{n+m}\cdot (\tup r'_m-\tup r_m)+\tup r_m)))).
\end{multline}
Let $g\colon (S^1)^{n+m}\to Y$ be the unique lift of 
\begin{multline*}
    (\tup s',\tup r')\mapsto f((k(y,(\tup s,\tup r),f)_1\cdot (\tup s'_1-\tup s_1)+\tup s_1,\dotsc, k(y,(\tup s,\tup r),f)_n\cdot (\tup s'_n-\tup s_n)+\tup s_n),
    \\ (k(y,(\tup s,\tup r),f)_{n+1}\cdot (\tup r'_1-\tup r_1)+\tup r_1, \dotsc, k(y,(\tup s,\tup r),f)_{n+m}\cdot (\tup r'_m-\tup r_m)+\tup r_m))
\end{multline*}
from $X$ to $Y$ starting at $y$, hence in particular $g(\tup s,\tup r')=u(\tup r')$. Then the element of \eqref{eqIteratedWind} in $(S^1)^m\x \freeO m(\PBlift n(\pi))$ is equivalent to the following point of $(S^1)^m\x \freeO m((S^1)^n\x\freeO n(Y))$:
\[(\tup r, \tup r' \mapsto (\tup s,\tup s'\mapsto g(\tup s',\tup r'))) = \iota_Y((\tup s,\tup r),g)
= (\iota_Y\circ \wind^{n+m}(\pi))((\tup s,\tup r),f).
\]
Since $((\tup s,\tup r),f)\in \PB{n+m}(\pi)$ was arbitrary, we conclude that the first square of \eqref{eqIterationDiagramRight} commutes.

Commutativity of the second square in \eqref{eqIterationDiagramRight} is basically stating that $\iota\colon (S^1)^{n+m}\x \freeO{n+m}(-)\into (S^1)^m\x \freeO m((S^1)^n\x\freeO n(-))$ is natural in $\Top$, and as such is significantly easier to confirm. Let $((\tup s,\tup r),f)\in (S^1)^{n+m}\x \freeO{n+m}(Y)$. We then do the straighforward calculation
\begin{align*}
    {}& ((S^1)^m\x \freeO m((S^1)^n\x \freeO n (h)))(\iota_Y((\tup s,\tup r),f))
    \\ ={}& ((S^1)^m\x \freeO m((S^1)^n\x \freeO n (h)))(\tup r, \tup r'\mapsto (\tup s, \tup s'\mapsto f(\tup s',\tup r')))
    \\ ={}& (\tup r, \tup r'\mapsto (\tup s,\tup s'\mapsto (h\circ f)(\tup s',\tup r')))
    \\ ={}& \iota_Z((\tup s,\tup r),h\circ f)
    \\ ={}& \iota_Z(((S^1)^{n+m}\x\freeO{n+m}(h))((\tup s,\tup r),f)).
\end{align*}
The second square of \eqref{eqIterationDiagramRight} commutes, and since the entire diagram \eqref{eqIterationDiagram} thus commutes in $\Cov$, we conclude that $\iota\colon \twistO{n+m}(-) \Rightarrow \twistO m(\twistO n(-))$ is a natural transformation of endofunctors $\Cov\to \Cov$.

\ref{itemCovLnPartialEvaluation}: If we combine \ref{itemCovIterateLn} with \ref{itemCovEvalSquare}, we get a natural transformation 
\[\twistO{n+1}(-)\overset{\iota}\Longrightarrow \twistO 1(\twistO n(-))\overset{\ev}\Longrightarrow \twistO n(-)\]
where the second transformation evaluates $S^1 \x \freeO 1(\twistO n(-))\to \twistO n(-)$. Given any space $X$, the combined map $(S^1)^{n+1}\x \freeO{n+1}(X)\to (S^1)^n\x \freeO n(X)$ takes a pair $((\tup s, r),f)\in (S^1)^{n+1}\x \freeO{n+1}(X)$ to
\begin{align*}
    & \ev_{(S^1)^n\x \freeO n(X)}^1(\iota_X((\tup s,r), f))
    \\ ={}& \ev_{(S^1)^n\x \freeO n(X)}^1(r,r'\mapsto (\tup s, \tup s'\mapsto f(\tup s',r')))
    \\ ={}& (\tup s, \tup s'\mapsto f(\tup s,r)).
\end{align*}
This is precisely the identity map on $(S^1)^n$ times the partial evaluation map $S^1\x\freeO{n+1}(X)\to \freeO n(X)$. Hence $(S^1)^n \x \pev$ coincides with the natural transformation \[\twistO{n+1}(-)\overset{\iota}\Longrightarrow \twistO 1(\twistO n(-))\overset{\ev}\Longrightarrow \twistO n(-).\qedhere\]
\end{proof}

\begin{remark}\label{remarkPassesToHoCov}
The constructions and results of this section pass to the homotopy category of $\Cov$. A homotopy class of maps $X \to Z$ in $\Cov$ is given by a choice of finite cover $Y$ of $X$ (up to isomorphism) and a homotopy class of maps from $Y$ to $Z$. The claim then follows from the fact that $\twistO n(-)$ is equal to $(S^1)^n \times \freeO n(-)$ on forward maps and $(S^1)^n \times \freeO n(-)$ passes to the homotopy category of spaces.
\end{remark}

%% file: groups_and_bisets_I.tex
\section{Bisets and free loop spaces for finite groups}\label{secGroups}
The goal of this section is to apply the results of the previous section to give explicit formulas for $\twistO{n}$ on the full subcategory of $\Ho(\Cov)$ consisting of spaces equivalent to a finite coproduct of classifying spaces of finite groups. This is possible because of the close relationship between this full subcategory of $\Ho(\Cov)$ and the Burnside category of finite groups \cite{MillerGroupoids}. In fact, these formulas for $\twistO{n}$ extend to give an endofunctor of the Burnside category of (formal unions of) finite groups.

If $G$ is a finite group then every connected covering space over the usual model for $BG$ is homeomorphic to $EG \times_G G/R \tto BG$ for some subgroup $R\leq G$ determined up to $G$-conjugation. The total space of this cover is homotopy equivalent to $BR$. Given a third finite group $H$, any map $BR\to BH$ is homotopic to $B(\ph)$ for some group homomorphism $\ph\colon R\to H$ determined up to $R$- and $H$-conjugation.
Thus any span in $\Cov$ from $BG$ to $BH$ is homotopic to a union of connected spans of the form $BG\ffrom EG/R \simeq BR\to BH$. Such spans can be modeled algebraically by $H$-free $(G,H)$-bisets (see Proposition \ref{propEmbedInHoCov}).
Since we wish to model free loop spaces $\freeO nBG$, we shall also consider disjoint unions of classifying spaces and model these by formal unions of finite groups.

\subsection{Burnside modules for finite groups and formal unions of finite groups} \label{sec:burnside}
\begin{definition}\label{defGroupBisetComposition}
Let $G$ and $H$ be finite groups. We let $\AG_+(G,H)$ denote the commutative monoid of isomorphism classes of finite $(G,H)$-bisets with free $H$-action and disjoint union as addition. Furthermore, let $\AG(G,H)$ denote the Grothendieck group of the commutative monoid $\AG_+(G,H)$. We will refer to the elements in $\AG(G,H)$ as virtual bisets and to $\AG(G,H)$ as the Burnside module from $G$ to $H$.

Given a third group $K$,  we have a composition map
\begin{alignat*}{2}
\circ\colon \AG(H,K) &\times \AG(G,H) &\ \rightarrow\ & \AG(G,K),
\\ Y &\ ,\  X &\ \mapsto\ & X\x_H Y,
\end{alignat*}
both on the level of commutative monoids and the Grothendieck groups.

The composition map is bilinear. The composition map necessarily switches the order of $X$ and $Y$ so that the acting groups match up. In an attempt to lessen the confusion when we compose many bisets at once, we shall instead exclusively use ``right-composition'' for bisets:
\begin{alignat*}{2}
\cmp\colon \AG(G,H) &\times \AG(H,K) &\ \rightarrow\ & \AG(G,K),
\\ X &\ ,\  Y &\ \mapsto\ & X\x_H Y.
\end{alignat*}
\end{definition}

If we think of $\AG(G,H)$ as denoting morphisms from $G$ to $H$ (in the biset category $\AG$ defined below), ``right-composition'' is the commuting diagram
\[
\begin{tikzpicture}
  \node (M) [matrix of math nodes] {
    G &[1cm] H \\[1cm]
     & K. \\
  };
  \path[->,arrow,auto]
    (M-1-1) edge node{$X$} (M-1-2)
            edge node[swap]{$X\cmp Y$} (M-2-2)
    (M-1-2) edge node{$Y$} (M-2-2)
  ;
\end{tikzpicture}
\]

The commutative monoid $\AG_+(G,H)$ is free with generators given by the isomorphism classes of transitive $H$-free $(G,H)$-bisets -- these also form a canonical basis of $\AG(G,H)$ as a $\Z$-module. The transitive $H$-free $(G,H)$-bisets are of the form
\[
G \times_{R}^{\ph} H = (G \times H) \big / (gr,h) \sim (g,\ph(r)h),
\]
where $R \leq G$ is a subgroup of $G$ (taken up to conjugacy in $G$) and $\ph \colon R\rightarrow H$ is a group homomorphism (taken up to conjugacy in $G$ and $H$). We will denote these $(G,H)$-bisets by $[R,\ph]_{G}^{H}$ or just $[R,\ph]$ when $G$ and $H$ are clear from context. It is also common to denote a virtual biset $X \in \AG(G,H)$ as $X_G^H$ when $G$ and $H$ are not clear from context.

\begin{convention}\label{conventionOrbitDecomposition}
Given $X\in \AG(G,H)$, we can write $X$ as a linear combination of transitive bisets:
\[\sum_{(R,\ph)} c_{R,\ph}\cdot [R,\ph]_G^H.\]
The summation runs over all $R\leq G$ and $\ph\colon R\to H$ (not taken up to conjugacy). The coefficient function $c_{(-)}$ is a choice of function from the set of all pairs $(R,\ph)$ to $\Z$ such that the sum of coefficients $c_{R',\ph'}$ over the pairs $(G,H)$-conjugate to $(R,\ph)$ is the number of copies of the orbit $[R,\ph]_G^H$ in $X$.

In particular the linear combination above is \emph{not} unique and we allow isomorphic biset orbits to be part of the sum for several different but conjugate pairs $(R,\ph)$. If we require $c_{(-)}$ to be concentrated on chosen representatives for the conjugacy classes of pairs, then the linear combination is unique.
\end{convention}

We will now extend the notion of Burnside module to a formal union of finite groups. A formal union of finite groups is an ordered finite tuple of finite groups (repetition is permitted).

\begin{definition}
  Let $G$ and $H$ be formal unions of finite groups. Suppose $G= G_1 \sqcup \dotsb \sqcup G_n$ and $H=H_1\sqcup \dotsb \sqcup H_m$ are the decompositions into groups -- the components of $G$ and $H$.

  We define $\AG(G,H)$ to be the set of matrices of virtual bisets between the components, where the rows correspond to the components of $G$ and the columns correspond to the components of $H$. Hence $X\in \AG(G,H)$ is an $(n\x m)$-matrix with entries
  \[X_{i,j} \in \AG(G_i,H_j)\quad\text{for $1\leq i\leq n$, $1\leq j\leq m$.}\]
  The composition $X\cmp Y$ for $X\in \AG(G,H)$ and $Y\in \AG(H,K)$ is just the usual composition of matrices:
  \[(X\cmp Y)_{i,k} = \sum_{j=1}^{m} X_{i,j}\cmp Y_{j,k}.\]
  Similarly, we let $\AG_+(G,H)\subseteq\AG(G,H)$ denote the subset of matrices $X$ in which every matrix entry is a biset $X_{i,j}\in\AG_+(G_i,H_j)$.
\end{definition}

It follows from the definition above that changing the order in a formal union of finite groups results in a canonically isomorphic finite union of finite groups. 

\begin{definition}
We let $\AG$ denote the category in which the objects are formal unions of finite groups and, for any pair of objects $G,H$, the set of morphisms from $G$ to $H$ is the Burnside module $\AG(G,H)$ of virtual biset matrices.

Similarly, the category $\AG_+$ has morphism sets $\AG_+(G,H)$.
\end{definition}

\begin{remark}
The category $\Grp$ of finite groups and homomorphisms maps to $\AG_+$ (and $\AG$) by sending a group homomorphism $\ph\colon G\to H$ to the transitive $(G,H)$-biset $[G,\ph]_G^H\in \AG_+(G,H)$, which is isomorphic to the biset $H_{G,\ph}^H$, where $G$ acts on $H$ from the left through $\ph$.
\end{remark}

The category of finite groupoids has a well-behaved notion of a finite covering map \cite[Definition 1.6]{MillerGroupoids}. The classifying space functor takes a finite cover of finite groupoids to a finite cover of spaces. It is convenient to use the homotopy category of the category of spans of finite groupoids in which the wrong way map is a finite covering map as an intermediate category between $\AG_+$ and $\Ho(\Cov)$. Any finite groupoid is equivalent to a coproduct of finite groups, where we view a group as being a groupoid with a single object. 

Given a finite group $G$ and a $G$-action on a set $X$, we will write $X \mmod G$ for the associated action groupoid. We will use $B$ to denote the realization of a groupoid as well as the map $\AG_+\to \Ho(\Cov)$ defined below.

\begin{definition}\label{defBisetsToCov}
We define a functor $B\colon \AG_+\to \Ho(\Cov)$ as follow. For a formal union of groups $G=G_1\sqcup\dotsb\sqcup G_n$, we let $BG = BG_1\sqcup\dotsb\sqcup BG_n$, which is the realization of the groupoid $(\ast\mmod G_1)\sqcup\dotsb\sqcup (\ast\mmod G_n)$.

Suppose $G$ and $H$ are finite groups and $X\in \AG_+(G,H)$ is a biset. The group $G\x H$ acts on $X$ from the right by $x(g,h) = g^{-1}xh$. Because $X$ is $H$-free, we have an equivalence of groupoids $X\mmod (G\x H)\simeq (X/H)\mmod G$. The equivalence depends on a choice of representatives for the $H$-orbits of $X$, but different choices result in naturally isomorphic equivalences. Hence we have a canonical well-defined equivalence $B(X\mmod (G\x H))\simeq B((X/H)\mmod G)$ in $\Ho(\Top)$.

The functor $B\colon \AG_+\to \Ho(\Cov)$ takes the biset $X\in \AG_+(G,H)$ to the realization of the span of groupoids
\[
\begin{tikzpicture}
\node (M) [matrix of math nodes] {
 &[1cm] (X/H)\mmod G &[1cm]  X \mmod (G\x H) &[1cm] \\[1cm]
 \ast\mmod G &&& \ast\mmod H. \\
};
\path [auto, ->, arrow]
    (M-2-1) edge[<<-] (M-1-2)
    (M-1-2) edge node{$\simeq$} (M-1-3)
    (M-1-3) edge (M-2-4)
;
\end{tikzpicture}
\]
The map of groupoids $(X/H)\mmod G\tto \ast\mmod G$ realizes to a finite sheeted covering map, so we get a well-defined span from $BG$ to $BH$ in $\Ho(\Cov)$.

For formal unions of groups $G=G_1\sqcup\dotsb\sqcup G_n$ and $H=H_1\sqcup\dotsb\sqcup H_m$ and a matrix $X\in\AG_+(G,H)$ of bisets, all the spans 
\[\ast\mmod G_i \ffrom (X_{i,j}/H_j)\mmod G_j \simeq X_{i,j} \mmod (G_i\x H_j) \to \ast\mmod H_j\]
combine to give a span from $BG$ to $BH$.

This construction is closely related to the functor defined in \cite[Section 5]{MillerGroupoids}.
\end{definition}

\begin{definition}
The functor $B\colon \AG_+\to \Ho(\Cov)$ takes a transitive $(G,H)$-biset $[R,\ph]_G^H$ and gives the span $B\ph\circ \tr_R^G\colon BG\to BH$, i.e. the span 
\[
\begin{tikzpicture}
\node (M) [matrix of math nodes] {
 &[1cm] EG/R &[1cm]  BG &[1cm] \\[1cm]
 BG &&& BH \\
};
\path [auto, ->, arrow]
    (M-2-1) edge[<<-] (M-1-2)
    (M-1-2) edge node{$\simeq$} (M-1-3)
    (M-1-3) edge (M-2-4)
;
\end{tikzpicture}
\]
in $\Ho(\Cov)$. Since there is a canonical isomorphism in $\Ho(\Cov)$ between $EG/R$ and $BG$, we will abuse notation and write $\tr_R^G$ for either the span $BG \ffrom EG/R \xto{\simeq} BR$ or for the wrong-way map $BG\ffrom EG/R$ by itself. 
\end{definition}

\begin{prop}\label{propEmbedInHoCov}
The functor $B\colon \AG_+\to \Ho(\Cov)$ is fully faithful and embeds $\AG_+$ as the full subcategory of $\Ho(\Cov)$ spanned by disjoint unions of classifying spaces for finite groups. The addition of bisets in $\AG_+(G,H)$ correspond to the disjoint union of spans between $BG$ and $BH$.
\end{prop}

\begin{proof}
This is the homotopy category version of \cite{MillerGroupoids}*{Theorem 5.2}.
\end{proof}

\subsection{Free loop spaces for classifying spaces of finite groups}\label{subsecFreeLoopGroups}
Since $S^1\simeq B\Z$, the $n$-fold free loop space $\freeO n BG$ is equivalent to $\Map(B\Z^n,BG)$ for any formal union of finite groups $G$. Every homomorphism $\Z^n\to G$ factors through $(\Z/\ell)^n$ for large enough $\ell$, so we have
\[\freeO n BG \simeq \colim_{\ell\in(\Z_{>0},\text{divisibility})} \Map(B(\Z/\ell)^n, BG)\]
for any finite group $G$. From now on we replace $S^1$ with the classifying space $B(\Z/\ell)$ for any $\ell$ sufficiently large (e.g. $\ell$ divisible by the order of $G$).

A homomorphism $(\Z/\ell)^n\to G$ picks out a commuting $n$-tuple in $G$, in the case of a formal union of groups, each commuting $n$-tuple lies in a single component of $G$. The free loop space $\freeO n BG$ is modeled algebraically by the mapping groupoid $\Map(\ast\mmod (\Z/\ell)^n, \ast\mmod G)$ of functors and natural transformations. This mapping groupoid is isomorphic as groupoids to the coproduct of action groupoids
\[\bigsqcup_i (\text{commuting $n$-tuples in $G_i$})\mmod G_i,\]
where each component $G_i$ of $G$ acts on the commuting $n$-tuples in $G_i$ by conjugation.

\begin{definition}
For a finite group $G$, let $\ntuples nG$ denote the set of commuting $n$-tuples of elements in $G$ and let $\cntuples nG$ denote the collection of $G$-conjugacy classes of commuting $n$-tuples in $G$, where we consider $\ntuples 0G$ to consist of the unique empty/trivial $0$-tuple. For a formal union $G=G_1\sqcup\dotsb\sqcup G_n$, we will write $\ntuples nG$ for the disjoint union of the finite sets $\ntuples nG_{i}$ and $\cntuples nG$ for the disjoint union of the finite sets $\cntuples nG_{i}$.
\end{definition}

Abusing notation, we will occasionally write $\tup a\in G$ when $\tup a$ is a tuple of elements in $G$.

\begin{definition}
For an $n$-tuple $\tup a$ in a finite group $G$, the centralizer $C_G(\tup a)$ is defined as
\begin{align*}
C_G(\tup a) :={}& \{g\in G\mid g^{-1}a_ig = a_i \text{ for $1\leq i\leq n$}\}
\\ ={}& \bigcap_{i=1}^n C_G(a_i)=C_G(\gen{a_1,\dotsc,a_n}).
\end{align*}
When $G$ is a formal union, and $\tup a$ is an $n$-tuple in $G$, we define the centralizer $C_G(\tup a)$ to be the centralizer inside the component of $G$ containing $\tup a$. As such, the centralizer $C_G(\tup a)$ is always a finite group and not a formal union.
\end{definition}

The component of $\freeO nBG$ corresponding to a particular commuting $n$-tuple $\tup a$ is homotopy equivalent to $BC_G(\tup a)$.
\begin{definition}\label{defGroupLoop}
For ease of notation, we will write $\freeO n G$ for the formal union
\[\freeO n G := \coprod_{[\tup a]\in \cntuples nG} C_G(\tup a),\]
for some choice of representatives for the conjugacy classes of commuting $n$-tuples in $G$. 

Note that if $G$ is a formal union, $G=G_1\sqcup\dotsb \sqcup G_m$, then $\freeO n G$ decomposes componentwise $\freeO n G = \freeO n (G_1) \sqcup \dotsb \sqcup \freeO n (G_m)$.
\end{definition}

For a finite group $G$, if we view the centralizers $C_G(\tup a)$ as categories with single objects, then the formal union $\freeO n G$ is equivalent as groupoids to 
\[G^{(n)}\mmod G.\]
Hence the formal union $\freeO n G$ models the $n$-fold free loop space of $BG$ up to homotopy:
\[B(\freeO n G) = \coprod_{[\tup a]\in \cntuples nG} BC_G(\tup a)\simeq \freeO n(BG).\]
There is a possible ambiguity in how the algebraic model $\freeO n G$ depends on the choice of representatives for the conjugacy classes of commuting $n$-tuples. However, in $\AG$ (and $\Ho(\Top)$) there are canonical equivalences between any two choices of representatives:
\begin{lemma}\label{lemmaGroupCentralizerIso}
Let $\tup a$ be an $n$-tuple of commuting elements in $G$,  and suppose that $\reptup a$ is another $n$-tuple conjugate to $\tup a$.
Any element $g\in G$ that conjugates $\tup a$ to $\reptup a$ also induces $c_g\colon C_G(\tup a) \to C_G(\reptup a)$, and any such conjugation gives rise to a biset
\[\change_{\tup a}^{\reptup a} \in \AG(C_G(\tup a), C_G(\reptup a)),\]
which is independent of the choice of conjugating element.

If $\reptup a'$ is a third $n$-tuple conjugate to $\tup a$ (and therefore $\reptup a$), then
 the chosen bisets are compatible with composition:
\[\change_{\tup a}^{\reptup a'} = \change_{\tup a}^{\reptup a} \cmp \change_{\reptup a}^{\reptup a'}.\]
\end{lemma}

\begin{proof}
Given any two elements $g,h\in G$ with $g^{-1}\tup a g=h^{-1} \tup a h =\reptup a$, we have $g^{-1} h\in C_G(\reptup a)$. Hence the two conjugation maps $c_g, c_{h}\colon C_G(\tup a)\to C_G(\reptup a)$ are themselves conjugate as maps into $C_G(\reptup a)$, and as such they give rise to the same biset
\[\change_{\tup a}^{\reptup a} = [C_G(\tup a), c_g]_{C_G(\tup a)}^{C_G(\reptup a)} = [C_G(\tup a),c_h]_{C_G(\tup a)}^{C_G(\reptup a)}.\]
Suppose $g'\in G$ satisfies $(g')^{-1}\reptup a g' = \reptup a'$. Then $gg'$ conjugates $\tup a$ to $\reptup a'$ and we have
\[\change_{\tup a}^{\reptup a} \cmp \change_{\reptup a}^{\reptup a'} = [C_G(\tup a), c_g]_{C_G(\tup a)}^{C_G(\reptup a)}\cmp [C_G(\reptup a), c_{g'}]_{C_G(\reptup a)}^{C_G(\reptup a')} = [C_G(\tup a), c_{gg'}]_{C_G(\tup a)}^{C_G(\reptup a')} = \change_{\tup a}^{\reptup a'}\]
as required.
\end{proof}
In view of Lemma \ref{lemmaGroupCentralizerIso}, given two $G$-conjugate $n$-tuples $\tup a,\reptup a$ of commuting elements in a finite group $G$ the diagram of groupoids
\[
\begin{tikzpicture}
\node (Loops) {$(\text{commuting $n$-tuples in $G$})\mmod G$};
\node (A) [above left=of Loops] {$C_G(\tup a)$};
\node (B) [below left=of Loops] {$C_G(\reptup a)$};
\path [auto,->,arrow]
    (A) edge (Loops)
        edge node{$c_g$} (B)
    (B) edge (Loops)
;
\end{tikzpicture}
\]
commutes in $\Ho(\Cov)$ and $\Ho(\Top)$. Hence we have a well-defined equivalence
\begin{lemma}\label{lemmaGroupsAlgModelForFreeLoops}
For every formal union of groups $G$, and $n\geq 0$, there is a canonical equivalence 
\[B(\freeO n G)= \coprod_{[\tup a]\in \cntuples nG} BC_G(\tup a)\xto\simeq \freeO n(BG)\]
in $\Ho(\Cov)$ and $\Ho(\Top)$.
\end{lemma}

\begin{convention}\label{conventionGroupTupleReps}
While Lemma \ref{lemmaGroupCentralizerIso} states that different choices of representatives for the conjugacy classes of tuples are equivalent, it will be necessary later on to fix the choice of representatives once and for all and to explicitly include the isomorphisms of Lemma \ref{lemmaGroupCentralizerIso} in formulas. When necessary, we will therefore assume that a choice of representatives has been made for all $G$-conjugacy classes of commuting $n$-tuples for all $n\geq 0$. 

We will furthermore assume that each representative $(n+1)$-tuple $(a_1,\dotsc,a_{n+1})$ is chosen such that the first $n$ elements $(a_1,\dotsc,a_n)$ form one of the previously chosen representative $n$-tuples.

The convention will be particularly relevant when we apply homomorphisms to $n$-tuples: if we apply $\ph \colon G \to H$ to a representative $\tup a$, we cannot assume that $\ph(\tup a)$ will always be a representative in the codomain as well. Similarly, if we consider the action of $\Sigma_n$ on $\cntuples n G$ permuting the coordinates of the tuple, we cannot ensure that all permutations of a representative $\tup a$ are also representatives -- e.g. $\tup a$ and a permutation $\sigma(\tup a)$ might be conjugate in $G$ but different, so we cannot choose both to be representatives of the conjugacy class.
\end{convention}

Next, the endofunctor $\freeO n\colon \Cov \to \Cov$ has a simple description in the realm of classifying spaces and bisets as follows. The functor $\freeO n$ applied to a union of classifying spaces is again a union of classifying spaces up to homotopy. The equivalence of $\AG_+$ with the full subcategory of $\Ho(\Cov)$ spanned by disjoint unions of classifying spaces therefore makes $\freeO n$ an endofunctor on $\AG_+$ as well since $\freeO n$ is an endofunctor on $\Ho(\Cov)$ by Remark \ref{remarkPassesToHoCov}.
The endofunctor $\freeO n$ on $\Cov$ preserves the monoid structure on $\Cov(X,Y)$, i.e. disjoint union of spans. Hence $\freeO n\colon \AG_+\to \AG_+$ extends to an endofunctor $\freeO n\colon \AG\to \AG$ by linearity on virtual bisets.

Suppose $G$ and $H$ are formal unions of finite groups and suppose $M\in \AG_+(G,H)$ is a matrix of bisets. By the embedding of categories in Proposition \ref{propEmbedInHoCov}, $M$ induces a span $BM\in[BG,BH]$ in $\Ho(\Cov)$. Applying $\freeO n\colon \Cov\to \Cov$ produces a span between $n$-fold free loop spaces $\freeO n(BM)\in [\freeO n BG, \freeO n BH]$. By Lemma \ref{lemmaGroupsAlgModelForFreeLoops} and Proposition \ref{propEmbedInHoCov}, we have canonical isomorphisms of morphism sets
\[\freeO n(BM)\in [\freeO n BG, \freeO n BH]\cong [B(\freeO n G),B(\freeO n H)] \cong \AG_+(\freeO n G,\freeO n H).\]
The resulting matrix in $\AG_+(\freeO n G,\freeO n H)$ is the one we denote by $\freeO n M$.

Given two formal unions of groups $G$ and $H$, the functor $\freeO n\colon \AG\to \AG$ acts componentwise on $G$ and $H$. For every biset matrix $M\in \AG(G,H)$ the resulting matrix $\freeO n(M)\in \AG(\freeO n G,\freeO nH)$ is therefore a block matrix with a block corresponding to each pair $\freeO n(G_i)$ and $\freeO n(H_j)$ and the block equals $\freeO n(M_{i,j})\in  \AG(\freeO n (G_i),\freeO n(H_j))$.
To describe $\freeO n\colon \AG\to \AG$ it is therefore sufficient to consider groups (rather than formal unions of groups).
\begin{prop}\label{propSimpleGroupLoop}
Let $G$ and $H$ be finite groups and suppose $M\in \AG_+(G,H)$ is a $(G,H)$-biset. The matrix $\freeO n(M)\in\AG( \freeO n G , \freeO n H)$ has entries parametrised by our chosen representatives of conjugacy classes of commuting $n$-tuples $\tup a$ in $G$ and commuting $n$-tuples $\tup b$ in $H$, and the entry $\freeO n(M)_{\tup a,\tup b}\in \AG(C_G(\tup a),C_H(\tup b))$ has the form
\[\freeO n(M)_{\tup a,\tup b} = \lc {\tup a} M ^{\tup b} = \{m\in M\mid a_i m = m b_i \text{ for all } 1\leq i\leq n\}.\]
The centralizer $C_G(\tup a)$ acts on the left of the fixed point set $\lc {\tup a} M ^{\tup b}$ and the centralizer $C_H(\tup b)$ acts on the right. The formula for $\freeO n(M)$ extends linearly to all virtual bisets in $\AG(G,H)$.
\end{prop}
\begin{remark}
In the terminology of \cite{GelvinReeh}, we have $\lc {\tup a} M ^{\tup b}=N_M(\tup a,\tup b)$. That is, we can think of $\lc {\tup a} M ^{\tup b}$ as consisting of the elements $m\in M$ ``conjugating'' $m^{-1} \tup a m = \tup b$. The inverse $m^{-1}$ doesn't make sense, so we require $\tup a m = m \tup b$ instead.
\end{remark}

\begin{proof}
By Definition \ref{defBisetsToCov} the span $BM$ from $BG$ to $BH$ in $\Ho(\Cov)$ takes the following form on the level of groupoids:
\[\ast\mmod G \ffrom (M/H)\mmod G \simeq M \mmod (G\x H) \to \ast\mmod H.\]
Applying $\freeO n\colon \Ho(\Cov)\to \Ho(\Cov)$ corresponds to applying $\Map(\ast\mmod(\Z/\ell)^n,-)$ on the level of groupoids to give us
\begin{multline*}
\Map(\ast\mmod(\Z/\ell)^n, \ast\mmod G) \ffrom \Map(\ast\mmod(\Z/\ell)^n,(M/H)\mmod G)\\ \simeq  \Map(\ast\mmod(\Z/\ell)^n,M\mmod (G\x H)) \to \Map(\ast\mmod(\Z/\ell)^n,\ast\mmod H).
\end{multline*}
Each functor from $\ast\mmod(\Z/\ell)^n$ into one of the other groupoids picks out an object and an $n$-tuple of commuting automorphisms of that object. A natural transformation between functors can then be described as acting on the objects by any group element, and conjugating the $n$-tuple by the same element. The mapping groupoids in the span above are therefore isomorphic to the action groupoids below:
\begin{align*}
    &(\text{commuting $n$-tuples in $G$})\mmod G
    \\ \ffrom{}& \left(\begin{tabular}{c}pairs of a comm. $n$-tuple $\tup a$ in $G$\\ and an $x\in M/H$ fixed by $\tup a$\end{tabular}\right) \mmod[\Big] G
    \\ \simeq{}&\left(\begin{tabular}{c}comm. $n$-tuple $\tup a$ in $G$,\\ comm. $n$-tuple $\tup b$ in $H$,\\ and $x\in \lc {\tup a} M ^{\tup b}$\end{tabular}\right)\mmod[\Big] (G\x H)
    \\ \to {}& (\text{commuting $n$-tuples in $H$})\mmod H.
\end{align*}
Let $L\in\AG_+(\freeO n G, \freeO n H)$ be the matrix of bisets with entries
\[L_{\tup a,\tup b} = \lc {\tup a} M ^{\tup b}\in \AG_+(C_G(\tup a),C_H(\tup b)),\]
where $\tup a$ and $\tup b$ run through our chosen representatives of conjugacy classes of commuting $n$-tuples in $G$ and $H$. The span $BL$ from $B\freeO n G$ to $B\freeO n H$ then has the following form on the level of groupoids:
\begin{multline*}
    \coprod_{[\tup a]\in \cntuples n G} \ast \mmod C_G(\tup a)
     \ffrom \coprod_{[\tup a]\in \cntuples n G} \Bigl(\coprod_{[\tup b]\in \cntuples n H} \lc {\tup a} M ^{\tup b} / C_H(\tup b)\Bigr) \mmod[\Big] C_G(\tup a)
    \\ \simeq \coprod_{\substack{[\tup a]\in \cntuples n G\\ [\tup b]\in\cntuples n H}}\lc {\tup a} M ^{\tup b}\mmod (C_G(\tup a)\x C_H(\tup b))
    \to  \coprod_{[\tup a]\in \cntuples n G} \ast\mmod C_H(\tup b).
\end{multline*}
The map of $C_G(\tup a)$-sets $\coprod_{[b]\in\cntuples n H} \lc {\tup a} M ^{\tup b} / C_H(\tup b)\to \lc{\tup a} M/H$ is a bijection because, for each coset $xH\in \lc{\tup a} M/H$ being fixed by $\tup a$, there is a unique tuple $\tup b$ with $\tup a x=x\tup b$, and changing the representative of the coset $xH$ changes $\tup b$ up to $H$-conjugation.

We now have a commutative diagram of groupoids where each row is an equivalence:
\[
\begin{tikzpicture}[ampersand replacement=\&, baseline=(BASE).base]
\node (M) [matrix of math nodes, nodes={text height=16pt, text depth=14pt}] {
    \displaystyle\coprod_{[\tup a]\in \cntuples n G} \ast \mmod C_G(\tup a) \&[1cm] (\text{comm. $n$-tuples in $G$})\mmod G \\[1cm]
    \displaystyle\coprod_{[\tup a]\in \cntuples n G} \Bigl(\coprod_{[\tup b]\in \cntuples n H} \lc {\tup a} M ^{\tup b} / C_H(\tup b)\Bigr) \mmod[\Big] C_G(\tup a) \& \left(\begin{tabular}{c}comm. $n$-tuple $\tup a$ in $G$,\\ and $x\in M/H$ fixed by $\tup a$\end{tabular}\right) \mmod[\Big] G \\[1cm]
    \displaystyle\coprod_{\substack{[\tup a]\in \cntuples n G\\ [\tup b]\in\cntuples n H}}\lc {\tup a} M ^{\tup b}\mmod (C_G(\tup a)\x C_H(\tup b)) \&  \left(\begin{tabular}{c}comm. $n$-tuple $\tup a$ in $G$,\\ comm. $n$-tuple $\tup b$ in $H$,\\ and $x\in \lc {\tup a} M ^{\tup b}$\end{tabular}\right)\mmod[\Big] (G\x H) \\[1cm]
    \displaystyle\coprod_{[\tup a]\in \cntuples n G} \ast\mmod C_H(\tup b) \& (\text{comm. $n$-tuples in $H$})\mmod H.\\
};
\path [auto, arrow, ->] 
    (M-1-1) edge node{$\simeq$} (M-1-2)
    (M-2-1) edge[->>] (M-1-1)
            edge node{$\simeq$} (M-3-1)
            edge node{$\simeq$} (M-2-2)
    (M-3-1) edge node{$\simeq$} (M-3-2)
            edge (M-4-1)
    (M-4-1) edge node{$\simeq$} (M-4-2)
    (M-2-2) edge[->>] (M-1-2)
            edge node{$\simeq$} (M-3-2)
    (M-3-2) edge (M-4-2)
;
\end{tikzpicture}
\]
The level-wise equivalences of groupoids show that the span $BL$ (the left column) from $B\freeO n G$ to $B\freeO n H$ corresponds to the span $\freeO n BM$ (the right column) from $\freeO n BG$ to $\freeO n BH$ under the bijection of morphism sets
\[[B\freeO n G, B\freeO n H] \cong [\freeO n BG, \freeO n BH]\]
coming from the equivalence in Lemma \ref{lemmaGroupsAlgModelForFreeLoops}.
Finally, the span $BL$ (and hence $\freeO n BM$) corresponds to the matrix of bisets $L\in \AG_+(\freeO nG, \freeO n H)$ as required.
\end{proof}

\begin{remark}\label{remarkSimpleLoopFunctor}
We note two special cases of the description of $\freeO n(M)$ above:
First consider the case in which $M$ comes from an actual group homomorphism $\ph\colon G\to H$, i.e. $M=[G,\ph]_G^H=\prescript{}{\ph}H_{\id}$. The set of fixed points $\lc {\tup a} M^{\tup b}$ is empty unless $\tup b$ is $H$-conjugate to $\ph(\tup a)$. In this case, choose some $h_0\in H$ that conjugates $\ph(\tup a)$ to $\tup b$. We then have $\change_{\ph(\tup a)}^{\tup b}=[C_H(\ph(\tup a)),c_{h_0}]\in \AG(C_H(\ph(\tup a)),C_H(\tup b))$ by the isomorphism of Lemma \ref{lemmaGroupCentralizerIso}, and we get
\begin{multline*}
\lc {\tup a} M^{\tup b} = \lc{\tup a} M^{\ph(\tup a)}\cmp \change_{\ph(\tup a)}^{\tup b}=\lc {\tup a} (\prescript{}{\ph}H_{\id}) ^{\ph(\tup a)}\cmp \change_{\ph(\tup a)}^{\tup b}= \{h\in H\mid \ph(a_i) h = h \ph(a_i)\}\cmp \change_{\ph(\tup a)}^{\tup b}\\ = \prescript{}{\ph}{}\bigl(C_H(\ph(\tup a))\bigr)_{\id}\cmp \change_{\ph(\tup a)}^{\tup b} = [C_G(\tup a),\ph]_{C_G(\tup a)}^{C_H(\ph(\tup a))}\cmp \change_{\ph(\tup a)}^{\tup b} = [C_G(\tup a),c_{h_0}\circ \ph]_{C_G(\tup a)}^{C_H(\tup b)}.
\end{multline*}
Thus we simply get the biset version of the homomorphism $\ph\colon C_G(\tup a) \to C_H(\ph(\tup a))$ restricted from $\ph\colon G\to H$, followed by mapping $\ph(\tup a)$ to the representative, $\tup b$, of its conjugacy class.

Next consider the case of a transfer map $\tr^G_H$ from $G$ to a subgroup $H$; this is the case $M=[H,\id]_G^H = {\prescript{}{G}G}_H$. The fixed point set $\lc {\tup a} ({\prescript{}{G}G}_{H})^{\tup b}$ is the collection of $g\in G$ that conjugate $\tup a$ to $\tup b$, i.e. $g^{-1}a_ig = b_i$ for $1\leq i\leq n$. Any two $g$'s that conjugate $\tup a$ to $\tup b$ differ by an element in the centralizer $C_G(\tup a)$, so $\lc {\tup a} ({\prescript{}{G}G}_{H})^{\tup b}$ is a transitive $(C_G(\tup a), C_H(\tup b))$-biset. Choosing a representative $g_0\in \lc {\tup a} ({\prescript{}{G}G}_{H})^{\tup b}$, we have
\[\lc {\tup a} ({\prescript{}{G}G}_{H})^{\tup b} = [C_G(\tup a)\cap g_0Hg_{0}^{-1}, c_{g_{0}}]_{C_G(\tup a)}^{C_H(\tup b)}.\]
This is the transfer from $C_G(\tup a)$ to $C_G(\tup a)\cap g_{0}Hg_{0}^{-1}=C_{g_0 H g_0^{-1}}(\tup a)$ followed by conjugation onto $C_H(\tup b)$.
\end{remark}

Both of the examples in Remark \ref{remarkSimpleLoopFunctor} are special cases of the following corollary that expresses the orbit decomposition of $\freeO n(M)_{\tup a,\tup b}$ in general.
\begin{cor}\label{corSimpleGroupLoopOrbits}
Let $G$ and $H$ be finite groups and suppose $M\in \AG(G,H)$ is a virtual $(G,H)$-biset. Furthermore, let $\tup a$ in $G$ and  $\tup b$ in $H$ be chosen representatives for conjugacy classes of commuting $n$-tuples. Consider the restriction of $M$ to the centralizer of $\tup a$, and write $M_{C_G(\tup a)}^H$ as a linear combination of basis elements (recalling Convention \ref{conventionOrbitDecomposition}):
\[
M_{C_G(\tup a)}^H =\sum_{(R,\ph)}c_{R,\ph} \cdot [R,\ph]_{C_G(\tup a)}^H,
\]
where $R\leq C_G(\tup a)$ and $\ph\colon R\to H$.

The matrix entry $\freeO n(M)_{\tup a,\tup b}$ then satisfies the formula
\[
\freeO n(M)_{\tup a,\tup b} = \sum_{\substack{(R,\ph)\text{ s.t. $\tup a\in R$ and}\\\text{$\ph(\tup a)$ is $H$-conjugate to $\tup b$}}}c_{R,\ph}\cdot [R,\ph]_{C_G(\tup a)}^{C_H(\ph(\tup a))}\cmp \change_{\ph(\tup a)}^{\tup b},
\]
with $R$, $\ph$, and $c_{R,\ph}$ as in the linear combination above.
\end{cor}

\begin{proof}
By Proposition \ref{propSimpleGroupLoop} we have $\freeO n(M)_{\tup a,\tup b} = \lc{\tup a} M^{\tup b}$, and since $\tup a\in C_G(\tup a)$, the restriction $M_{C_G(\tup a)}^H$ has the same set of fixed points $\lc{\tup a} M^{\tup b}=\lc{\tup a}(M_{C_G(\tup a)}^H)^{\tup b}$. 

We then study the fixed points for each of the basis elements $[R,\ph]_{C_G(\tup a)}^H$ in the decomposition of $M_{C_G(\tup a)}^H$. An element of $\lc{\tup a}([R,\ph]_{C_G(\tup a)}^H)^{\tup b}$ is represented by a pair $(g,h)\in C_G(\tup a)\x H$ such that $g^{-1}\tup a g\in R$ and $h^{-1}\ph(g^{-1} \tup a g)h = \tup b$.
The element $g$ centralizes the $n$-tuple $\tup a$, so we just require $\tup a\in R$ and $h^{-1}\ph(\tup a) h= \tup b$. We can change $g$ to be any element of $C_G(\tup a)$, and $h$ can vary by any element in $C_H(\ph(\tup a))$ on the left (or by any element of $C_H(\tup b)$ on the right) and the requirements will still be satisfied.

Consequently, the set of fixed points $\lc{\tup a}([R,\ph]_{C_G(\tup a)}^H)^{\tup b}$ is empty unless $\tup a\in R$ and $\ph(\tup a)$ is $H$-conjugate to $\tup b$, in which case we get a single $(C_G(\tup a),C_H(\tup b))$-orbit of the form
\[[R,c_h\circ \ph]_{C_G(\tup a)}^{C_H(\tup b)} = [R,\ph]_{C_G(\tup a)}^{C_H(\ph(\tup a))}\circ \change_{\ph(\tup a)}^{\tup b}.\]
Summing over all the pairs $(R,\ph)$ gives the formula as claimed.
\end{proof}

\begin{remark}\label{remarkTransferInPi0}
The formula of Proposition \ref{propSimpleGroupLoop} relates to class functions for groups as follows. Let $H\leq G$ finite groups. Consider the free abelian group $\ZZ[\pi_0\freeO n BG] = H_0(\freeO n BG)$. The ring of abelian group homomorphisms $\Hom(\ZZ[\pi_0\freeO n BG], \CC)$ is isomorphic to the ring of $\CC$-valued generalized class functions for $G$ (in the sense of \cite{HKR}). The transfer $\tr_H^G$ for $\ZZ[\pi_0\freeO n (-)]$ from the component $C_G(\tup a)$ to $C_G(\tup a)\cap gHg^{-1}$ is just multiplication with the index $\abs{C_G(\tup a)/C_G(\tup a)\cap gHg^{-1}}$. The induced map $(\tr_H^G)_*\colon\ZZ[\pi_0\freeO n BG]\to \ZZ[\pi_0\freeO n BH]$ therefore takes the form
\begin{align*}
(\tr_H^G)_*([\tup a]_G) &= \sum_{\substack{[\tup b]_H\text{ s.t.}\\ \exists g: g^{-1}\tup ag=\tup b}} \abs{C_G(\tup a)/C_G(\tup a)\cap gHg^{-1}} \cdot [\tup b]_H
\\ &= \sum_{\substack{[\tup b]_H\text{ s.t. }\\ \exists g: g^{-1}\tup ag=\tup b}} \sum_{z\in C_G(\tup a)/C_G(\tup a)\cap gHg^{-1}} [(zg)^{-1} \tup a (zg)]_H
\\ &= \sum_{\substack{g\in G/H\text{ s.t. }\\ g^{-1}\tup a g\in H}} [g^{-1}\tup a g]_H.
\end{align*}

Characters for complex $G$-representations are in particular $\CC$-valued class functions for $G$.
Since the components of $\freeO 1 BG$ correspond to conjugacy classes of elements in $G$, we can think of $\CC$-valued characters as giving class functions in $\Hom(\ZZ[\pi_0\freeO 1 BG],\CC)$. The biset structure on the functors $\freeO 1 (-)$ and $\ZZ[\pi_0\freeO 1 B(-)]$ induce a biset structure on class functions for finite groups. The calculations above and in Remark \ref{remarkSimpleLoopFunctor} show that the restriction and transfer maps considered here coincide with the usual formulas for restriction and induction of characters. In the case of $\freeO n (-)$ and $\ZZ[\pi_0\freeO n B(-)]$, these formulas are closely related to those that appear in \cite[Theorem D]{HKR}.
\end{remark}

\subsection{Evaluation maps and $\twistO n$ for finite groups}
In this subsection we will give explicit formulas for the evaluation map $(S^1)^n \times \freeO n (BG) \to BG$ and for $\twistO n X$, when $X$ is a $(G,H)$-biset. After this, we reformulate the results of Theorem \ref{thmCovLiftingFunctor} in terms of these formulas for bisets. 

The evaluation map $(S^1)^n\x \freeO n BG\xto\ev BG$, or for $\ell$ large enough, the map 
\[B(\Z/\ell)^n\x \Map(B(\Z/\ell)^n, BG) \xto\ev BG,\]
is modeled by groupoids as the functor
\[(\ast\mmod (\Z/\ell)^n)\x \Map(\ast\mmod (\Z/\ell)^n, \ast\mmod G)\to \ast \mmod G.\]
A morphism on the left is a pair consisting of a morphism $\tup t=(t_1,\dotsc,t_n)$ in $\ast\mmod (\Z/\ell)^n$ and a natural transformation $g\colon \tup a\Rightarrow \tup b$ between commuting $n$-tuples in $G$. The evaluation functor sends such a pair $(\tup t, g)$ to the group element/morphism $a_1^{t_1}\dotsm a_n^{t_n} g = g b_1^{t_1}\dotsm b_n^{t_n}$ in $\ast\mmod G$. Note that this is the evaluation map in the $2$-category of groupoids.

Via the equivalence $\freeO n BG\simeq B\freeO n G$ of Lemma \ref{lemmaGroupsAlgModelForFreeLoops}, we get the following algebraic model $\ev\colon (\Z/\ell)^n\x \freeO n G\to G$ for the evaluation map. For an $n$-tuple $\tup a$ of group elements  and an $n$-tuple $\tup t$ of integers, we let $\tup a^{\tup t}$ denote the pointwise exponentiation $(a_1^{t_1},\dotsc, a_n^{t_n})$.
\begin{lemma}\label{lemmaGroupEvaluation}
The evaluation map $\ev\colon (\Z/\ell)^n\x \freeO n G\to G$, restricted to the component $(\Z/\ell)^n\x C_G(\tup a)$ of $(\Z/\ell)^n\x \freeO n G$ corresponding to a commuting $n$-tuple $\tup a$, is given by the group homomorphism $\ev_{\tup a}\colon (\Z/\ell)^n\x C_G(\tup a)\to C_G(\tup a)\leq G$ with
\[\ev_{\tup a}(\tup t,z) = a_1^{t_1}\dotsm a_n^{t_n} \cdot z = (\prod_{i=1}^{n}a_i^{t_i})\cdot z\]
for all $\tup t\in (\Z/\ell)^n$ and $z\in C_G(\tup a)$.
\end{lemma}

In order to replace $\freeO n BG$ with $\Map(B(\Z/\ell)^n, BG)$ we need $\ell$ to be large enough, e.g. a multiple of $\abs{G}$. 
However the same $\ell$ does not work for all finite groups at the same time, and so the algebraic evaluation map $\ev\colon (\Z/\ell)^n\x \freeO n G \to G$ only makes sense on a subclass of the objects in $\AG$.

If $\AG_\ell$ is the full subcategory of $\AG$ spanned by the unions of finite groups where all elements have order dividing $\ell$, then we have a well defined functor $(\Z/\ell)^n\x \freeO n (-)\colon \AG_\ell\to \AG_\ell$ that does model $S^1\times \freeO n(B(-))$. When $\ell$ divides $m$ we have a fully faithful inclusion $\AG_\ell\subseteq \AG_m$ as well as a natural map $(\Z/m)^n\x \freeO n G \rightarrow (\Z/\ell)^n\x \freeO n G$ in $\AG_m$ that is compatible with the evaluation map, for $G\in \AG_\ell$. 

\begin{convention}\label{conventionOrderOfCyclics}
Given any finite collection of objects in $\AG$ it is possible to find some $\ell$ such that all of the objects are contained in $\AG_\ell$. As the value of $\ell$ does not play a critical role in what follows, we shall ignore this ambiguity from now on and suppose that we have chosen $\ell$ large enough for all the groups involved in the discussions.
\end{convention}

To construct the functors $\twistO n\colon \Cov\to \Cov$ in Section \ref{secCoveringMaps}, we took the pullback along $\ev\colon (S^1)^n\x \freeO n X\to X$ of a finite sheeted covering map $X\ffrom Y$ to form the square
\[
\begin{tikzpicture}[ampersand replacement=\&, baseline=(BASE).base]
\node (M) [matrix of math nodes] {
    \&[.5cm]\&[.5cm] \PB n(\pi) \&[.75cm]\&[.5cm] \\[1cm]
    \& (S^1)^n\x\freeO nX \&\& Y \& \\[1cm]
    \&\& X. \&\&  \\
};
\path [auto, arrow, ->] 
    (M-1-3) edge[->>] node[swap]{$\ev_X^*\pi$} (M-2-2)
            edge node{$\cover{\ev_X}$} (M-2-4)
    (M-2-2) edge node[swap]{$\ev_X$} (M-3-3)
    (M-2-4) edge[->>] node(BASE){$\pi$} (M-3-3)
;
\end{tikzpicture}
\]
For finite groups, the covering maps have the form $BG\xffrom{\tr_R^G} EG/R \simeq BR$ for subgroups $R\leq G$. By the equivalence of categories $\AG_+\simeq \Ho(\Cov)$ of Proposition \ref{propEmbedInHoCov}, the backward map $\tr_R^G$ is represented by the $(G,R)$-biset $[R,\id]_G^R = G_G^R$. Determining the pullback $\PB n(\tr_R^G)$ is equivalent to asking how the composite $[(\Z/\ell)^n\x \freeO n G, \ev]\cmp [R,\id]_G^R$ decomposes into transitive bisets in $\AG_+((\Z/\ell)^n\x \freeO n G, R)$.

Consider a commuting $n$-tuple $\tup a$ in $G$ and the corresponding component $(\Z/\ell)^n\x C_G(\tup a)$ of $(\Z/\ell)^n\x \freeO n G$. If we decompose $[(\Z/\ell)^n\x C_G(\tup a), \ev_{\tup a}]\cmp [R,\id]_G^R$ into orbits as a linear combination
\[[(\Z/\ell)^n\x C_G(\tup a), \ev_{\tup a}]\cmp [R,\id]_G^R=\sum_{\substack{D\leq (\Z/\ell)^n\x C_G(\tup a)\\\ph\colon D\to R}} c_{D,\ph}\cdot [D,\ph]_{(\Z/\ell)^n\x C_G(\tup a)}^R\]
with coefficients $c_{D,\ph}\in \Z_{\geq 0}$,
then the part of $\PB n(\tr_R^G)$ that sits over the component $B(\Z/\ell)^n\x BC_G(\tup a)$ is homotopy equivalent to the collection of $c_{D,\ph}$ copies of the classifying space $BD$ for each pair $(D,\ph)$, with the maps $B\ph\colon BD\to BR$ forming the components of $\PB n(\tr_R^G)\to BR$ in the pullback.

The same approach works for a general $(G,H)$-biset $M$. The composite $[(\Z/\ell)^n\x \freeO n G, \ev] \cmp M \in \AG_+((\Z/\ell)^n\x \freeO n G, H)$
models the composed span
\begin{equation}\label{eqBisetComposedWithEval}
\begin{tikzpicture}[ampersand replacement=\&, baseline=(BASE).base]
\scriptsize 
\matrix (M) [matrix of math nodes] {
  \&[0cm] \PB n(\pi) \&[0cm]   \&[-.1cm] \&[0.3cm]     \\[1cm]
 (S^1)^n\x B\freeO n G \&\& \abs{(M/H)\mmod G} \&\& \\[1cm]
 \& BG \&\& \abs{M\mmod G\x H} \& \\[1cm]
 \&\&\&\& BH. \\
};
\path [auto, arrow, ->]
    (M-1-2) edge[->>] node[swap]{$\ev^*\pi$} (M-2-1)
            edge  (M-2-3)
    (M-2-1) edge node[swap]{$\ev$} (M-3-2)
    (M-2-3) edge[->>] node{$\pi$} (M-3-2)
            edge node{$\simeq$} (M-3-4)
    (M-3-4) edge (M-4-5)
;
\end{tikzpicture}
\end{equation}
Hence if we decompose $[(\Z/\ell)^n\x C_G(\tup a), \ev_{\tup a}]\cmp M$ into orbits
\[
[(\Z/\ell)^n\x C_G(\tup a), \ev_{\tup a}]\cmp M=\sum_{\substack{D\leq (\Z/\ell)^n\x C_G(\tup a)\\\ph\colon D\to H}} c_{D,\ph}\cdot [D,\ph]_{(\Z/\ell)^n\x C_G(\tup a)}^H,
\]
then the part of $\PB n(\pi)$ that lies over the component $(S^1)^n\x B C_G(\tup a)$ is homotopy equivalent to $c_{D,\ph}$ copies of $B(D)$, and the maps $B\ph\colon BD \to BH$, put together, make up the map $\PB n(\pi)\to BH$.

The following lemma gives us some control on the subgroups $D\leq (\Z/\ell)^n\x C_G(\tup a)$ that appear in such pullbacks.
\begin{lemma}\label{lemmaGroupUpperPathOrbits}
Let $G$ and $H$ be finite groups and let $M\in \AG(G,H)$ be a virtual biset. Suppose $\tup a$ is an $n$-tuple of commuting elements in $G$.

Consider the restriction of $M$ to the centralizer of $\tup a$:
\[M_{C_G(\tup a)}^H=[C_G(\tup a),\incl]_{C_G(\tup a)}^G \cmp M\in \AG(C_G(\tup a), H).\]
Recalling Convention \ref{conventionOrbitDecomposition}, write $M_{C_G(\tup a)}^H$ as a linear combination of basis elements
\begin{equation} \label{eqDecomposeBisetOrbits}
M_{C_G(\tup a)}^H = \sum_{(R,\ph)} c_{R,\ph} \cdot [R,\ph]_{C_G(\tup a)}^H,
\end{equation}
where $R\leq C_G(\tup a)$ and $\ph\colon R\to H$.

The composite $\ev\cmp M \in \AG((\Z/\ell)^n\x \freeO nG, H)$ on the component $(\Z/\ell)^n\x C_G(\tup a)$ then satisfies the formula
\[(\ev\cmp M)_{\tup a} = \ev_{\tup a}\cmp M = \sum_{(R,\ph)} c_{R,\ph} \cdot [\ev_{\tup a}^{-1}(R), \ph\circ \ev_{\tup a}]_{(\Z/\ell)^n\x C_G(\tup a)}^H,\]
with the same pairs of $R\leq C_G(\tup a)$ and $\ph\colon R\to H$ as in Equation \eqref{eqDecomposeBisetOrbits}.
\end{lemma}

\begin{remark}
Following Convention \ref{conventionOrbitDecomposition}, the decomposition
\[M_{C_G(\tup a)}^H = \sum_{(R,\ph)} c_{R,\ph} \cdot [R,\ph]_{C_G(\tup a)}^H\]
has the sum run over all pairs $(R,\ph)$, with $R\leq C_G(\tup a)$ and $\ph\colon R\to H$, and not just $(C_G(\tup a),H)$-conjugacy classes of pairs. 

Part of the statement of the lemma is the claim that choosing a different linear expression for $M_{C_G(\tup a)}^H$, with different coefficients $c_{R,\ph}$ that are possibly non-zero for several pairs in the same $(C_G(\tup a),H)$-conjugacy class, will result in equal sums
\[\sum_{(R,\ph)} c_{R,\ph} \cdot [\ev_{\tup a}^{-1}(R), \ph\circ \ev_{\tup a}]_{(\Z/\ell)^n\x C_G(\tup a)}^H.\]
\end{remark}

\begin{proof}[Proof of Lemma \ref{lemmaGroupUpperPathOrbits}]
The evaluation map $\ev_{\tup a}\colon (\Z/\ell)^n\x C_G(\tup a)\to C_G(\tup a)\leq G$ has image $C_G(\tup a)$ on the component $(\Z/\ell)^n\x C_G(\tup a)$ of $(\Z/\ell)^n\x \freeO nG$. 
We can then proceed by the following calculation with bisets:
\begin{align*}
&[(\Z/\ell)^n\x C_G(\tup a),\ev_{\tup a}]_{(\Z/\ell)^n\x C_G(\tup a)}^G\cmp M_G^H
\\ ={}& [(\Z/\ell)^n\x C_G(\tup a), \ev_{\tup a}]_{(\Z/\ell)^n\x C_G(\tup a)}^{C_G(\tup a)} \cmp [C_G(\tup a),\incl]_{C_G(\tup a)}^G\cmp M_G^H
\\ ={}& [(\Z/\ell)^n\x C_G(\tup a), \ev_{\tup a}]_{(\Z/\ell)^n\x C_G(\tup a)}^{C_G(\tup a)} \cmp M_{C_G(\tup a)}^H
\\ ={}& [(\Z/\ell)^n\x C_G(\tup a), \ev_{\tup a}]_{(\Z/\ell)^n\x C_G(\tup a)}^{C_G(\tup a)} \cmp \Bigl(\sum_{(R,\ph)} c_{R,\ph} \cdot [R,\ph]_{C_G(\tup a)}^H\Bigr)
\\ ={}& \sum_{(R,\ph)} c_{R,\ph} \cdot \bigl([(\Z/\ell)^n\x C_G(\tup a), \ev_{\tup a}]_{(\Z/\ell)^n\x C_G(\tup a)}^{C_G(\tup a)} \cmp [R,\ph]_{C_G(\tup a)}^H\bigr)
\\ ={}& \sum_{(R,\ph)} c_{R,\ph} \cdot [\ev_{\tup a}^{-1}(R), \ph\circ \ev_{\tup a}]_{(\Z/\ell)^n\x C_G(\tup a)}^H.
\end{align*}
The last equality follows from the double coset formula for bisets \cite{Bouc}*{Lemma 2.3.24} in the special case $[R,\ph]_A^B\cmp [T,\psi]_ B^C = [\ph^{-1}(T), \psi\circ \ph]_A^C$ when $\ph\colon R\to B$ is surjective.
\end{proof}

For $c\in \Z_{\geq 0}$ and $X\in \Top$ let $c\cdot X$ denote $c$ disjoint copies of the space $X$.
Suppose $M\in\AG_+(G,H)$ is a biset with the associated covering map $\pi\colon (M/H)\mmod G \tto \ast\mmod G$.
The components of $\PB n(\pi)$ that lie over the component $(\Z/\ell)^n\x C_G(\tup a)$ of $(\Z/\ell)^n\x \freeO n G$ are then, according to Lemma \ref{lemmaGroupUpperPathOrbits}, homotopy equivalent to the disjoint union
\[\coprod_{(R,\ph)} c_{R,\ph}\cdot B(\ev_{\tup a}^{-1} R),\]
where the coefficients $c_{R,\ph}$ come from the decomposition of $M_{C_G(\tup a)}^H$.

The classifying space of the action groupoid $(M/H)\mmod C_G(\tup a)$ is equivalent to $\coprod_{(R,\ph)} c_{R,\ph}\cdot BR$ according to the decomposition of the biset $M_{C_G(\tup a)}^H\in \AG_+(C_G(\tup a),H)$. These decompositions into orbits fit into the following diagram of groupoids:
\[
\begin{tikzpicture}[ampersand replacement=\&, baseline=(BASE).base]
\matrix (M) [matrix of math nodes] {
 \displaystyle\coprod_{(R,\ph)} c_{R,\ph} \cdot \ev_{\tup a}^{-1} R \&[2.5cm] \coprod_{(R,\ph)} c_{R,\ph}\cdot R \&[2.5cm]     \\[1cm]
 (M/H)\mmod ((\Z/\ell)^n\x C_G(\tup a)) \& (M/H)\mmod C_G(\tup a) \& (M/H)\mmod G \\[1cm]
 (\Z/\ell)^n \x C_G(\tup a) \& C_G(\tup a) \& G. \\
};
\path [auto, arrow, ->]
    (M-1-1) edge node{$\simeq$} (M-2-1)
            edge node{$\ev_{\tup a}$} (M-1-1 -| M-1-2.west)
    (M-1-2) edge node{$\simeq$} (M-2-2)
    (M-2-1) edge[->>] node[swap]{$\ev_{\tup a}^*\pi$} (M-3-1)
            edge node{$(M/H)\mmod \ev_{\tup a}$} (M-2-2)
    (M-2-2) edge[->>] node[swap]{$\incl^*\pi$} (M-3-2)
            edge node{$(M/H)\mmod \incl$} (M-2-3)
    (M-2-3) edge[->>] node{$\pi$} (M-3-3)
    (M-3-1) edge node{$\ev_{\tup a}$} (M-3-2)
    (M-3-2) edge node{$\incl$} (M-3-3)
;
\end{tikzpicture}
\]
The two lower squares are pullbacks, and the coefficients $c_{R,\ph}$ are as in Lemma \ref{lemmaGroupUpperPathOrbits}. We use the coproduct in the top left as an algebraic model for $\PB n(\pi)$.

\begin{definition}\label{defGroupPB}
For a biset $M\in \AG_+(G,H)$, we let $\PB n(M)$ denote the following algebraic model for $\PB n(\pi)$, with $\pi\colon (M/H)\mmod G \tto BG$:
\[
\PB n(M) = \coprod_{[\tup a]\in \cntuples n G} \PB n(M)_{\tup a}, 
\]
where
\[\PB n(M)_{\tup a}= \coprod_{(R,\ph)} c_{R,\ph}\cdot \ev_{\tup a}^{-1}(R).\]
Here, the coefficients $c_{R,\ph}$ come from the decomposition of $M_{C_G(\tup a)}^H$ in Lemma \ref{lemmaGroupUpperPathOrbits} and $\tup a \in [\tup a]$ is our chosen representative of the conjugacy class $[\tup a]\in \cntuples n G$. Note that $\PB n(M)_{\tup a}$ is the part of $\PB n (M)$ lying over the component $(\Z/\ell)^n\x C_G(\tup a)$ of $(\Z/\ell)^n\x \freeO nG$.
\end{definition}

The diagram \eqref{eqBisetComposedWithEval} involving $\PB n(\pi)$ in $\Cov$ has the following algebraic model involving $\PB n(M)_{\tup a}$ when we restrict to the component $(\Z/\ell)^n\x C_G(\tup a)$:
\[
\begin{tikzpicture}[ampersand replacement=\&, baseline=(BASE).base]
\scriptsize 
\matrix (M) [matrix of math nodes] {
  \&[0cm] \PB n(M)_{\tup a} \&[0cm]   \&[-.1cm] \&[0.3cm]     \\[1cm]
 (\Z/\ell)^n\x C_G(\tup a) \&\& \displaystyle\coprod_{(R,\ph)}c_{R,\ph}\cdot R \&\& \\[1cm]
 \& C_G(\tup a) \&\& (M/H)\mmod G \& \\[1.1cm]
 \&\& G\& \& H. \\
};
\path [auto, arrow, ->]
    (M-1-2) edge  (M-2-1)
            edge node{$\ev_{\tup a}$} (M-2-3)
    (M-2-1) edge node[swap]{$\ev_G$} (M-3-2)
    (M-2-3) edge (M-3-2)
            edge (M-3-4)
    (M-3-2) edge node[swap]{$\incl_{C_G(\tup a)}^G$} (M-4-3)
    (M-3-4) edge (M-4-5)
            edge node{$\pi$} (M-4-3)
    (M-2-3) edge[bend left=45] node{$\ph$} (M-4-5)
;
\end{tikzpicture}
\]
Here the backward maps are inclusions of subgroups (corresponding to finite covering maps).

Next we shall give an algebraic version  of the map $k\colon \PB n(\pi)\to (\Z_{>0})^n$ from  Definition \ref{defCovKMap}. Fix a component $\ev_{\tup a}^{-1}(R)$ of $\PB n(M)$ lying over $(\Z/\ell)^n\x C_G(\tup a)$. 
The corresponding component $E((\Z/\ell)^n\x C_G(\tup a))/\ev^{-1}_{\tup a} R$ of $\PB n(\pi)$ sits over the component of $(S^1)^n\x \freeO n G$ indexed by $\tup a$. Thus the corresponding $n$-fold loop $(S^1)^n\to BG$ is just picking out the commuting $n$-tuple $\tup a$ in $\pi_1(BG)\cong G$. By Remark \ref{remarkKMap} the exponent $k_i$ for each $1\leq i\leq n$ is the smallest positive integer such that $a_i^{k_i}\in \pi_1(BC_G(\tup a))\leq \pi_1(BG)$ lies in the further subgroup $\pi_1(BR)\cong R$. Also, the exponent $k_i$ is constant on the entire component $E((\Z/\ell)^n\x C_G(\tup a))/\ev^{-1}_{\tup a} R$ of $\PB n(\pi)$. Therefore, we have the following algebraic analogue of the map $k\colon \PB n(\pi)\to (\Z_{>0})$:

\begin{definition}\label{defGrpKMap}
Recall that $\cntuples nG$ denotes the set of $G$-conjugacy classes of commuting $n$-tuples of $G$. For a group $H$, let $\SubGrp(H)$ be the set of subgroups of $H$. We define 
\[
k\colon \coprod_{[a]\in \cntuples n G} \SubGrp(C_G(\tup a)) \to (\Z_{>0})^n
\]
by 
\[k([\tup a]\in \cntuples n G, R\leq C_G(\tup a))_i = (\text{the smallest $k_i>0$ such that $a_i^{k_i}\in R$})\]
for each $1\leq i\leq n$.
\end{definition}

\begin{remark}
The $k$-exponents of Definitions \ref{defCovKMap} and \ref{defGrpKMap} agree:
Given a component $E((\Z/\ell)^n\x C_G(\tup a))/\ev^{-1}_{\tup a} R$ of $\PB n(\pi)$ the map $k\colon \PB n(\pi)\to (\Z_{>0})^n$ takes the constant value $k(\tup a, R)\in (\Z_{>0})^n$ on this component.
\end{remark}

In $\Cov$, $\PBlift n(\pi)$ consists of those components in $\PB n(\pi)$ where the $n$-fold loop lifts to an $n$-fold loop through $\pi$, or equivalently those components where $k_i\colon \PB n(\pi)\to \Z_{>0}$ equals $1$ for all $1\leq i\leq n$.
The following lemma characterizes when a component of $\PB n(M)$ is the algebraic model for a component of $\PBlift n(\pi)$:

\begin{lemma}\label{lemGroupTwistedSubgroups}
Let $\tup a$ be a commuting $n$-tuple in $G$, and let $R\leq C_G(\tup a)$ be a subgroup. The following four statements about the preimage $\ev_{\tup a}^{-1}(R)$ in $(\Z/\ell)^n\x C_G(\tup a)$ are equivalent:
\begin{enumerate}
\renewcommand{\theenumi}{$(\roman{enumi})$}\renewcommand{\labelenumi}{\theenumi}
\item\label{itemGroupExponents} $k(\tup a,R)_i=1$ for all $1\leq i\leq n$.
\item\label{itemGroupTupleInSubgroup} All elements of $\tup a$ lie in the subgroup $R$.
\item\label{itemGroupSubgroupProductAB} $\ev_{\tup a}^{-1}(R)$ is a product of subgroups $A\x B$ with $A\leq (\ZZ/\ell)^n$ and $B\leq C_G(\tup a)$.
\item\label{itemGroupSubgroupProductR} $\ev_{\tup a}^{-1}(R) = (\Z/\ell)^n\x R$.
\end{enumerate}
\end{lemma}

\begin{proof}
Note first that $\tup a\in C_G(\tup a)$ since the elements of the $n$-tuple commute.

\ref{itemGroupExponents} $\Leftrightarrow$ \ref{itemGroupTupleInSubgroup}: Immediate by the definition of $k(\tup a, R)\in (\Z_{>0})^n$.

\ref{itemGroupTupleInSubgroup} $\Rightarrow$ \ref{itemGroupSubgroupProductR}: If all elements of $\tup a$ lie in the subgroup $R\leq C_G(\tup a)$, then so do their powers, and $\ev_{\tup a}^{-1}(R) = (\Z/\ell)^n\x R$ is immediate.

\ref{itemGroupSubgroupProductR} $\Rightarrow$ \ref{itemGroupSubgroupProductAB}: Property \ref{itemGroupSubgroupProductAB} is just a weakening of \ref{itemGroupSubgroupProductR}.

\ref{itemGroupSubgroupProductAB} $\Rightarrow$ \ref{itemGroupTupleInSubgroup}: If $(\tup t,z)\in \ev_{\tup a}^{-1}(R)$, and if $\tup t'\in (\Z/\ell)^n$ is any tuple, then the pair
$(\tup t - \tup t',(\prod_{i=1}^n a_i^{t'_i})\cdot z)$ is also in the preimage $\ev_{\tup a}^{-1}(R)\leq (\Z/\ell)^n\x C_G(\tup a)$. Hence if $\ev_{\tup a}^{-1}(R)=A\x B$, we must have $A=(\Z/\ell)^n$, and
in particular we have $(\Z/\ell)^n\x 1\leq \ev_{\tup a}^{-1}(R)$. The elements of $\tup a$, being the images under $\ev_{\tup a}$ of the generators from $(\Z/\ell)^n$, therefore lie in $R$.
\end{proof}

\begin{definition}\label{defGroupPBlift}
Given a biset $M\in \AG_+(G,H)$ and a commuting $n$-tuple $\tup a$ in $G$, we let $\PBlift n(M)_{\tup a}$ denote the collection of components of $\PB n(M)_{\tup a}$ satisfying the properties of Lemma \ref{lemGroupTwistedSubgroups}, so that
\[
\PBlift n(M)_{\tup a} = \coprod_{\substack{\tup a\in R\leq C_G(\tup a)\\\ph\colon R\to H}} c_{R,\ph}\cdot \ev_{\tup a}^{-1}(R) = \coprod_{\substack{\tup a\in R\leq C_G(\tup a)\\\ph\colon R\to H}} c_{R,\ph}\cdot ((\Z/\ell)^n\x R),
\]
where the coefficients $c_{R,\ph}$ come from the decomposition of $M_{C_G(\tup a)}^H$ in Lemma \ref{lemmaGroupUpperPathOrbits}.

We also let $\PBlift n(M)$ denote the union $\coprod_{[\tup a]\in\cntuples nG} \PBlift n(M)_{\tup a}$ over all our chosen representatives $\tup a$ of conjugacy classes of commuting $n$-tuples.
\end{definition}

The final piece needed to describe the functors $\twistO n\colon \AG\to \AG$ is the winding map $\wind(\pi)\colon \PB n(\pi)\to \PBlift n(\pi)$. The map $\wind(\pi)$ winds an $n$-fold loop $(S^1)^n\to BG$ around itself $k_i$ times in the $i$th coordinate direction. 
Given a tuple $\tup a$, we cannot assume that the coordinate-wise power $\tup a^{k(\tup a,R)}=(a_1^{k(\tup a,R)_1},\dotsc, a_n^{k(\tup a,R)_n})$ is always the representative of its conjugacy class in $G$. We therefore suppose that a choice of representatives has been made in advance as in Convention \ref{conventionGroupTupleReps}.

The algebraic model $\wind(M)\colon \PB n(M) \to \PBlift n(M)$ for $\wind(\pi)$ then takes a component $\ev_{\tup a}^{-1}(R)$ of $\PB n(M)_{\tup a}$ and sends it to $\PBlift n(M)_{\tup b}$, where $\tup b$ is the chosen representative of the $G$-conjugacy class of the tuple $\tup a^{k(\tup a,R)}=(a_1^{k(\tup a,R)_1},\dotsc, a_n^{k(\tup a,R)_n})$. 

Choose some $g\in G$ with $c_g(\tup a^{k(\tup a,R)}) = \tup b$ as in Lemma \ref{lemmaGroupCentralizerIso}. The part $\PBlift n(M)_{\tup b}$ for the tuple $\tup b$ has the form
\[\PBlift n(M)_{\tup b} = \coprod_{\substack{\tup b\in R'\leq C_G(\tup b)\\ \ph'\colon R'\to H}} c'_{R',\ph'}\cdot\ev_{\tup b}^{-1}(R'),\]
where the coefficients $c'_{R',\ph'}$ come from a decomposition of $M_{C_G(\tup b)}^H$ instead of $M_{C_G(\tup a)}^H$. In particular, since $C_G(\tup a)\leq C_G(\tup a^{k(\tup a,R)})\xto[c_g]{\cong} C_G(\tup b)$, each orbit of $M_{C_G(\tup a)}^H$ is contained in a unique orbit of $M_{C_G(\tup a^{k(\tup a,R)})}^H$ though several orbits of $M_{C_G(\tup a)}^H$ might combine to give a single orbit of $M_{C_G(\tup a^{k(\tup a,R)})}^H\cong c_g^* (M_{C_G(\tup b)}^H)$.

Given a component $\ev_{\tup a}^{-1} R$ of $\PB n(M)_{\tup a}$, this component corresponds to the orbit of a point $x\in M_{C_G(\tup a)}^H$ with orbit of the form $[R,\ph]_{C_G(\tup a)}^H$. The winding map $\wind(\pi)$ in $\Cov$ only changes the $n$-fold loop and keeps the same point in the fiber of $\pi$. This means that the algebraic version of $\wind(\pi)$, $\wind(M)$, has to send $\ev_{\tup a}^{-1} R$ to $\PBlift n(M)_{\tup b}$ but hit the orbit of the same point $x\in M_{C_G(\tup a^{k(\tup a,R)})}^H\cong c_g^*(M_{C_G(\tup b)}^H)$. Under the isomorphism $c_g\colon M_{C_G(\tup a^{k(\tup a,R)})}^H \xto{\cong} M_{C_G(\tup b)}^H$ the orbit of $x$ in $M_{C_G(\tup a^{k(\tup a,R)})}^H$ corresponds to the orbit of $gx$ in $M_{C_G(\tup b)}^H$. Note that a different choice of $g$ conjugating $\tup a^{k(\tup a,R)}$ to $\tup b$ still gives the same orbit of $gx$ in $M_{C_G(\tup b)}^H$, since different $g$'s only differ by elements of $C_G(\tup b)$.

The corresponding orbit of $M_{C_G(\tup b)}^H$ has the form $[R',\ph']_{C_G(\tup b)}^H$, where $c_g(R)\leq R'$ and $(\ph'\circ c_g)|_R = \ph$, and $\wind(M)$ takes the component $\ev_{\tup a}^{-1} R$ of $\PB n(M)_{\tup a}$ to the component $\ev_{\tup b}^{-1} R'$ of $\PBlift n(M)_{\tup b}$.
By Lemma \ref{lemGroupTwistedSubgroups} it follows that $\ev_{\tup b}^{-1} R' = (\Z/\ell)^n\x R'$  since $\tup b=c_g(\tup a^{k(\tup a,R)}) \in c_g(R)\leq R'$.

In order to work out a formula for the map $\ev_{\tup a}^{-1} R\to (\Z/\ell)^n\x R'$, recall that $\wind(M)$ commutes with evaluation maps (Remark \ref{remarkWindRespectsEval}) so that we must have a commutative diagram
\[
\begin{tikzpicture}[ampersand replacement=\&, baseline=(BASE).base]
\matrix (M) [matrix of math nodes] {
 \PB n(M)_{\tup a}\geq \ev_{\tup a}^{-1} R  \&[2.5cm] (\Z/\ell)^n\x g R'g^{-1} \&[2.5cm] (\Z/\ell)^n\x R' \leq \PBlift n(M)_{\tup b}     \\[2cm]
  \& G. \&  \\
};
\path [auto, arrow, ->]
    (M-1-1) edge node{$?$} (M-1-2)
    (M-1-1.south east) edge node[swap]{$\ev_{\tup a}$} (M-2-2)
    (M-1-1.north east) edge[bend left=30] node{$\wind(M)$} (M-1-3.north west)
    (M-1-2) edge node[swap]{$\cong$} node{$(\Z/\ell)^n\x c_g$} (M-1-3)
            edge node{$\ev_{\tup a^{k(\tup a,R)}}$} (M-2-2)
    (M-1-3.south west) edge node{$\ev_{\tup b}$} (M-2-2)
;
\end{tikzpicture}
\]
In $\Cov$, the map $\wind(\pi)$ leaves the $(S^1)^n$-coordinate unchanged and only uses it to shift the $n$-fold loop. Consequently, $\wind(M)$ must take a pair $(\tup t,z)\in \ev_{\tup a}^{-1} R$ to a pair of the form $(\tup t,z')\in (\Z/\ell)^n\x gR'g^{-1}$ satisfying 
\[(\prod_{i=1}^n a_{i}^{t_i})\cdot z = \ev_{\tup a}(\tup t,z) =\ev_{\tup a^{k(\tup a,R)}}(\tup t, z') = (\prod_{i=1}^n a_i^{k(\tup a,R)_i\cdot t_i})\cdot z'.\]
It follows that $z'$ is uniquely determined with 
\[z'=(\prod_{i=1}^n a_i^{t_i-k(\tup a,R)_i\cdot t_i})\cdot z.\] 
Thus the definition of $\wind(\pi)$ has determined the following definition of the map $\wind(\tup a,R)\colon \ev_{\tup a}^{-1} R \to (\Z/\ell)^n\x R$, which we will use to give a formula for $\wind(M)$ and $\twistO n(M)$.

\begin{lemma}\label{lemmaUntwisting}
Consider a commuting $n$-tuple $\tup a$ in $G$ and a subgroup $R\leq C_G(\tup a)$. Recall that the tuple of exponents $k(\tup a,R)\in (\Z_{>0})^n$ are minimal such that the tuple $\tup a^{k(\tup a,R)}$ lies in the subgroup $R$.

We then have an isomorphism $\wind(\tup a,R)\colon \ev_{\tup a}^{-1}(R) \xto{\cong} \ev^{-1}_{\tup a^{k(\tup a,R)}} R=(\Z/\ell)^n\x R$ given by
\[\wind(\tup a,R)(\tup t,z) = (\tup t, a_1^{t_1 - k(\tup a,R)_1 t_1}\dotsm a_n^{t_n-k(\tup a,R)_n t_n}\cdot z),\]
for pairs $(\tup t,z)\in \ev_{\tup a}^{-1}(R)$. Furthermore, $\wind(\tup a,R)$ commutes with evaluation maps:
\[\ev_{\tup a^{\tup k(\tup a,R)}} \circ \wind(\tup a,R) = \ev_{\tup a}|_{\ev_{\tup a}^{-1}(R)}.\]
\end{lemma}

\begin{proof}
A pair $(\tup t,z)\in (\Z/\ell)^n\x C_G(\tup a)$ lies in the preimage $ \ev_{\tup a}^{-1}(R)$ if and only if $(\prod_{i=1}^n a_i^{\tup t_i})\cdot z\in R$. At the same time we also have $a_i^{k(\tup a,R)_i}\in R$, for every $1\leq i \leq n$. It follows that, for each pair $(\tup t,z)\in \ev_{\tup a}^{-1}(R)$, we have $(\tprod_{i=1}^n (a_i^{k(\tup a,R)_i})^{-t_i}\cdot a_i^{t_i})\cdot z \in R$. Thus $\wind(\tup a,R)$ lands in the product $(\Z/\ell)^n\x R$. Since $z$ along with all powers of each $a_i$ commute, $\wind(\tup a,R)$ is a group homomorphism.

The inverse homomorphism $\wind(\tup a,R)^{-1}$ has the straightforward expression
\[\wind(\tup a,R)^{-1}(\tup t,z) = (\tup t,(\prod_{i=1}^n a_i^{k(\tup a,R)_i t_i - t_i})\cdot z),\]
and for all $(\tup t,z)\in (\Z/\ell)^n\x R$ we have $(\tup t,(\tprod_{i=1}^n a_i^{k(\tup a,R)_i t_i - t_i})\cdot z) \in \ev_{\tup a}^{-1}(R)$ since $z$ and $a_i^{k(\tup a,R)_i t_i}$ lie in $R$ already.

Finally, $\wind(\tup a,R)$ commutes with evaluation maps when we use $\ev_{\tup a^{\tup k(\tup a,R)}}$ for the codomain of $\wind(\tup a,R)$:
\begin{multline*}
{}\ev_{\tup a^{\tup k(\tup a,R)}}(\wind(\tup a,R)(\tup t,z))
={}\ev_{\tup a^{\tup k(\tup a,R)}}(\tup t,(\prod_{i=1}^n a_i^{t_i-k(\tup a,R)_i t_i})\cdot z) 
={}\\  (\prod_{i=1}^n (a_i^{k(\tup a,R)_i})^{t_i}\cdot a_i^{t_i-k(\tup a,R)_i t_i})\cdot z  
={}(\prod_{i=1}^n a_i^{t_i})\cdot z 
={} \ev_{\tup a}(\tup t,z)
\end{multline*}
for all $(\tup t,z)\in \ev_{\tup a}^{-1}(R)$.
\end{proof}

The discussion over the last few pages gives us the following expression for the map 
\[
\wind(M)\colon \PB n(M)\to \PBlift n(M)
\]
in terms of the maps $\wind(\tup a,R)$ above:

\begin{lemma}\label{lemGroupWindingMap}
Let $M\in \AG_+(G,H)$ be a biset for finite groups $G$ and $H$. Given a component $\ev_{\tup a}^{-1} R$ of $\PB n(M)_{\tup a}$, let $\tup b$ be the chosen representative in the $G$-conjugacy class of $\tup a^{k(\tup a,R)}$, and let $g\in G$ be an element conjugating $\tup a^{k(\tup a,R)}$ to $\tup b$.

If the component $\ev_{\tup a}^{-1} R$ of $\PB n(M)_{\tup a}$ corresponds to the orbit of a point $x\in M_{C_G(\tup a)}^H$, then the map $\wind(M)\colon \PB n(M)\to \PBlift n(M)$ takes $\ev_{\tup a}^{-1} R$ to the component $\ev_{\tup b}^{-1} R'=(\Z/\ell)^n\x R'$ of $\PBlift n(M)_{\tup b}$ corresponding to the orbit of the point $gx\in M_{C_G(\tup b)}^H$. Furthermore, the map $\wind(M)|_{\ev_{\tup a}^{-1} R} \colon \ev_{\tup a}^{-1} R\to \ev_{\tup b}^{-1} R'$ equals the composite
\[
\ev_{\tup a}^{-1} R\xto{\wind(\tup a,R)} (\Z/\ell)^n\x R \xto{(\Z/\ell)^n\x c_g} (\Z/\ell)^n\x g^{-1}Rg \leq (\Z/\ell)^n\x R'.
\]
\end{lemma}

We now have all the ingredients needed to express the endofunctor $\twistO n\colon \AG_+\to \AG_+$  as well as the linear extension $\twistO n\colon \AG\to \AG$. Recall that $\twistO n\colon \AG_+\to \AG_+$ is induced by $\twistO n\colon \Cov\to \Cov$ via the fully faithful inclusion $\AG_+ \hookrightarrow \Ho(\Cov)$.

\begin{prop}\label{propGroupTwistedLoopFunctor}
The functor $\twistO n\colon \AG_+\to \AG_+$ (or $\twistO n\colon \AG\to \AG$) sends a $(G,H)$-biset $M$ to a biset matrix $\twistO n(M)\in \AG_+((\Z/\ell)^n \times \freeO nG, (\Z/\ell)^n \times \freeO nH)$, with entries indexed by representatives $\tup a$ and $\tup b$ for the conjugacy classes of commuting $n$-tuples in $G$ and $H$ respectively. The matrix entries satisfy the formula
\begin{multline*}
\twistO n(M)_{\tup a,\tup b} \\= \hspace{-.4cm} \sum_{\substack{(R,\ph) \text{ s.t. } \ph(\tup a^{k(\tup a,R)})\\\text{is $H$-conj. to $\tup b$}}}\hspace{-.4cm} c_{R,\ph}\cdot 
\Bigl( [\ev_{\tup a}^{-1}(R),(\id_{(\Z/\ell)^n}\x\ph)\circ \wind(\tup a, R)]_{(\Z/\ell)^n\x C_G(\tup a)}^{(\Z/\ell)^n\x C_H(\ph(\tup a^{k(\tup a,R)}))} \\ \cmp ((\Z/\ell)^n\x \change_{\ph(\tup a^{k(\tup a,R)})}^{\tup b}) \Bigr).
\end{multline*}
Here the sum runs over those $R\leq C_G(\tup a)$ and $\ph\colon R\to H$ such that $\ph(\tup a^{k(\tup a,R)})$ is $H$-conjugate to $\tup b$. The coefficients $c_{R,\ph}$ arise from the decomposition 
\[M_{C_G(\tup a)}^H = \sum_{\substack{R\leq C_G(\tup a)\\ \ph\colon R\to H}} c_{R,\ph} \cdot [R,\ph]_{C_G(\tup a)}^H\]
as in Lemma \ref{lemmaGroupUpperPathOrbits}.
\end{prop}

\begin{proof}
Given a biset $M\in \AG_+(G,H)$ the corresponding span in $\Cov$ is given by the span of groupoids $G\xffrom{\pi} (M/H)\mmod G\xto{f} H$ according to Definition \ref{defBisetsToCov}.
The span $\twistO n(\pi,f)$ in $\Cov$ is then
\[
\begin{tikzpicture}[ampersand replacement=\&, baseline=(BASE).base]
\scriptsize 
\matrix (M) [matrix of math nodes] {
  \&[.2cm] \PB n(\pi)  \&[.2cm]   \&[-.5cm] \&[-.5cm]     \\[1cm]
 (S^1)^n\x \freeO n BG \&\& \PBlift n(\pi) \&\& \\[1cm]
 \&  \&\& (S^1)^n\x \freeO n(\abs{(M/H)\mmod G}) \& \\[1.1cm]
 \&\& \& \& (S^1)^n\x \freeO nBH. \\
};
\path [auto, arrow, ->]
    (M-1-2) edge node[swap]{$\ev_G^*\pi$} (M-2-1)
            edge node{$\wind(\pi)$} (M-2-3)
    (M-2-3) edge node{$\simeq$} (M-3-4)
    (M-3-4) edge node{$(S^1)^n\x\freeO n(f)$} (M-4-5)
;
\end{tikzpicture}
\]
With the algebraic models for $\PB n(\pi)$, $\PBlift n(\pi)$ and $\wind(\pi)$ given in Definitions \ref{defGroupPB}, \ref{defGroupPBlift}, and Lemma \ref{lemGroupWindingMap}, we have the following span of unions of groups modeling $\twistO n(\pi,f)$:
\[
\begin{tikzpicture}[ampersand replacement=\&, baseline=(BASE).base]
\scriptsize 
\matrix (M) [matrix of math nodes] {
  \&[-.5cm] \displaystyle\coprod_{\tup a\in \cntuples nG}\coprod_{\substack{R\leq C_G(\tup a)\\ \ph\colon R\to H}} c_{R,\ph}\cdot \ev_{\tup a}^{-1}R \&[-1cm]   \&[-2.3cm] \&[-2cm]     \\[1cm]
 \displaystyle\coprod_{\tup a\in \cntuples nG}(\Z/\ell)^n\x C_G(\tup a) \&\& \displaystyle\coprod_{\tup a'\in \cntuples nG}\coprod_{\substack{\tup a'\in R'\leq C_G(\tup a')\\ \ph'\colon R'\to H}} c_{R',\ph'}\cdot (\Z/\ell)^n\x R' \&\& \\[1cm]
 \&  \&\& (\ast\mmod (\Z/\ell)^n)\x \Map(\ast\mmod(\Z/\ell)^n ,(M/H)\mmod G) \& \\[1.1cm]
 \&\& \& \& \displaystyle\coprod_{\tup b\in \cntuples nH}(\Z/\ell)^n\x  C_H(\tup b). \\
};
\path [auto, arrow, ->]
    (M-1-2) edge node[swap]{$\incl$} (M-2-1)
            edge node{$\wind(M)$} (M-2-3)
    (M-2-3) edge node{$\simeq$} (M-3-4)
    (M-3-4) edge (M-4-5)
    (M-2-3.south) edge[bend right=55] node[swap]{$\id_{(\Z/\ell)^n}\x (\change_{\ph'(\tup a')}^{\tup b}\circ \ph')$} (M-4-5.west)
;
\end{tikzpicture}
\]
The last isomorphism $\change_{\ph'(\tup a')}^{\tup b}\colon C_H(\ph'(\tup a')) \xto{\cong} C_H(\tup b)$ is needed because the chosen representative $\tup b$ for the $H$-conjugacy class of $\ph'(\tup a')$ might be different from $\ph'(\tup a')$.
This is similar to how the formula for $\wind(M)$ involves a conjugation $c_g\colon C_G(\tup a^{k(\tup a,R)}) \xto{\cong} C_G(\tup a')$, so that $\wind(M)$ on the component $\ev_{\tup a}^{-1} R$ corresponding to the orbit of $x\in M_{C_G(\tup a)}^H$ takes the form
\[\ev_{\tup a}^{-1} R\xto[\cong]{\wind(\tup a,R)} (\Z/\ell)^n\x R \xto[\cong]{\id_{(\Z/\ell)^n}\x c_g} (\Z/\ell)^n \x g^{-1} R g \leq (\Z/\ell)^n \x R'\]
for the component $\ev_{\tup a'}^{-1} R'= (\Z/\ell)^n \x R'$ corresponding to the (possibly larger) orbit of the same point $x\in M_{C_G(\tup a')}^H$, with $\tup a'$ representing the $G$-conjugacy class of $\tup a^{k(\tup a,R)}$.

Because the orbit of $x\in  M_{C_G(\tup a)}^H$ has the form $[R,\ph]_{C_G(\tup a)}^H$ and the orbit of $x\in M_{C_G(\tup a^{k(\tup a,R)})}^H$ has the form $[gR'g^{-1},\ph'\circ c_g]_{C_G(\tup a^{k(\tup a,R)})}^H$, we have $R\leq gR'g^{-1}$ and $(\ph'\circ c_g)|_R = \ph$.
Consequently $\ph'(\tup a')=\ph(\tup a^{k(\tup a,R)})$, and the composite
\[ (\Z/\ell)^n\x R \xto[\cong]{\id_{(\Z/\ell)^n}\x c_g} (\Z/\ell)^n \x g^{-1} R g \leq (\Z/\ell)^n \x R'\xto{\id_{(\Z/\ell)^n}\x\ph'} (\Z/\ell)^n\x C_H(\ph'(\tup a'))\]
is just $\id_{(\Z/\ell)^n}\x \ph\colon (\Z/\ell)^n \x R\to (\Z/\ell)^n\x C_H(\ph(\tup a^{k(\tup a,R)}))$. The span above simplifies to
\[
\begin{tikzpicture}[ampersand replacement=\&, baseline=(BASE).base]
\scriptsize 
\matrix (M) [matrix of math nodes] {
  \&[-.5cm] \displaystyle\coprod_{\tup a\in \cntuples nG}\coprod_{\substack{R\leq C_G(\tup a)\\ \ph\colon R\to H}} c_{R,\ph}\cdot \ev_{\tup a}^{-1}R \&[-.7cm]        \\[1cm]
 \displaystyle\coprod_{\tup a\in \cntuples nG}(\Z/\ell)^n\x C_G(\tup a) \&\& \displaystyle\coprod_{\tup b\in \cntuples nH}(\Z/\ell)^n\x  C_H(\tup b). \\
};
\path [auto, arrow, ->]
    (M-1-2) edge node[swap]{$\incl$} (M-2-1)
            edge node{$\id_{(\Z/\ell)^n}\x(\change_{\ph(\tup a^{k(\tup a,R)})}^{\tup b}\circ \ph)\circ \wind(\tup a,R)$} (M-2-3)
;
\end{tikzpicture}
\]
Finally the matrix entry $\twistO n(M)_{\tup a,\tup b}$ picks out the part of the span above that lies over the components $(\Z/\ell)^n\x C_G(\tup a)$ and $(\Z/\ell)^n\x C_H(\tup b)$, and thus only those $\ev_{\tup a}^{-1}$ where $\tup b$ represents the $H$-conjugacy class of $\ph(\tup a^{k(\tup a,R)})$. Those components form the span of groups corresponding to the biset
\begin{multline*}
    \sum_{\substack{(R,\ph) \text{ s.t. } \ph(\tup a^{k(\tup a,R)})\\\text{is $H$-conj. to $\tup b$}}}\hspace{-.4cm} c_{R,\ph}\cdot 
\Bigl( [\ev_{\tup a}^{-1}(R),(\id_{(\Z/\ell)^n}\x\ph)\circ \wind(\tup a, R)]_{(\Z/\ell)^n\x C_G(\tup a)}^{(\Z/\ell)^n\x C_H(\ph(\tup a^{k(\tup a,R)}))} \\ \cmp ((\Z/\ell)^n\x \change_{\ph(\tup a^{k(\tup a,R)})}^{\tup b}) \Bigr)
\end{multline*}
as claimed. Because the formulas in Proposition \ref{propGroupTwistedLoopFunctor} are linear in the biset $M\in \AG_+(G,H)$, the formulas hold for the endofunctor $\twistO n\colon \AG\to \AG$ as well.
\end{proof}

\begin{remark}\label{remarkGroupUntwisted}
Instead of applying $\wind(\tup a,R)$ to turn $\ev^{-1}_{\tup a}(R)$ into $(\Z/\ell)^n\x R$, we can alternatively take only those summands of Proposition \ref{propGroupTwistedLoopFunctor} for which $\ev^{-1}_{\tup a}(R)$ is already equal to $(\Z/\ell)^n\x R$. By Lemma \ref{lemGroupTwistedSubgroups}, this happens if and only if the subgroup $R\leq C_S(\tup a)$ contains $\tup a$.
Taking the sum over pairs $(R,\ph)$ with $\tup a\in R$ corresponds to Definition \ref{defCovUntwistedLn}, where we restrict to $\PBlift n(\pi) \subseteq \PB n(\pi)$, which corresponds to those components of the free loop space for which the loops in the base space lift to actual loops in the total space. For finite groups $G$ and $H$, and a virtual biset matrix $M\in \AG(G,H)$, we define $\untwistO n(M)\in \AG((\Z/\ell)^n\x G, (\Z/\ell)^n\x H)$ to be the matrix with entries as in Proposition \ref{propGroupTwistedLoopFunctor} except leaving out all summands where $\tup a$ is not in $R$, that is
\begin{align*}
\untwistO n(M)_{\tup a,\tup b} ={}& \hspace{-.4cm} \sum_{\substack{(R,\ph) \text{ s.t. }\tup a\in R\text{ and} \\ \ph(\tup a^{k(\tup a,R)})\text{ is $H$-conj. to $\tup b$}}}\hspace{-.4cm} c_{R,\ph}\cdot 
\Bigl( [\ev_{\tup a}^{-1}(R),(\id_{(\Z/\ell)^n}\x\ph)\circ \wind(\tup a, R)]_{(\Z/\ell)^n\x C_G(\tup a)}^{(\Z/\ell)^n\x C_H(\ph(\tup a^{k(\tup a,R)}))} 
\\ &\cmp ((\Z/\ell)^n\x \change_{\ph(\tup a^{k(\tup a,R)})}^{\tup b}) \Bigr)
\\ ={}& \hspace{-.4cm} \sum_{\substack{(R,\ph) \text{ s.t. } \tup a\in R\text{ and}\\ \ph(\tup a)\text{ is $H$-conj. to $\tup b$}}}\hspace{-.4cm} c_{R,\ph}\cdot 
\Bigl( [(\Z/\ell)^n\x R,(\id_{(\Z/\ell)^n}\x\ph)]_{(\Z/\ell)^n\x C_G(\tup a)}^{(\Z/\ell)^n\x C_H(\ph(\tup a))} 
\\ &\cmp ((\Z/\ell)^n\x \change_{\ph(\tup a)}^{\tup b}) \Bigr).
\end{align*}
Comparing with Corollary \ref{corSimpleGroupLoopOrbits}, we see that $\untwistO n(M)$ coincides with $(\Z/\ell)^n\x \freeO n(M)$. We observed this for the category $\Cov$ in Remark \ref{remarkCovUntwisted}.
\end{remark}

 We proceed to state Theorem \ref{thmCovMain} in terms of the category $\AG$, recalling Convention \ref{conventionOrderOfCyclics}, with only minor alterations:
\begin{theorem}\label{thmGroupLiftingFunctor}\label{thmGroupMain}
The endofunctors $\twistO n \colon \AG\to \AG$ for $n\geq 0$ of Proposition \ref{propGroupTwistedLoopFunctor} have the following properties:
\begin{enumerate}
\renewcommand{\theenumi}{$(\roman{enumi})$}\renewcommand{\labelenumi}{\theenumi}
\item[$(\emptyset )$] Let $L^{\dagger,\Cov}_n\colon \Cov \to \Cov$ be the functor constructed in Section \ref{secCoveringMaps}. The functor $\twistO n\colon \AG \to \AG$ is the linearization of the restriction of $\Ho(L^{\dagger,\Cov}_n)$ to the full subcategory of $\Ho(\Cov)$ spanned by finite unions of classifying spaces of finite groups.
\item\label{itemGroupLnZero} $\twistO 0$ is the identity functor on $\AG$.
\item\label{itemGroupLnObjects} On objects, $\twistO n$ takes a formal union of groups $G$ to the formal union of groups
\[\twistO n(G) = (\Z/\ell)^n\x \freeO n (G)=\coprod_{\tup a\in \cntuples nG} (\Z/\ell)^n\x C_G(\tup a).\]
\item\label{itemGroupEquivariant} The group $\Sigma_n$ acts on $\freeO nG=\coprod_{\tup a\in \cntuples nG} C_G(\tup a)$ by permuting the coordinates of the $n$-tuples $\tup a$. Explicitly, if $\sigma\in \Sigma_n$ and if $\widetilde{\sigma(\tup a)}$ is the representative for the $G$-conjugacy class of $\sigma(\tup a)$, then $\sigma\colon \freeO nG\to \freeO nG$ maps $C_G(\tup a)=C_G(\sigma(\tup a))$ to $C_G(\widetilde{\sigma(\tup a)})$ via the isomorphism $\change_{\sigma(\tup a)}^{\widetilde{\sigma(\tup a)}}\in \AG(C_G(\tup a),C_G(\widetilde{\sigma(\tup a)}))$. 

The functor $\twistO n$ is equivariant with respect to the $\Sigma_n$-action on $(\Z/\ell)^n\x \freeO n(-)$ that permutes the coordinates of both $(\Z/\ell)^n$ and $\freeO n(-)$, i.e. for every $\sigma\in\Sigma_n$ the diagonal action of $\sigma$ on $(\Z/\ell)^n\x \freeO n(-)$ induces a natural isomorphism $\sigma\colon \twistO n\underset\cong\Rightarrow \twistO n$.
\item\label{itemGroupLnForwardMaps} On forward maps, i.e. transitive bisets $[G,\ph]_G^H\in \AG(G,H)$ with $\ph\colon G\to H$, the functor $\twistO n$ coincides with $(\Z/\ell)^n \x \freeO n(-)$ so that 
\[\twistO n([G,\ph]_G^H)=(\Z/\ell)^n\x \freeO n([G,\ph]_G^H)=(\Z/\ell)^n\x \freeO n(\ph).\]
In particular, $\freeO n([G,\ph]_G^H)$ is the biset matrix that takes a component $C_G(\tup a)$ of $\freeO nG$ to the component $C_H(\tup b)$ by the map $c_h\circ\ph\colon C_G(\tup a)\to C_H(\ph(\tup a))$, where $\tup b$ represents the $H$-conjugacy class of $\ph(\tup a)$ and $h\in H$ conjugates $\ph(\tup a)$ to $\tup b$.
\item\label{itemGroupEvalSquare} For all $n \geq 0$, the functor $\twistO n$ commutes with evaluation maps, i.e. the evaluation maps $\ev_G\colon (\Z/\ell)^n\x \freeO n(G)\to G$ form a natural transformation $\ev\colon \twistO n \Rightarrow \Id_{\AG}$.
\item\label{itemGroupLnPartialEvaluation} For all $n \geq 0$, the partial evaluation maps $\pev_G\colon \Z/\ell\x \freeO {n+1}(G) \to \freeO n (G)$ given by 
\[
\pev(t,z)= (a_{n+1})^t\cdot z\in C_G(a_1,\dotsc,a_n), \quad \text{for } t\in \Z/\ell,  z\in \coprod_{\tup a\in \cntuples{n+1}G} C_G(\tup a),
\] form natural transformations $(\Z/\ell)^n\x \pev\colon \twistO {n+1} \Rightarrow \twistO n$.
\item\label{itemGroupIterateLn} For all $n,m\geq 0$, and any formal union of groups $G$, the formal union $(\Z/\ell)^{n+m}\x \freeO {n+m} G$ embeds into $(\Z/\ell)^m\x \freeO m((\Z/\ell)^n\x \freeO n G)$ as the components corresponding to the commuting $m$-tuples in $(\Z/\ell)^n\x \freeO n G$ that are zero in the $(\Z/\ell)^n$-coordinate, i.e. the embedding is given by
\begin{multline*}
    \bigl((\tup s,\tup r), z\bigr)\in (\Z/\ell)^{n+m}\x C_G((\tup x,\tup y))
    \\\mapsto \bigl(\tup r, (\tup s,z) \bigr) \in (\Z/\ell)^m\x C_{(\Z/\ell)^n\x C_G(\tup x)}(\tup 0\x \tup y)
\end{multline*}
for $\tup s\in (\Z/\ell)^n$, $\tup r\in (\Z/\ell)^m$, $\tup x\in \cntuples n G$, $\tup y\in \cntuples m G$, and $z\in C_G(\tup x,\tup y)$, and where $C_{(\Z/\ell)^n\x C_G(\tup x)}(\tup 0\x \tup y)$ is a component of 
\[\freeO m((\Z/\ell)^n\x C_G(\tup x)) \subseteq \freeO m((\Z/\ell)^n\x \freeO n G).\]

These embeddings $(\Z/\ell)^{n+m}\x \freeO {n+m} G\to (\Z/\ell)^m\x \freeO m((\Z/\ell)^n\x \freeO n G)$ then form a natural transformation $\twistO {n+m}(-)\Rightarrow \twistO m(\twistO n(-))$.
\end{enumerate}
\end{theorem}

\begin{proof}
\ref{itemGroupLnZero}-\ref{itemGroupLnObjects}: Immediate from Theorem \ref{thmCovMain}\ref{itemCovLnZero}-\ref{itemCovLnObjects} and our choice of modelling $S^1\x \freeO nBG$ by the formal union of groups $(\Z/\ell)^n\x \freeO nG$ as in Lemma \ref{lemmaGroupsAlgModelForFreeLoops}.

\ref{itemGroupEquivariant}: Given a permutation $\sigma\in \Sigma_n$ and the representative $\tup a$ of some $G$-conjugacy class of commuting $n$-tuples, we need to be careful that the permuted tuple $\sigma(\tup a)$ might not be the chosen representative of its conjugacy class. Permuting the coordinates of $\freeO n BG$ just corresponds to the identity map between centralizers $C_G(\tup a)=C_G(\sigma(\tup a))$, but then we need to exchange $\sigma(\tup a)$ with its representative $\widetilde{\sigma(\tup a)}$, hence the need for the isomorphism $\change_{\sigma(\tup a)}^{\widetilde{\sigma(\tup a)}}$ between $C_G(\tup a)=C_G(\sigma(\tup a))$ and $C_G(\widetilde{\sigma(\tup a)})$ in $\AG$.

Now that we have established the action of $\Sigma_n$ on $\freeO n G$ as above, the invariance of $\twistO n(-)$ follows from Theorem \ref{thmCovMain}\ref{itemCovEquivariant}.

\ref{itemGroupLnForwardMaps}: The transitive biset $[G,\ph]_G^H\in \AG(G,H)$ corresponds to the span $BG\xffrom\id BG \xto{B\ph} BH$ so Theorem \ref{thmCovMain}\ref{itemCovLnForwardMaps} implies that
\[\twistO n([G,\ph]_G^H)=(\Z/\ell)^n\x \freeO n([G,\ph]_G^H).\]

\ref{itemGroupEvalSquare}: Follows from Theorem \ref{thmCovMain}\ref{itemCovEvalSquare}.

\ref{itemGroupLnPartialEvaluation}: The partial evaluation map $\pev_X\colon S^1\x \freeO{n+1} X\to \freeO nX$ in Theorem \ref{thmCovMain}\ref{itemCovLnPartialEvaluation} evaluates an $(n+1)$-fold loop at its last coordinate. Comparing with the algebraic evaluation map $\ev_G\colon (\Z/\ell)^n\x \freeO nG\to G$ of Lemma \ref{lemmaGroupEvaluation}, the partial evaluation map $\pev_G\colon \Z/\ell\x \freeO{n+1} G \to \freeO n G$ must take the component $\Z/\ell \x C_G(\tup a)$ corresponding to an $(n+1)$-tuple $\tup a$ and send it to the component $C_G(a_1,\dotsc,a_n)$ of $\freeO nG$ by taking the product with the $t$th power of the last coordinate of $\tup a$ for $t \in \Z/\ell$:
\[\pev_{\tup a}(t,z) = (a_{n+1})^t\cdot z\in C_G(a_1,\dotsc,a_n)\]
for $t\in\Z/\ell$ and $z\in C_G(\tup a)$.
That we get a natural transformation this way is due to Theorem \ref{thmCovMain}\ref{itemCovLnPartialEvaluation}.

\ref{itemGroupIterateLn}: This follows from Theorem \ref{thmCovMain}\ref{itemCovIterateLn} once we translate the embedding $(S^1)^{n+m}\x \freeO{n+m}(BG)\to (S^1)^m\x \freeO m((S^1)^n\x \freeO n(BG))$ to the algebraic setting. The topological embedding has the formula
\begin{multline*}((\tup s,\tup r), f)\in (S^1)^{n+m}\x\freeO{n+m}(BG)
\\ \mapsto \Bigl(\tup r, \tup r'\mapsto \bigl(\tup s, \tup s'\mapsto f(s',r')\bigr)\Bigr)\in (S^1)^m\x \freeO m((S^1)^n\x \freeO n(BG)),\end{multline*}
for $\tup s,\tup s'\in (S^1)^n$, $\tup r,\tup r'\in (S^1)^m$ and $f\in \freeO {n+m} X$.

Consider a component $(\Z/\ell)^{n+m}\x C_G(\tup x,\tup y)$ in $(\Z/\ell)^{n+m}\x \freeO{n+m}G$. The embedding has to land in the component corresponding to some $m$-tuple $\tup w$ of $(\Z/\ell)^n\x \freeO nG$ that models the topological map $(S^1)^m\to (S^1)^n\x \freeO n(BG)$. The $m$-tuple corresponds to an $m$-fold loop of $(S^1)^n\x \freeO n(BG)$ that is constant in the $(S^1)^n$-coordinate, thus the $(\Z/\ell)^n$-coordinate of $\tup w$ is just the zero-tuple $\tup 0\in (\Z/\ell)^n$. The $\freeO n(G)$-coordinate of the $m$-tuple $\tup w$ corresponds to the topological map $(S^1)^m\to \freeO n(BG)$ given by
\[\tup r'\mapsto (\tup s' \mapsto f(s',r')).\]
This is an $m$-fold loop in the component of $\tup s'\mapsto f(\tup s', \tup 0)$, so lies in the component of $\freeO n(BG)$ corresponding to the first $n$ coordinates of $f\in \freeO{n+m}(BG)$. Consequently, the $\freeO n(G)$-coordinate of $\tup w$ lies in the component of $\freeO n(G)$ corresponding to $\tup x$, i.e. the centralizer $C_G(\tup x)$. Now the $\freeO n(G)$-coordinate of $\tup w$ is a commuting $m$-tuple in $C_G(\tup x)$, and considering how $(S^1)^m\to \freeO n(BG)$ moves the basepoint $\freeO n(BG)\xto{\gamma\mapsto \gamma(\tup 0)} BG$ as an $m$-fold loop in $BG$ matching the last $m$-coordinates of $f\in \freeO{n+m}(BG)$, the $\freeO n(G)$-coordinate of $\tup w$ must be $\tup y$ as an $m$-tuple in $C_G(\tup x)$.

In total $\tup w= \tup 0\x \tup y$ as a commuting $m$-tuple in $(\Z/\ell)^n\x C_G(\tup x)$. The embedding therefore takes the component $(\Z/\ell)^{n+m}\x C_G(\tup x,\tup y)$ in $(\Z/\ell)^{n+m}\x \freeO{n+m}G$ to the component $(\Z/\ell)^m\x C_{(\Z/\ell)^n\x C_G(\tup x)}(\tup 0\x \tup y)$.

The map $(\Z/\ell)^{n+m}\x C_G(\tup x,\tup y)\to (\Z/\ell)^m$ is just the projection $((\tup s,\tup r),z)\mapsto \tup r$. Furthermore, evaluating the basepoint of $\freeO m(-)$ topologically as for $\freeO n(-)$ above, we can see that $(\Z/\ell)^{n+m}\x C_G(\tup x,\tup y)\to C_{(\Z/\ell)^n\x C_G(\tup x)}(\tup 0\x \tup y)$ is the map $((\tup s,\tup r),z)\mapsto (\tup s,z)\in C_{(\Z/\ell)^n\x C_G(\tup x)}(\tup 0\x \tup y)$.

In total, $(\Z/\ell)^{n+m}\x C_G(\tup x,\tup y)$ is mapped to $(\Z/\ell)^m\x C_{(\Z/\ell)^n\x C_G(\tup x)}(\tup 0\x \tup y)$ by the formula $((\tup s,\tup r),z)\mapsto (\tup r,(\tup s, z))$.
\end{proof}

\begin{remark}
The above proof depends on Theorem \ref{thmCovMain}. Of course, it is possible to give a purely algebraic proof of Theorem \ref{thmGroupMain} instead.  
\end{remark}

Suppose $\cC$ is a subclass of all abelian groups such that $\cC$ is closed under isomorphism, subgroups, and quotients. We then denote by $\freeO n_{\cC}(G)$ the union of those components $C_G(\tup a)$ in $\freeO n (G)$ where the commuting $n$-tuple $\tup a$ generates a subgroup $\gen{\tup a}$ in $\cC$.

The constructions in this chapter always map a tuple $\tup a$ to tuples $\tup b$ where we take powers of elements in $\tup a$ and map them by homomorphisms. Consequently, whenever $\gen{\tup a}\in \cC$ and $\gen{\tup b}\not\in \cC$, then $\freeO n(X)_{\tup a,\tup b}= 0$ and $\twistO n(X)_{\tup a,\tup b}=0$ for all virtual bisets $X$. This implies the following:
\begin{prop}\label{propSubclassOfGroups}
Suppose $\cC$ is a subclass of all abelian groups, such that $\cC$ is closed under isomorphism, subgroups, and quotients. Furthermore suppose $\ell$ is big enough that every $\tup a\in \ntuples nG$ with $\gen a\in \cC$ has all element orders dividing $\ell$.

We then have functors $\freeO n_{\cC} \colon \AG \to \AG$ and $\twistO {n,\cC}\colon \AG\to \AG$ analogous to $\freeO n$ and $\twistO n$, respectively.
On objects we have \[\twistO {n,\cC}(G) = (\Z/\ell)^n\x \freeO n_{\cC}(G),\]
and $\twistO {n,\cC}$ satisfies a list of properties analogous to Theorem \ref{thmGroupMain}.
\end{prop}

\begin{remark}
We do not need $(\Z/\ell)^n$ to be a member of the class $\cC$ as long as it projects onto all  $\gen{\tup a}$ for tuples $\tup a \in \ntuples nG$ with $\gen{\tup a} \in \cC$. The functor $\freeO n_{\cC}(G)$ simply restricts attention to those homomorphisms $(\Z/\ell)^n \to G$ that factor through quotients of $(\Z/\ell)^n$ lying in $\cC$.
\end{remark}

\begin{proof}
We prove this in two steps. 

First assume that $\ell$ is large enough that any map $\Z \to G$ factors through $\Z/\ell$. We define $\tilde{\twistO {n,\cC}}$ by restricting $\twistO n$ to the subset of components $ (\Z/\ell)^n\x \freeO n_{\cC}(G) \subseteq  (\Z/\ell)^n\x \freeO n G$. On morphisms $M\in \AG(G,H)$,  the matrix $\tilde{\twistO {n,\cC}}(M)\in \AG((\Z/\ell)^n\x \freeO n_{\cC}(G), (\Z/\ell)^n\x \freeO n_{\cC}(H))$ is well-defined since $\twistO n(M)_{\tup a,\tup b} = 0$ whenever $\gen{\tup a}\in \cC$ and $\gen{\tup b}\not\in \cC$.

Now let $\ell'$ be a divisor of $\ell$ satisfying that every map $\Z^n \to A\subseteq G$ with $A\in \cC$ factors through $(\Z/\ell')^n$.
The evaluation map $(\Z/\ell)^n\times \freeO n_{\cC}(G)\to G$ factors via the projection $\Z/\ell\to \Z/\ell'$ as
\[
(\Z/\ell)^n\times \freeO n_{\cC}(G)\to (\Z/\ell')^n\times \freeO n_{\cC}(G)\xto{\ev} G.
\]
Given $M\in \AG(G,H)$ we will fill in the diagram
\[
\begin{tikzpicture}
\node (M) [matrix of math nodes] {
    (\Z/\ell)^n\times \freeO n_{\cC}(G) &[2cm] (\Z/\ell')^n\times \freeO n_{\cC}(G) &[2cm] G \\[2cm]
    (\Z/\ell)^n\times \freeO n_{\cC}(H) &[2cm] (\Z/\ell')^n\times \freeO n_{\cC}(H) &[2cm] H. \\
};
\path[arrow, auto, ->] 
    (M-1-1) edge (M-1-2)
            edge node[swap]{$\tilde{\twistO {n,\cC}}(M)$} (M-2-1)
    (M-1-2) edge[dashed] node[swap]{$\twistO{n,\cC}(M)$} (M-2-2)
            edge node{$\ev$} (M-1-3)
    (M-1-3) edge node[swap]{$M$} (M-2-3)
    (M-2-1) edge (M-2-2)
    (M-2-2) edge node{$\ev$} (M-2-3)
;
\end{tikzpicture}
\]
Being a restriction of $\twistO n(M)$, the matrix $\tilde{\twistO {n,\cC}}(M)$ satisfies the formula of Proposition \ref{propGroupTwistedLoopFunctor}. We can simply define $\twistO{n,\cC}(M)$ according to the same formula, replacing each occurrence of $\ell$ with $\ell'$. The formula for $\wind(\tup a, R)$ does not involve $\ell$ or $\ell'$ at all, and so the diagram commutes.

The surjective group homomorphism $(\Z/\ell)^n\times \freeO n_{\cC}(G) \to (\Z/\ell')^n\times \freeO n_{\cC}(G)$ is an epimorphism in $\AG$, and so $\twistO{n,\cC}(M)$ is the unique map making the diagram commute. From this uniqueness it follows that $\twistO{n,\cC}$ is functorial and satisfies the properties of Theorem \ref{thmGroupMain}.
\end{proof}

Given a prime $p$, a particularly relevant collection $\cC$ is the collection of all finite abelian $p$-groups. If $G$ is a finite group, then we write $\freeO n_p(G) = \freeO n_{\cC} (G)$ for the free loops with respect to the collection of finite abelian $p$-groups. This consists of all components of $\freeO n (G)$ corresponding to commuting $n$-tuples of $p$-power order elements. When dealing exclusively with $p$-power order elements we only consider maps $\Z/\ell\to G$ with images of $p$-power order, hence we may suppose that $\ell$ equals a sufficiently high power $p$. As such, for $e$ large enough, we have
\[
\twistO {n,p}(G) = (\Z/p^e)^n\x \freeO n_p(G)
\]
on objects in $\AG$.

\subsection{$\twistO n$ does not extend to all classifying spectra}

Let $\Sp$ be the category of spectra. Since the homotopy classes of maps between two spectra form an abelian group, composing the suspension spectrum functor $\Sigma_+^\infty\colon \Ho(\Cov)\to \Ho(\Sp)$ and $B \colon \AG_+ \to \Ho(\Cov)$, provides us with a functor $\AG\to \Ho(\Sp)$. A finite group $G$ in $\AG$ is sent to the classifying spectrum $\Sigma_+^\infty BG$.

The Burnside ring $A(G)$ of finite $G$-sets acts on $\AG(G,H)$ by cartesian products, or equivalently $G/R\in A(G)$ acts on $\AG(G,H)$ by precomposing with $\incl_R^G\circ\tr_R^G=[R,\id]_G^G\in \AG(G,G)$. A variant of the Segal conjecture \cite{LewisMayMcClure} states that on morphism sets the functor $\AG\to \Ho(\Sp)$ completes $\AG(G,H)$ at the augmentation ideal of the Burnside ring $A(G)$.

\begin{theorem} [\cite{CarlssonSegal, AGM, LewisMayMcClure}]
Let $G$ and $H$ be finite groups. The map $\AG(G,H)\to [\Sigma_+^\infty BG, \Sigma_+^\infty BH]$ taking $[R,\ph]_G^H$ to $\Sigma^\infty_+B\ph \circ \tr_R^G$ induces an equivalence
\[ [\Sigma_+^\infty BG, \Sigma_+^\infty BH] \cong \AG(G,H)^\wedge_{I(G),}\]
where $I(G)=\ker(A(G)\xto{\abs{-}} \Z)$ is the augmentation ideal of $A(G)$.
\end{theorem}

We conclude Section \ref{secGroups} with a concrete calculation of $\twistO 1$ applied to a particular biset for the group $C_2$. This example will demonstrate why it is impossible to extend $\twistO 1\colon \AG \to \AG$ to a functor on the homotopy category of the category of classifying spectra of finite groups. The sequel \cite{RSS_Bold2} will remedy this by completing at a prime $p$ and using the theory of fusion systems.

\begin{example}\label{exampleDoesNotPassToSpectra} 
Let $C_2$ denote the cyclic group of order two consisting of the neutral element $\iota$ and the involution $\tau$. We will denote the trivial subgroup of $C_2$ by $1$ so that $1 = \{\iota\}$.

Throughout this example we shall consider the $(C_2,C_2)$-biset 
\[X = [1,\id]_{C_2}^{C_2} - 2[C_2,\id]_{C_2}^{C_2} = C_2\x_1 C_2 - 2( C_2)\in \AG(C_2,C_2).\]
Recall that the Burnside ring $A(C_2)$ acts on $\AG(C_2,C_2)$ via the ring homomorphism $A(C_2)\to \AG(C_2,C_2)$ taking a left $C_2$-set $Y\in A(C_2)$ to $Y\x C_2\in \AG(C_2,C_2)$, where $C_2$ acts on $Y\x C_2$ diagonally on the left. The biset $X$ is the image of $Y=C_2 - 2(C_2/C_2)$ under this ring homomorphism.

The virtual biset $Y=C_2 - 2(C_2/C_2)$ has augmentation $0$, hence $Y$ lies in the augmentation ideal $I_{C_2}\leq A(C_2)$ and $X\in I_{C_2}\cdot \AG(C_2,C_2)$.

Let $n=1$ and consider $\twistO 1(X)\in \AG(\Z/2^e\x \freeO 1 (C_2), \Z/2^e\x \freeO 1 (C_2))$. To distinguish the two components of $\freeO 1(C_2)$ we shall denote them $C_\iota= C_{C_2}(\iota)$ the centralizer of the neutral element, and $C_\tau = C_{C_2}(\tau)$ the centralizer of the involution.

The biset $[C_2,\id]_{C_2}^{C_2}$ is the identity morphism in $\AG(C_2,C_2)$. Since $\twistO 1$ is a functor, it follows that $\twistO 1([C_2,\id]_{C_2}^{C_2})$ is the identity matrix 
\[
\twistO 1([C_2,\id]_{C_2}^{C_2}) = \bordermatrix{
&\iota & \tau \cr
\iota & [\Z/2^e \x C_\iota, \id] & 0 \cr
\tau & 0 & [\Z/2^e \x C_\tau, \id]
}.
\]
The biset $[1,\id]_{C_2}^{C_2}$ decomposes as the transfer $[1,\id]_{C_2}^1$ followed by the inclusion $[1,\id]_1^{C_2}$, hence by functoriality $\twistO 1([1,\id]_{C_2}^{C_2})$ decomposes as well.
If we apply $\twistO 1$ to the inclusion $[1,\id]_1^{C_2}$, Theorem \ref{thmGroupMain}\ref{itemGroupLnForwardMaps} tells us that $\twistO 1([1,\id]_1^{C_2})$ simply maps $C_1(\iota)=1$ into $C_{C_2}(\iota) = C_\iota$. That is, we get the matrix
\[
\twistO 1([1,\id]_1^{C_2}) = \bordermatrix{
&\iota & \tau \cr
\iota & [\Z/2^e \x 1, \id]_{\Z/2^e \x 1}^{\Z/2^e \x C_\iota} & 0 
}.
\]
Lastly, we cannot determine $\twistO 1([1,\id]_{C_2}^{1})$ from the properties in Theorem \ref{thmGroupMain} alone. We instead have to refer to Proposition \ref{propGroupTwistedLoopFunctor}, which simplifies since $C_\iota = C_\tau = C_2$:
\[
\twistO 1([1,\id]_{C_2}^1) = \bordermatrix{
&\iota  \cr
\iota & [\ev_\iota^{-1}(1), \wind(\iota, 1)]_{\Z/2^e \x C_\iota}^{\Z/2^e \x 1}   \cr
\tau & [\ev_\tau^{-1}(1), \wind(\tau, 1)]_{\Z/2^e \x C_\tau}^{\Z/2^e \x 1} 
} = 
\bordermatrix{
&\iota  \cr
\iota & [\Z/2^e \x 1, \id]_{\Z/2^e \x C_\iota}^{\Z/2^e \x 1}   \cr
\tau & [\ev_\tau^{-1}(1), \wind(\tau, 1)]_{\Z/2^e \x C_\tau}^{\Z/2^e \x 1} 
}.
\]
Composing the two previous matrices, we get
\[
\twistO 1([1,\id]_{C_2}^{C_2}) = \bordermatrix{
&\iota & \tau \cr
\iota & [\Z/2^e \x 1, \id]_{\Z/2^e \x C_\iota}^{\Z/2^e \x C_\iota} & 0 \cr
\tau  & [\ev_\tau^{-1}(1), \wind(\tau, 1)]_{\Z/2^e \x C_\tau}^{\Z/2^e \x C_\iota}  & 0
}.
\]
Finally, we substract the identity matrix twice to calculate $\twistO 1(X)$:
\begin{align*}
\twistO 1(X) &= \twistO 1([1,\id]_{C_2}^{C_2}-2[C_2,\id]_{C_2}^{C_2}) \\ & = \bordermatrix{
&\iota & \tau \cr
\iota & [\Z/2^e \x 1, \id] - 2[\Z/2^e \x C_\iota, \id] & 0 \cr
\tau  & [\ev_\tau^{-1}(1), \wind(\tau, 1)]  & -2[\Z/2^e \x C_\iota, \id]
}.
\end{align*}
Our original biset $X$ lies in the submodule $I_{C_2}\AG(C_2,C_2)$ and is the image of $Y \in I_{C_2} \subset A(C_2)$, so if we look at the powers of $X$, the sequence $X, X^2, X^3, \dotsc$ converges to zero in the completion $\AG(C_2,C_2)^\wedge_{I_{C_2}}\cong [\Sinfp BC_2, \Sinfp BC_2]$.
Hence if we want $\twistO 1\colon \AG\to \AG$ to pass to a functor on classifying spectra and homotopy classes of maps, we need the sequence $\twistO 1(X), \twistO 1(X)^2, \twistO 1(X)^3,\dotsc$ to converge to $0$ as well. If the entries of $\twistO 1(X)$ (or some power of $\twistO 1(X)$) each lie in the respective submodules generated by the augmentation ideals, this would be true.

We have augmentation homomorphisms $\e\colon \AG(G,H)\to \Z$ preserving composition for any finite groups given by $\e(Z) = \abs Z / \abs H$. The elements of the submodule $I_G \AG(G,H)$ are mapped to $0$ by $\e$. Applying $\e$ entrywise to $\twistO 1(X)$ we get the integer matrix
\begin{align*}
\e(\twistO 1(X)) &=  \bordermatrix{
&\iota & \tau \cr
\iota & \e([\Z/2^e \x 1, \id]) - 2\e([\Z/2^e \x C_\iota, \id]) & 0 \cr
\tau  & \e([\ev_\tau^{-1}(1), \wind(\tau, 1)])  & -2\e([\Z/2^e \x C_\iota, \id])
} \\&=  \bordermatrix{
&\iota & \tau \cr
\iota & 0 & 0 \cr
\tau  & 2  & -2
}.
\end{align*}
The powers of the matrix are 
\begin{align*}
\e(\twistO 1(X))^m &=  \bordermatrix{
&\iota & \tau \cr
\iota & 0 & 0 \cr
\tau  & -(-2)^m  & (-2)^m
}.
\end{align*}

From this calculation we conclude two things. First, the entries of the $\tau$-row in $\twistO 1(X)^m$ are never going to be in the submodule generated by the augmentation ideal, so $\twistO 1(X)^m$ cannot converge to the $0$-matrix in the $I$-adic completions, and consequently $\twistO 1$ does not induce a functor on classifying spectra of finite groups in general. Second, $\twistO 1(X)^m$ does in fact converge to the $0$-matrix in the $2$-adic topologies! Hence if we follow $\twistO 1$ by $2$-completion, there is hope that the composite
\[
\AG \xto{\twistO 1} \AG \to \Ho(\Sp) \xto{(-)^\wedge_2} \Ho(\Sp_2)
\]
might pass to a functor from classifying spectra of finite groups to $2$-complete classifying spectra. 
\end{example}

\begin{remark}\label{remarkDoesNotPassToSpectra}
One might wonder if another definition of $\twistO 1$ might pass to $I$-adic completions. However, in Example \ref{exampleDoesNotPassToSpectra} we mainly use property \ref{itemGroupLnForwardMaps} of Theorem \ref{thmGroupMain} -- except to calculate the transfer $\twistO 1([1,\id]_{C_2}^{1})$.

Suppose $L\colon \AG\to \AG$ is a functor satisfying properties \ref{itemGroupLnObjects} and \ref{itemGroupLnForwardMaps}, and let 
\[
L([1,\id]_{C_2}^1) =  
\bordermatrix{
&\iota  \cr
\iota & A \in \AG(\Z/2^e \x C_\iota,\Z/2^e \x 1)   \cr
\tau & B \in \AG(\Z/2^e \x C_\tau,\Z/2^e \x 1) 
}
\]
be unspecified.
Then we can still calculate the augmentation 
\begin{align*}
\e(L(X)) &=  \bordermatrix{
&\iota & \tau \cr
\iota & \abs{A} - 2 & 0 \cr
\tau  & \abs{B}  & -2
}.
\end{align*}
Because of the second diagonal entry, the powers of $L(X)$ still do not converge to the $0$-matrix in the $I$-adic topologies, so $L$ cannot pass to a functor on classifying spectra and homotopy classes of maps between them.
\end{remark}

Fix a prime $p$. It is easy to extend the functor $\twistO n$ from the full subcategory of finite $p$-groups in $\AG$ to the homotopy category of $p$-completed classifying spectra of finite $p$-groups. For a space $X$, we will write $\Sinfpp X = (\Sinfp X)^{\wedge}_{p}$ for the $p$-completion of the suspension spectrum of $X_{+}$. From \cite[Lemma 2.4 and Proposition 2.18]{RSS_p-completion}, it follows that if $R$ and $S$ are $p$-groups, then the Segal conjecture implies that 
\[
[\Sinfpp BR, \Sinfpp BS] \cong \AG(R,S)^\wedge_p \cong \AG(R,S) \otimes \Z_p.
\]
Thus we may restrict $\twistO n$ to the full subcategory of finite $p$-groups and base change to $\Z_p$ (as $\twistO n$ is additive) to extend to the homotopy category of $p$-completed classifying spectra of finite $p$-groups. In the sequel \cite{RSS_Bold2}, we will use the theory of fusion systems to extend this further to the homotopy category of $p$-completed classifying spectra of finite groups.

%% file: covering_appendix_I.tex
\section{Mapping spaces for finite sheeted coverings}\label{appendixCovering}
In this appendix we prove that the free loop space of a finite cover of topological spaces is again a finite cover, with no assumptions about the spaces involved. Further, we show that this result remains true in the category of compactly generated weakly Hausdorff topological spaces when we replace the compact-open topology on the free loop space with its compactly generated refinement.

Given a finite sheeted covering map $p\colon E\to B$ and topological space $X$, we can ask whether the induced map $\Map(X,E)\to\Map(X,B)$ is again a finite sheeted covering map.
In this appendix we provide conditions on $X$ that ensure that we always get a finite sheeted covering of mapping spaces with no further assumptions on the covering map $p\colon E\to B$ or the spaces involved. As a special case $X=(S^1)^n$ satisfies the conditions in question so that the $n$-fold free loop space functor $\freeO n(-)$ takes covering maps to covering maps (Corollary \ref{corFreeLoopSpaceCoveringMap}).

\begin{definition}
We say that an open covering $\{U_i\}_{i\in I}$ of a space $X$ is \emph{finite intersection connected} if every finite intersection 
\[
U_{i_1} \cap \dotsb \cap U_{i_n}
\]
is connected (possibly empty). In particular, each $U_i$ is itself connected.
\end{definition}

\begin{definition}
We say that a locally connected space $X$ is \emph{locally finite intersection connected} if every open covering of $X$ has a finite intersection connected open refinement.
\end{definition}

Examples of locally finite intersection connected spaces are the tori $(S^1)^n$ for $n\geq 0$. The circle $S^1\cong \R/\Z$ has a basis consisting of open intervals of small length, and the intersection of any finite number of open intervals is again an open interval or empty and thus connected. The same is true for the product basis in $(S^1)^n$.

\begin{lemma}\label{lemmaEvenlyCoveredFiniteCovering}
Let $p\colon E\to B$ be a finite sheeted covering map, and let $X$ be a locally finite intersection connected, compact, Hausdorff space. Suppose $f\colon X\to B$ is a continuous map, then $X$ has a finite, finite intersection connected, open covering $\{U_i\}_{i\in I}$ such that each $f(\overline{U_i})$ is contained in an evenly covered open set of $B$.
\end{lemma}

\begin{proof}
First note that $\{f^{-1}(V) \mid \text{$V\subseteq B$ is open and evenly covered}\}$ provides an open cover of $X$. Next, a compact Hausdorff space is regular, so each $x\in f^{-1}(V)$ has a small neighborhood $x\in U\subseteq \overline U \subseteq f^{-1}(V)$.
This gives us a second covering $\{U_i\}_{i\in I}$ of $X$ such that each $f(\overline{U_i})$ is contained in an evenly covered open set.
We then replace $\{U_i\}_{i\in I}$ by a finite intersection connected open refinement, and finally we take a finite subcover by the compactness of $X$.
\end{proof}

We now prove that $\Map(X,E)\to\Map(X,B)$ is a finite sheeted covering. The proof is based on the ideas of \cite{mappingspaces}*{Theorem 3.8}. The difference is that \cite{mappingspaces} requires a covering of $X$ by contractible open sets and only proves that $\Map(X,E)\to\Map(X,B)$ is a covering over each homotopy class of maps in $\Map(X,B)$, i.e. over each path-component.
\begin{prop}\label{propMappingSpaceCoveringMap}
Let $p\colon E\to B$ be a finite sheeted covering map, and let $X$ be a locally finite intersection connected, compact, Hausdorff space. The induced map $p_*\colon \Map(X,E)\to\Map(X,B)$ in the compact-open topology is then a finite sheeted covering map as well.
\end{prop}

\begin{proof}
The compact-open topology on $\Map(X,B)$ has a subbasis consisting of the subsets
\[
C_B(K,V) = \{f\in \Map(X,B) \mid f(K)\subseteq U\},
\]
where $K\subseteq X$ is compact and $V\subseteq B$ is open. Similarly, we have a subbasis for the topology on $\Map(X,E)$ given by subsets $C_E(K,V)$ with $K\subseteq X$ compact and $V\subseteq E$  open.

Let $f\in \Map(X,B)$ be given. We then have to produce an evenly covered neighborhood of $f$ in the compact-open topology. By Lemma \ref{lemmaEvenlyCoveredFiniteCovering} there is a covering of $X$ by finitely many open sets $U_1,\dotsc,U_n$ such that each $f(\overline{U_i})$ is contained in an evenly covered neighborhood $f(\overline{U_i})\subseteq V_i\subseteq B$.
Since $X$ is compact Hausdorff, each $\overline{U_i}$ is also compact. Hence $f$ is contained in the open intersection 
\[f\in C_B(\overline{U_1},V_1) \cap \dotsb \cap C_B(\overline{U_n},V_n),\]
though this is not yet the evenly covered neighborhood we seek.

Each $V_i$ is evenly covered by $p$, so we can choose a trivialization of $p^{-1}(V_i)\cong W_{i,1} \sqcup \dotsb \sqcup W_{i,r}$ in $E$ such that $p$ maps each $W_{i,k}$ homeomorphically onto $V_i$. The trivialization might not be unique, so suppose from now on that we have made a choice of trivialization for each $V_i$. Furthermore, note that the number of sheets $r$ might depend on which $V_i$ we a considering, but whenever $V_i$ and $V_j$ have non-empty intersection, they necessarily have the same number of sheets.

At the intersection of the open sets $V_{ij}:=V_i\cap V_j$, with $i\neq j$, we might run into the problem that the chosen sheets over $V_i$ and $V_j$ do not agree on the intersection so that some sheet $W_{i,k}$ over $V_i$ intersects several sheets over $V_j$. This can happen whenever the intersection $V_{ij}$ is not connected.
To account for this problem, we will partition $V_{ij}$ into a disjoint union of open subsets as follows: For every point $b\in V_{ij}$, each point in the fiber $p^{-1}(b)$ lies in a unique pair of sheets $W_{i,k}$ over $V_i$ and $W_{j,l}$ over $V_j$. Given $b\in V_{ij}$ we thus get an associated bijection $\sigma_b\in \Sigma_r$, where $r$ is the size of the fiber over $b$, such that each intersection $W_{i,k}\cap W_{j,\sigma_b(k)}$ contains a point of the fiber.
We now say that two points $b$ and $b'$ in $V_{ij}$ are equivalent if and only if they give rise to the same bijection $\sigma_b=\sigma_{b'}$, and let $V_{ij}^\sigma$ denote the subset of points giving rise to the bijection $\sigma\in \Sigma_r$.

The subset $V_{ij}^\sigma$ can be written as
\[
V_{ij}^\sigma = p(W_{i,1}\cap W_{j,\sigma(1)})\cap p(W_{i,2}\cap W_{j,\sigma(2)}) \cap \dotsb \cap p(W_{1,r}\cap W_{j,\sigma(r)}),
\]
which is a finite intersection of open sets and therefore open in $B$. The open sets $V_{ij}^\sigma$ for $\sigma\in \Sigma_r$ provides a partition of $V_{ij}$ into finitely many disjoint open sets, many of which might be empty.

Since $U_{ij}:= U_i\cap U_j\subseteq X$ is assumed to be connected, and supposing $U_{ij}$ is non-empty, $f$ maps $\overline {U_{ij}}$ to a unique set $V_{ij}^{\sigma_{ij}}$ in the partition of $V_{ij}$. We now claim that the following neighborhood of $f$ in $\Map(X,B)$ is evenly covered:
\[
D := \bigcap_{i=1}^n C_B(\overline{U_i}, V_i) \cap \bigcap_{\substack{1 \leq i < j\leq n\\ U_{ij}\neq \emptyset}} C_B(\overline{U_{ij}}, V_{ij}^{\sigma_{ij}}).
\]
The key property of the subsets involved in defining $D$ is that whenever a finite intersection of $U_i$'s is non-empty, the corresponding open set in $B$ is evenly-covered by a well-defined set of sheets:
\begin{itemize}
    \item $V_i$ is covered by the sheets $W_{i,1},\dotsc,W_{i,r}$.
    \item If $U_{ij} = U_i\cap U_j$ is non-empty, it must be mapped to $V_{ij}^{\sigma_{ij}}\subseteq B$, which is evenly covered by the sheets
    \[
        p^{-1}(V_{ij}^{\sigma_{ij}})\cap W_{i,1}\cap W_{j,\sigma_{ij}(1)}\ ,\ \dotsc\ ,\ p^{-1}(V_{ij}^{\sigma_{ij}})\cap W_{i,r}\cap W_{j,\sigma_{ij}(r)}.
    \]
    \item Similar statements hold for non-empty intersections of three or more $U_i$'s, but we will not need these in this proof.
\end{itemize}
Now suppose some $g\in D\subseteq \Map(X,B)$ lifts to a map $\tilde g\colon X\to E$ with $p\circ \tilde g = g$. By assumption each $U_i$ (and thus also $\overline{U_i}$) is connected, so $\tilde g(\overline{U_i})$ lies in a unique sheet $W_{i,s_i}$ of $p^{-1}(V_i)$ for some $1\leq s_i\leq r$, where $r$ is the number of sheets over $V_i$. Taking the corresponding sheet $W_{i,s_i}\subseteq E$ over $V_i$ for each $U_i$, we see that 
\[
\tilde g\in \bigcap_{i=1}^n C_E(\overline{U_i},W_{i,s_i}) \cap \bigcap_{\substack{1 \leq i < j\leq n\\ U_{ij}\neq \emptyset}} C_E(\overline{U_{ij}}, p^{-1}(V_{ij}^{\sigma_{ij}})).
\]
Furthermore, since $g$ maps any non-empty $\overline{U_{ij}}$ to $V_{ij}^{\sigma_{ij}}$, we must have $\tilde g(\overline{U_{ij}})$ contained in one of the sheets over $V_{ij}^{\sigma_{ij}}$, and consequently $s_j=\sigma_{ij}(s_i)$ whenever $U_{ij}=U_i\cap U_j$ is non-empty.

We now define a ``compatible choice of sheets'' to be a choice of a sheet $W_{i,s_i}$ over $V_i$, for each $1\leq i \leq n$, satisfying that whenever $U_i\cap U_j$ is non-empty, we must have $s_j=\sigma_{ij}(s_i)$.
Note that once we choose a sheet over $V_i$ for a single $U_i$ in each connected component of $X$, there is at most one compatible choice of sheets extending the choices made -- the remaining sheets are forced upon us by non-empty intersections.

By the argument above, any $\tilde g\in p_*^{-1}(D)$ which is a lift of some map $g\in D$, has to be contained in
\[\bigcap_{i=1}^n C_E(\overline{U_i},W_{i,s_i}) \cap \bigcap_{\substack{1 \leq i < j\leq n\\ U_{ij}\neq \emptyset}} C_E(\overline{U_{ij}}, p^{-1}(V_{ij}^{\sigma_{ij}}))\]
for some compatible choice of sheets $W_{1,s_1},\dotsc,W_{n,s_n}$. Conversely, given a compatible choice of sheets, the subset 
\[\bigcap_{i=1}^n C_E(\overline{U_i},W_{i,s_i}) \cap \bigcap_{\substack{1 \leq i < j\leq n\\ U_{ij}\neq \emptyset}} C_E(\overline{U_{ij}}, p^{-1}(V_{ij}^{\sigma_{ij}}))\]
is mapped to $D$ by $p_*\colon \Map(X,E)\to \Map(X,B)$. It follows that $p_*^{-1}(D)$ decomposes as a disjoint union
\[
p_*^{-1}(D) = \coprod_{\substack{\text{compatible}\\\text{choices of sheets}}}\!\left( \bigcap_{i=1}^n C_E(\overline{U_i},W_{i,s_i}) \cap \bigcap_{\substack{1 \leq i < j\leq n\\ U_{ij}\neq \emptyset}} C_E(\overline{U_{ij}}, p^{-1}(V_{ij}^{\sigma_{ij}}))\right).
\]
It remains to prove that for each compatible choice of sheets, the open subset
\[
D' := \bigcap_{i=1}^n C_E(\overline{U_i},W_{i,s_i}) \cap \bigcap_{\substack{1 \leq i < j\leq n\\ U_{ij}\neq \emptyset}} C_E(\overline{U_{ij}}, p^{-1}(V_{ij}^{\sigma_{ij}}))
\]
maps homeomorphically to $D$. 

For each $1\leq i \leq n$, let $\eta_i= p|_{W_{i,s_i}}\colon W_{i,s_i}\xto\cong V_i$ be the homeomorphism between the sheet $W_{i,s_i}$ and $V_i$. 
We then provide a point-set inverse to the map $p_*\colon D' \to D$ as follows:
Let $g\in D$, and let us define a map $\tilde g\colon X\to E$ by the formula
\[
\tilde g(x) =\begin{cases}
\eta_1^{-1}(g(x)) & \text{if }x\in U_1,
\\ \qquad\vdots
\\ \eta_n^{-1}(g(x)) & \text{if }x\in U_n.
\end{cases}
\]
We have defined $\tilde g$ by a sequence of continuous maps on an open cover of $X$, so $\tilde g$ is continuous if it is well-defined. 

To see that $\tilde g$ is well-defined, consider a non-empty intersection $U_{ij}=U_i\cap U_j$. Because $g\in D$, we know that $g(U_{ij})\subseteq V_{ij}^{\sigma_{ij}}$. Furthermore, since the the choice of sheets $W_{i,s_i}$ is assumed to be compatible, we know that $s_j=\sigma_{ij}(s_i)$ so that 
\[p^{-1}(V_{ij}^{\sigma_{ij}})\cap W_{i,s_i}\cap W_{j,s_j}\]
is one of the sheets over $V_{ij}^{\sigma_{ij}}$. Consequently, $\eta_i^{-1}$ and $\eta_j^{-1}$ land in the same sheet when restricted to $g(U_{ij})\subseteq V_{ij}^{\sigma_{ij}}$, hence $\eta_i^{-1}\circ g$ and $\eta_j^{-1}\circ g$ agree on $U_{ij}$ and $\tilde g\colon X\to E$ is well-defined.
By construction, we have $\tilde g\in D'$. This construction provides a point-set inverse to $p_*\colon D'\to D$. We conclude that $p_*\colon D'\to D$ is a continuous bijection. 

It remains to check that $p_*\colon D' \to D$ is a homeomorphism. Consider the following commutative square:
\[
\begin{tikzpicture}
\node (M) [matrix of math nodes] {
D' &[.5cm]\displaystyle \bigcap_{i=1}^n C_E(\overline{U_i},W_{i,s_i})  &[1cm]\displaystyle \prod_{i=1}^n \Map(\overline{U_i},W_{i,s_i})\\[1cm]
D &\displaystyle \bigcap_{i=1}^n C_B(\overline{U_i},V_i) &\displaystyle \prod_{i=1}^n \Map(\overline{U_i},V_i).\\
};
\path [auto,arrow,->]
    (M-1-1) edge node{$p_*$} (M-2-1)
    (M-1-2) edge[right hook->] (M-1-3)
            edge node{$p_*$} (M-2-2)
    (M-1-3) edge node{$\cong$} (M-2-3)
    (M-2-2) edge[right hook->] (M-2-3)
;
\path
    (M-1-1) -- node{$\subseteq$} (M-1-2)
    (M-2-1) -- node{$\subseteq$} (M-2-2)
;
\end{tikzpicture}
\]
The horizontal maps are continuous embeddings: We only prove that the top map is an embedding, the other proof is completely analogous. Injectivity of the map is clear as the open sets $U_i$ cover $X$. Given any subbasis element $C_E(K,O)$ with $K\subseteq X$ compact and $O\subseteq E$ open, the intersection \[
C_E(K,O)\cap \bigcap_{i=1}^n C_E(\overline{U_i},W_{i,s_i})
\]
is simply the preimage of
\[
\prod_{i=1}^n C_{\Map(\overline{U_i},W_{i,s_i})}(K\cap \overline{U_i}, O\cap W_{i,s_i}),
\]
where $C_{\Map(\overline{U_i},W_{i,s_i})}(K\cap \overline{U_i}, O\cap W_{i,s_i})$ denotes a subbasis element in the compact-open topology of $\Map(\overline{U_i},W_{i,s_i})$. Thus the topology on $\bigcap_{i=1}^n C_E(\overline{U_i},W_{i,s_i})$ agrees with the subspace topology coming from $\prod_{i=1}^n \Map(\overline{U_i},W_{i,s_i})$, and this remains true when restricting to the further open subset $D'$.

Because the composite $D'\xto{p_*} D \hookrightarrow \prod_{i=1}^n \Map(\overline{U_i},V_i)$ is an embedding, it follows that the initial map $p_*\colon D' \to D$ is an embedding as well. Finally, an embedding $p_*\colon D' \to D$ which is bijective, is in fact a homeomorphism.
This concludes our proof that $D$ is an evenly covered neighborhood of $f\in \Map(X,D)$ (with finitely many sheets), and therefore $p_*\colon \Map(X,E)\to \Map(X,B)$ is a finite sheeted covering map.
\end{proof}

\begin{cor}\label{corFreeLoopSpaceCoveringMap}
Let $p\colon E\to B$ be any finite sheeted covering map. Then the induced map $\freeO n(p)\colon \freeO n(E) \to \freeO n(B)$ is also a finite sheeted covering map for all $n\geq 0$, when $\freeO n(-)$ is endowed with the compact-open topology.
\end{cor}

\begin{proof}
Apply Proposition \ref{propMappingSpaceCoveringMap} with $X=(S^1)^n$, which is both locally finite intersection connected, compact, and Hausdorff.
\end{proof}

When working with the category of compactly generated weak Hausdorff spaces, we usually want to replace the compact-open topology on $\freeO n(-) = \Map(S^1, -)$ with the compactly generated weak Hausdorff replacement.
We may ask whether a finite sheeted covering $\freeO n(E) \to \freeO n(B)$ from Corollary \ref{corFreeLoopSpaceCoveringMap} is still a covering map after taking compactly generated weak Hausdorff replacements. This is a consequence of the following general lemma concerning finite sheeted coverings and refined topologies.
\begin{lemma}\label{lemmaInducedCoveringMaps}
Let $p\colon E\to B$ be a finite sheeted covering map, and let $F\colon \Top \to \Top$ be functor on all topological spaces such that:
\begin{itemize}
    \item $F(X)$ has the same underlying set as $X$, and the topology of $F(X)$ contains the topology of $X$.
    \item $X\xto {F(f)} Y$ equals $X\xto{f} Y$ as maps of sets.
    \item $F$ sends open embeddings $U\into X$ to open embeddings $F(U)\into F(X)$.
    \item $F(-)$ preserves finite disjoint unions of topological spaces.
\end{itemize}
Then $F(p)=p\colon F(E)\to F(B)$ is also a finite sheeted covering.
\end{lemma}

\begin{proof}
Any evenly covered open set $U$ in $B$ is still open in $F(B)$ and evenly covered as a set by the continuous map $F(p)^{-1}(U) \xto{F(p)|} U$. It is therefore sufficient to prove that $F(p)=p\colon F(E)\to F(B)$ is an open map, in which case $F(p)^{-1}(U) \xto{F(p)|} U$ will be a homeomorphism on each sheet over $U$.

To prove that $F(p)=p\colon F(E)\to F(B)$ is an open map, suppose $A\subseteq E$ is open in $F(E)$. We then aim to prove that $p(A)$ is open in $F(B)$. We cover $p(A)$ by evenly covered open sets in $B$, so that
\[p(A) =\bigcup_{U \text{ evenly covered in }B} p(A)\cap U.\]
It is now sufficient to prove that each $p(A)\cap U$ is open in $F(B)$. Consider the following diagram for an evenly covered $U$ in $B$:
\[
\begin{tikzpicture}
\node[matrix of math nodes] (M) {
\displaystyle\coprod_{j=1}^r V_j &[2cm] E \\[2cm]
U & B \\
};
\path[->, auto, arrow]
    (M-1-1) edge node{$\coprod_{j=1}^r \iota_j$} (M-1-2)
            edge node{$p$} node[swap]{$\cong\text{on sheets}$} (M-2-1)
    (M-1-2) edge node{$p$} (M-2-2)
    (M-2-1) edge (M-2-2)
;
\end{tikzpicture}
\]
If we apply $F$ to this diagram, we get the following square, where the left map is still a homeomorphism on sheets by the assumption that $F$ preserves coproducts:
\[
\begin{tikzpicture}
\node[matrix of math nodes] (M) {
\displaystyle\coprod_{j=1}^r F(V_i) &[2cm] F(E) \\[2cm]
F(U) & F(B). \\
};
\path[->, auto, arrow]
    (M-1-1) edge node{$\coprod_{j=1}^r \iota_j$} (M-1-2)
            edge node{$p$} node[swap]{$\cong\text{on sheets}$} (M-2-1)
    (M-1-2) edge node{$p$} (M-2-2)
    (M-2-1) edge (M-2-2)
;
\end{tikzpicture}
\]
Now $p(A)\cap U = p(A\cap p^{-1}(U)) = \coprod_{j=1}^r p(A\cap V_j)$ as sets. Each $A\cap V_i$ is open in $F(V_j)$ as $A\cap V_j=\iota_j^{-1}(A)$ and $A$ is assumed to be open in $F(E)$.
By the homeomorphism $F(V_i)\xto{p} F(U)$, we get that $p(A\cap V_j)$ is open in $F(U)$.
Hence $p(A)\cap U$ is open in $F(U)$ and thus in $F(E)$ since $F(U)$ is an open subspace of $F(B)$.
\end{proof}

\begin{cor}
Let $p\colon E\to B$ be a finite sheeted covering map between compactly generated weak Hausdorff spaces, then the induced map $\freeO n(p) \colon \freeO n(E)\to \freeO n(B)$ in the category of compactly generated weak Hausdorff spaces is a finite sheeted covering map for all $n\geq 0$. 
\end{cor}

\begin{proof}
Let $\freeO n(p)_0\colon \freeO n(E)_0\to \freeO n(B)_0$ be the induced map between the free loop spaces equipped with the compact-open topologies, then $\freeO n(p)_0$ is a covering map by Corollary \ref{corFreeLoopSpaceCoveringMap}.

Next, let $K\colon \Top \to \Top$ be the functor taking any topological space to its compactly-generated refinement. If we apply $K$ to either free loop space, $\freeO n(E)_0$ or $\freeO n(B)_0$, the space we get is not only compactly generated, but also weak Hausdorff by \cite{Strickland}*{Proposition 2.24}. Hence $\freeO n(E) = K\freeO n(E)_0$ and $\freeO n(B) = K\freeO n(B)_0$.

The functor $K\colon \Top \to \Top$ satisfies the first two properties of Lemma \ref{lemmaInducedCoveringMaps} by definition. The third property is \cite{Strickland}*{Lemma 2.26}. The fourth property is \cite{Strickland}*{Proposition 2.2}. It follows by Lemma \ref{lemmaInducedCoveringMaps} that $\freeO n(p) = K\freeO n(p)_0$ is a finite sheeted covering map.
\end{proof}

\begin{prop}\label{propKContinuousGeneralProof}
Given a finite sheeted covering map $\pi\colon Y\to X$, the map\linebreak $k\colon \PB n(\pi)\to (\Z_{>0})^n$ given in Definition \ref{defCovKMap} is continuous when $\freeO n (-)$ is equipped with the compact-open topology.

When $X$ and $Y$ are compactly generated weak Hausdorff, the conclusion also holds when $\freeO n (-)$ is equipped with the compactly generated refinement of the compact-open topology.
\end{prop}

\begin{proof}
For a compactly generated weak Hausdorff space $X$ the map $K\freeO n(-)_0 \to \freeO n(-)_0$ to the compact-open topology from its compactly generated refinement is continuous. Hence it is sufficient to prove that the map $k$ is continuous on the compact-open topology for $\freeO n (-)$, and any topological spaces.

Let $\pi\colon Y\to X$ be a finite sheeted covering map, and consider any point $(y,\tup s,f)\in \PB n(\pi)$. Then $f\colon (S^1)^n\to X$ is an $n$-fold loop, $\tup s\in (S^1)^n$ a point in the torus, and $y\in \pi^{-1}(f(\tup s))$ a point in the fiber over $f(\tup s)\in X$.
As in the proof of Proposition \ref{propMappingSpaceCoveringMap}, we apply Lemma \ref{lemmaEvenlyCoveredFiniteCovering} to find a finite intersection connected, open covering of the torus $(S^1)^n$ consisting of finitely many open sets $U_1,\dotsc,U_m$ such that each $f(\overline {U_i})$ is contained in an evenly covered open set $V_i\subseteq X$.
Following the notation of Proposition \ref{propMappingSpaceCoveringMap}, we choose trivializations
\[\pi^{-1}(V_i)\cong W_{i,1}\sqcup \dotsb \sqcup W_{i,r},\] 
where the number of sheets $r$ might depend on the evenly covered set $V_i$.

We again let $U_{ij}$ denote $U_i\cap U_j$, and for each $V_{ij}=V_i\cap V_j\neq \emptyset$, we partition $V_{ij}$ into open subsets $V_{ij}^\sigma $, where $\sigma\in \Sigma_r$, such that $V_{ij}^\sigma $ is evenly covered by the sheets $W_{i,1}\cap W_{j,\sigma(1)},\dotsc, W_{i,r}\cap W_{j,\sigma(r)}$.

As in the proof of Proposition \ref{propMappingSpaceCoveringMap}, we next consider the open neighborhood $f\in D\subseteq \freeO n(X)$, where
\[
D = \bigcap_{i=1}^n C_X(\overline{U_i}, V_i) \cap \bigcap_{\substack{1 \leq i < j\leq n\\ U_{ij}\neq \emptyset}} C_X(\overline{U_{ij}}, V_{ij}^{\sigma_{ij}}).
\]
Since the $U_i$ cover the torus, the tuple $\tup s$ is contained in some $U_{i_0}$. It follows that $f(\tup s)\in f(U_{i_0})\subseteq V_{i_0}$, and so the point $y\in \pi^{-1}(f(\tup s))$ lies in a unique sheet $W_{i_0,h_0}$. We claim that $k\colon \PB n(\pi)\to (\Z_{>0})^n$ is constant when restricted to 
\[
D' := \PB n(\pi)\cap ( W_{i_0,h_0}\times U_{i_0} \times D) \subseteq Y\times (S^1)^n\times \freeO n(X),
\]
which is an open neighborhood of $(y,\tup s,f)$ in $\PB n(\pi)$.

The union $V_{\cup} = \bigcup_{i=1}^m V_i$ is an open subset of $X$, and so $\freeO n(V_{\cup})$ is an open subset of $\freeO n(X)$ as well. We can pull back the covering map $\pi$ to $V_\cup$ and get the diagram of pullback squares, with $W_\cup = \bigcup_{i,h} W_{i,h}$:
\[
\begin{tikzpicture}
\node (M) [matrix of math nodes, nodes={text height=8pt, text depth=.7pt}] {
    \PB n(\pi|_{V_\cup}) &[1.5cm]\displaystyle W_\cup  &[1.5cm] Y\\[1.5cm]
    (S^1)^n\x \freeO n(V_\cup) & V_\cup & X.\\
};
\path [auto, ->,arrow]
    (M-1-1) edge (M-1-2)
            edge (M-2-1)
    (M-2-1) edge (M-2-2)
    (M-1-2) edge (M-1-3)
            edge node{$\pi|_{V_\cup}$} (M-2-2)
    (M-2-2) edge (M-2-3)
    (M-1-3) edge node {$\pi$} (M-2-3)
;
\end{tikzpicture}
\]
Here $\PB n(\pi|_{V_\cup})$ is an open subset of $\PB n(\pi)$ and we have an inclusion $(y,\tup s,f)\in D'\subseteq \PB n(\pi|_{V_\cup})$.

We will now replace $V_\cup$ with a finite topological space $V$. The finite space $V$ has a point $v_i$ corresponding to each non-empty $V_i$, a point $v_{ij}^\sigma$ corresponding to each non-empty $V_{ij}^\sigma$, a point $v_{ij\ell}^{\sigma_1\sigma_2}$ corresponding to each nonempty intersection
\[V_{ij\ell}^{\sigma_1\sigma_2} = V_{ij}^{\sigma_1}\cap V_{j\ell}^{\sigma_2} = V_{ij}^{\sigma_1}\cap V_{j\ell}^{\sigma_2} \cap V_{i\ell}^{\sigma_2\, \circ\, \sigma_1}, \]
and so on for further non-empty intersections $V_{ij\ell\dotsm}^{\sigma_1\sigma_2\dotsm}$.
The points of $V$ form a poset corresponding to the inclusion of intersections, and the topology of $V$ is the opposite Alexandrov topology in which downward closed sets are open. This topology has a basis given by the downward intervals 
\[
\{v \mid v\leq v_{ij\ell\dotsm}^{\sigma_1\sigma_2\dotsm}\}
\]
corresponding to all further intersections contained in $V_{ij\ell\dotsm}^{\sigma_1\sigma_2\dotsm}$.

We define a map $\rho\colon V_\cup \to V$ by sending each point $x\in V$ to the smallest $v_{ij\ell\dotsm}^{\sigma_1\sigma_2\dotsm}$ such that $V_{ij\ell\dotsm}^{\sigma_1\sigma_2\dotsm}$ contains $x$ (i.e. the intersection of all $V_{ij}^\sigma$'s containing $x$), or send $x$ to $v_i$ if $V_i$ is the only subset containing $x$. The map $\rho$ is then continuous since the preimage of the interval $\{v \mid v\leq v_{ij\ell\dotsm}^{\sigma_1\sigma_2\dotsm}\}$ is simply the finite intersection $V_{ij\ell\dotsm}^{\sigma_1\sigma_2\dotsm}$ as an open set in $X$.

The finite space $V$ has a covering space $W$ constructed in a similar way. The points of $W$, denoted by $w_{ij\ell\dotsm,h}^{\sigma_1\sigma_2\dotsm}$ correspond to all the finite non-empty intersection of sheets
\[
W_{ij\ell\dotsm,h}^{\sigma_1\sigma_2\dotsm} = W_{i,h}\cap W_{j,\sigma_1(h)} \cap W_{\ell, \sigma_2(\sigma_1(h))}\cap \dotsb.
\]
And $W$ is again a poset by inclusion of intersections, and we equip it with the opposite Alexandrov topology.

Each non-empty intersection  $V_{ij\ell\dotsm}^{\sigma_1\sigma_2\dotsm}$ is evenly covered by the sheets 
\[
W_{ij\ell\dotsm,1}^{\sigma_1\sigma_2\dotsm}, W_{ij\ell\dotsm,2}^{\sigma_1\sigma_2\dotsm}, \dotsc, W_{ij\ell\dotsm,r}^{\sigma_1\sigma_2\dotsm}.
\]
We can therefore define a finite sheeted covering map $\pi'\colon W\to V$, where the fiber over each point $v_{ij\ell\dotsm}^{\sigma_1\sigma_2\dotsm}$ is given by 
\[
(\pi')^{-1}(v_{ij\ell\dotsm}^{\sigma_1\sigma_2\dotsm}) = \{w_{ij\ell\dotsm,1}^{\sigma_1\sigma_2\dotsm}, w_{ij\ell\dotsm,2}^{\sigma_1\sigma_2\dotsm}, \dotsc, w_{ij\ell\dotsm,r}^{\sigma_1\sigma_2\dotsm}  \}.
\]
The map $\pi'$ is a covering map because each downward interval in the basis for the topology of $V$ is evenly covered by the corresponding intervals in $W$.

The covering map $\pi|_{V_\cup}$ over $V_\cup$ now arises in a new way as the pullback of the covering map $\pi'$ over $V$:
\[
\begin{tikzpicture}
\node (M) [matrix of math nodes, nodes={text height=8pt, text depth=.7pt}] {
    \PB n(\pi|_{V_\cup}) &[1.5cm]W_\cup &[1.5cm] W\\[1.5cm]
    (S^1)^n\x \freeO n(V_\cup) & V_\cup & V.\\
};
\path [auto, ->,arrow]
    (M-1-1) edge (M-1-2)
            edge (M-2-1)
    (M-2-1) edge (M-2-2)
    (M-1-2) edge (M-1-3)
            edge node{$\pi|_{V_\cup}$} (M-2-2)
    (M-2-2) edge node{$\rho$} (M-2-3)
    (M-1-3) edge node {$\pi'$} (M-2-3)
;
\end{tikzpicture}
\]
The preimage of each downward interval $\{v \mid v\leq v_{ij\ell\dotsm}^{\sigma_1\sigma_2\dotsm}\}$ in $V$ is the evenly covered open set $V_{ij\ell\dotsm}^{\sigma_1\sigma_2\dotsm}$ in $V_\cup$, and the sheets over $V_{ij\ell\dotsm}^{\sigma_1\sigma_2\dotsm}$ are precisely the preimage in $W_\cup$ of the sheets over $\{v \mid v\leq v_{ij\ell\dotsm}^{\sigma_1\sigma_2\dotsm}\}$.

The continuous map $\rho\colon V_\cup\to V$ induces a continuous map $\rho_*\colon \freeO n(V_\cup)\to \freeO n(V)$ and a further map $\rho_*\colon \PB n(\pi|_{V_\cup})\to \PB n(\pi')$.
Given a loop $\gamma\colon S^1\to V_\cup$ and some $s'\in S^1$ the action of $\gamma$ on the fiber over $\gamma(s)$ is isomorphic to the action of the composite $\rho\circ \gamma\colon S^1\to V$ on the fiber over $\rho(\gamma(s))$. Since the map $k\colon \PB n(\pi|_{V_\cup})\to (\Z_{>0})^n$ is defined in terms of such actions, we conclude that $k$ factors through the map $\rho_*\colon \PB n(\pi|_{V_\cup})\to \PB n(\pi')$. 

The triple $(y,\tup s, f)\in \PB n(\pi|_{V_\cup})$ is mapped by $\rho_*$ to the triple $(w_{ij\ell\dotsm,h}^{\sigma_1\sigma_2\dotsm}, \tup s, \rho\circ f)\in \PB n(\pi')$, where $y$ lies in the intersection of sheets $W_{ij\ell\dotsm,h}^{\sigma_1\sigma_2\dotsm}$.
Furthermore $\rho_*(D')$ is a subset of the neighborhood
\[
D'' := \PB n(\pi')\cap ( \{ w\leq w_{i_0,h_0}\}\times U_{i_0} \times \rho_*(D))) \subseteq W\times (S^1)^n\times \freeO n(V).
\]
It suffices to show that $k\colon \PB n(\pi')\to (\Z_{>0})^n$ is constant on $D''$.

We claim that if we consider the $n$-fold loops $f'\in D\subseteq \freeO n(V_\cup)$ and postcompose with $\rho$, then the resulting $n$-fold loops $\rho\circ f'\in \rho_*(D)\subseteq \freeO n(V)$ are all homotopic to each other in $\freeO n(V)$. The consequence of this is that while $D''$ might not be path connected, it is still contained in a single path-component inside $\PB n(\pi')$. We have used that $U_{i_0}$ is a connected open set in $(S^1)^n$ and thus path connected. Because $D''$ is contained in a single path-component of $\PB n(\pi')$, we can apply the proof of Proposition \ref{propKContinuousWithPathsProof} to see that $k$ takes the same value on all points of $D''$.

It suffices to prove that all $n$-fold loops in $\rho_*(D)$ are homotopic inside $\freeO n(V)$.
Consider the particular map $f_0\colon (S^1)^n\to V$ given by
$f_0(\tup s') = v_{ij\ell\dotsm}^{\sigma_1\sigma_2\dotsm}$, when $\tup s'\in U_{ij\ell\dotsm}$, where $U_{ij\ell\dotsm}$ is the intersection of all $U_i$'s containing $\tup s'$, and where 
\[
V_{ij\ell\dotsm}^{\sigma_1\sigma_2\dotsm} = V_{ij}^{\sigma_1} \cap V_{j\ell}^{\sigma_2} \cap \dotsb
\]
is the intersection of the neighborhoods corresponding to $U_{ij}, U_{j\ell},\dotsc$ in the definition of $D$.
Firstly, $f_0$ is continuous since the preimage of the downward interval $\{v \mid v\leq v_{ij\ell\dotsm}^{\sigma_1\sigma_2\dotsm}\}$ is $U_{ij\ell\dots}$ or empty.
Secondly, while $f_0$ need not itself be an element of $\rho_*(D)$, any $n$-fold loop $f'\in D$ satisfies
for every $\tup s'\in (S^1)^n$ that if $\tup s'$ lies in some intersection $U_{ij\ell\dotsm}$, then $f'(\tup s')$ lies in the corresponding intersection $V_{ij\ell\dotsm}^{\sigma_1\sigma_2\dotsm}$ (and possibly further $V$'s not already in this intersection). Consequently we have $\rho(f'(\tup s'))\leq f_0(\tup s')$ as elements in the poset $V$ for all $\tup s'\in (S^1)^n$. It follows that $\rho\circ f'$ is homotopic to $f_0$, e.g. with the homotopy $H\colon (S^1)^n\x I \to V$ given by
\[
H(\tup s',t)= \begin{cases}
\rho(f'(\tup s')) & \text{if $0\leq t < 1$,}
\\ f_0(\tup s') & \text{if $t=1$.}
\end{cases}
\]
To see that $H$ is continuous, note that the preimage of a downward interval $\{v \mid v\leq v_{ij\ell\dotsm}^{\sigma_1\sigma_2\dotsm}\}$ under $H$ is
\[
H^{-1}(\{v \mid v\leq v_{ij\ell\dotsm}^{\sigma_1\sigma_2\dotsm}\}) = \begin{cases} (f')^{-1}(V_{ij\ell\dotsm}^{\sigma_1\sigma_2\dotsm}) \x [0,1) \cup U_{ij\ell\dotsm} \x [0,1] & \\ &\hspace{-4cm} \text{if $V_{ij}^{\sigma_1}, V_{j\ell}^{\sigma_2}, \dotsc $ correspond }\\ & \hspace{-4cm} \text{to $U_{ij}, U_{j\ell}, \dotsc$in the definition of $D$,}
\\ (f')^{-1}(V_{ij\ell\dotsm}^{\sigma_1\sigma_2\dotsm}) \x [0,1) &\hspace{-2cm} \text{otherwise.}
\end{cases}
\]
Either of these is an open set in $(S^1)^n\x I$. We conclude that $H$ is a continuous homotopy from $f'$ to $f_0$ in $\freeO n(V)$. It follows that $\rho_*(D)$ is contained in a single path-component of $\freeO n(V)$, hence $D''$ is contained in a single path-component of $\PB n(\pi')$, hence $k\colon \PB n(\pi') \to (\Z_{>0})^n$ is constant on $D''$, hence $k\colon \PB n(\pi|_{V_\cup}) \to (\Z_{>0})^n$ is constant on $D'$, and finally it follows that $k\colon \PB n(\pi) \to (\Z_{>0})^n$ is continuous.
\end{proof}
